\documentclass[a4paper,10pt,reqno]{amsart}

\usepackage{hyperref}
\usepackage{amssymb}
\usepackage{dsfont}

\newcommand{\inpr}[2]{\ensuremath{\langle #1 , #2 \rangle}}

\def\grint{\mathrm{G}}
\newcommand{\lic}[2]{\ensuremath{\grint^{#1}(#2)}}
\def\grintfalc{\mathcal{G}}

\def\R{\mathbb{R}}
\def\Q{\mathbb{Q}}
\def\N{\mathbb{N}}
\def\Z{\mathbb{Z}}
\def\S{\mathbb{S}}
\def\H{\mathbb{H}}

\def\Acal{\mathcal{A}}
\def\Bcal{\mathcal{B}}
\def\Ccal{\mathcal{C}}
\def\Ncal{\mathcal{N}}
\def\Ical{\mathcal{I}}
\def\Pcal{\mathcal{P}}

\def\ind{\mathds{1}}

\def\hau{\mathcal{H}}
\def\leb{\mathcal{L}}

\def\sfrak{\mathfrak{s}}
\def\Hfrak{\mathfrak{H}}
\def\Rfrak{\mathfrak{R}}
\def\Ffrak{\mathfrak{F}}
\def\Efrak{\mathfrak{E}}
\def\Bfrak{\mathfrak{B}}
\def\Lfrak{\mathfrak{L}}
\def\dfrak{\mathfrak{d}}

\def\netm{\mathcal{M}}
\newcommand{\diam}[1]{\ensuremath{\left| #1 \right|}}

\def\eps{\varepsilon}
\def\ph{\varphi}

\def\gauge{\mathfrak{D}}

\def\dist{\mathrm{d}}
\def\dd{\mathrm{d}}
\def\ee{\mathrm{e}}
\def\proj{\mathrm{p}}

\def\as{\mbox{a.s.}}
\def\esp{\mathbb{E}}
\def\prob{\mathbb{P}}
\def\event{\mathcal{E}}
\def\eqlaw{\stackrel{d}{=}}

\def\Rho{\mathrm{P}}
\def\Mu{\mathrm{M}}
\def\Nu{\mathrm{N}}
\def\Alpha{\mathrm{A}}
\def\Tau{\mathrm{T}}
\def\Kappa{\mathrm{K}}
\def\Zeta{\mathrm{Z}}

\def\Qrm{\mathrm{Q}}
\def\Hyp{\mathrm{Hyp}}
\def\Hom{\mathrm{Hom}}

\def\Hdim{\dim_{\rm H}}
\def\Pdim{\dim_{\rm P}}

\newcommand{\opball}[2]{\mathrm{B}_{#1}(#2)}
\newcommand{\clball}[2]{\overline{\mathrm{B}}_{#1}(#2)}

\theoremstyle{plain}
\newtheorem{thm}{Theorem}
\newtheorem{prp}{Proposition}
\newtheorem{lem}{Lemma}
\newtheorem{cor}{Corollary}

\theoremstyle{definition}
\newtheorem{df}{Definition}

\theoremstyle{remark}

\title{Multifractal analysis of L\'evy fields}
\author{Arnaud Durand \and St\'ephane Jaffard}
\date{\today}
\address{Universit\'e Paris-Sud 11\\ Laboratoire de Math\'ematiques, UMR8628\\ Orsay, F-91405}
\email{arnaud.durand@math.u-psud.fr}
\address{Universit\'e Paris-Est Cr\'eteil Val de Marne\\ Laboratoire d'Analyse et de Math\'ematiques Appliqu\'ees, UMR8050\\ 61 avenue du G\'en\'eral de Gaulle\\ Cr\'eteil, F-94010}
\email{jaffard@univ-paris12.fr}

\subjclass[2010]{Primary: 60G60, 60G51; Secondary: 60G17, 60D05, 28A78, 28A80}

\keywords{L\'evy random fields, multifractal analysis, Hausdorff measures and dimension, sets with large intersection, Diophantine approximation, ubiquity}

\begin{document}

\begin{abstract}
We study the pointwise regularity properties of the L\'evy fields introduced by T.~Mori; these fields are the most natural generalization  of L\'evy processes to the multivariate setting. We determine their spectrum of singularities, and we show that their H\"older singularity sets satisfy a large intersection property in the sense of K.~Falconer.
\end{abstract}

\maketitle

%%%%%%%%%%%%%%%%%%%%%%%%%%%%%%%%%%%%%%%%%%%%%%%%%%%%%%%%%%%%%%%%%%%%%%%%%%%%%%
%%%%%%%%%%%%%%%%%%%%%%%%%%%%%%%%%%%%%%%%%%%%%%%%%%%%%%%%%%%%%%%%%%%%%%%%%%%%%%
%%%%%%%%%%%%%%%%%%%%%%%%%%%%%%%%%%%%%%%%%%%%%%%%%%%%%%%%%%%%%%%%%%%%%%%%%%%%%%
\section{Introduction}\label{sec:intro}
%%%%%%%%%%%%%%%%%%%%%%%%%%%%%%%%%%%%%%%%%%%%%%%%%%%%%%%%%%%%%%%%%%%%%%%%%%%%%%
%%%%%%%%%%%%%%%%%%%%%%%%%%%%%%%%%%%%%%%%%%%%%%%%%%%%%%%%%%%%%%%%%%%%%%%%%%%%%%
%%%%%%%%%%%%%%%%%%%%%%%%%%%%%%%%%%%%%%%%%%%%%%%%%%%%%%%%%%%%%%%%%%%%%%%%%%%%%%

The determination of the uniform and pointwise regularity of stochastic processes has been a longstanding issue, starting with the discovery of the uniform modulus of continuity of Brownian motion by L\'evy and Wiener, the law of the iterated logarithm by Khintchine and Kolmogorov, and the study of the irregularity of the Brownian paths by Paley, Wiener and Zygmund, as well as Dvoretzky.

Remarkably, a connexion between stochastic processes and fractals was also first established in the case of Brownian motion: {\em slow points} where the modulus of continuity is smaller than almost everywhere (the $\sqrt{\log\log(1/| t-t_0 |)}$ term is replaced by a constant $a  >0$), and {\em fast points} (where this term is replaced by $a\sqrt{\log(1/| t-t_0 |)}$)  were shown  to occur on random fractal sets;  the dimensions of these collections of sets (indexed by the parameter $a$)   were determined by E.~Perkins and by S.~Orey and S.J.~Taylor, respectively. This connexion  did not remain confined to the specific case of Brownian motion: in the mid 80s, {\em multifractal analysis} was introduced in order to propose a general framework for the study of the local variations of the regularity of stochastic processes, using the mathematical tools supplied by fractal analysis. Let us be more specific.

\begin{df}\label{df:point}
Let $f:\R^{d}\to\R$ be a locally bounded function, $t_{0}\in\R^d$ and $\alpha>0$. The function $f$ belongs to $C^{\alpha}(t_0)$ if there exist $C>0$ and a polynomial $P_{t_{0}}$ of degree less than $\alpha$ such that, for all $t$ in a neighborhood of $t_0$,
\[
|f(t) -P_{t_0}(t)|\leq C \|t-t_{0}||^{\alpha}.
\] 
The {\em H\"older exponent} of $f$ at $t_0$ is 
\[
\alpha_{f}(t_0)=\sup\{ \alpha\geq 0 \:|\: f \in C^{\alpha}(t_0)\}\in [0,\infty].
\]  
\end{df}

Multifractal analysis is concerned with the determination of the Hausdorff dimension of the sets of points where the H\"older exponent takes a given value. Two collections of sets play a particular role: the {\em iso-H\"older sets} of $f$, defined by
\[
E_{f}(h)=\{ t\in\R^{d} \:|\: \alpha_{f}(t)=h \},
\]
and the  {\em singularity sets}   $E'_{f}(h)$ which consist of the points $t$ where $f$ is continuous and  satisfies $\alpha_{f}(t)\leq h$. Note that, though these sets can be studied for any random field, they are really pertinent only when the H\"older exponent is nonconstant (which excludes for instance Brownian motion, where it takes everywhere the value $1/2$). If the H\"older exponent changes form point to point, one is interested in determining the Hausdorff dimension of the iso-H\"older sets $E_{f}(h)$. The corresponding notion is supplied by the {\em local spectrum of singularities} defined as follows.

\begin{df}
Let $f:\R^d\to\R$ be a locally bounded function and let $W$ be a nonempty open subset of $\R^{d}$. The local spectrum of singularities of $f$ is
\[
d_{f}(h,W)=\Hdim (E_{f}(h)\cap W),
\]
where $\Hdim$ denotes Hausdorff dimension (with the convention that $\Hdim\emptyset=-\infty$).
\end{df}

For many examples of random fields $Y$, the spectrum of singularities does not depend on the particular region $W$ that is considered, and is actually a deterministic function. If  such is the case, i.e. if there is a deterministic function $d_{Y}(h)$ such that
\[
\as \quad \forall h\in[0,\infty] \quad \forall W\neq\emptyset \mbox{ open} \qquad d_{Y}(h,W)=d_{Y}(h),
\]
the random field $Y$ is called {\em homogeneous} (following the terminology of~\cite{Jaffard:2004fh}). In that case, the support of the spectrum is the set of $h\in[0,\infty]$ such that $d_{Y}( h)\geq 0$ or, equivalently, such that $\alpha_{f}(t)=h$ for some $t\in\R^{d}$. A spectrum is {\em degenerate} if its support is reduced to a single point. This means that the corresponding H\"older exponent occurs everywhere (as is the case for Brownian motion). In that situation, the field in called a {\em mono-H\"older field}. Otherwise, it is called {\em multifractal}.

As shown in~\cite{Durand:2007fk,Jaffard:1999fg}, L\'evy processes supply examples of multifractal processes  which are homogeneous. Note that this is not the case for all classical random processes, see e.g.~\cite{Barral:2010vn} where J.~Barral, N.~Fournier, S.~Jaffard and S.~Seuret show examples of Markov processes whose spectrum is random and depends on the region, or~\cite{Durand:2007kx} where A.~Durand studies a model of wavelet series based on a tree-indexed Markov chain whose spectrum is random too. Moreover, though many  processes with a fairly general spectrum have been constructed (such as the Gaussian processes built by A.~Ayache, Jaffard and M.~Taqqu~\cite{Ayache:2007fk}), the only natural large class of homogeneous random processes which has been studied is that of L\'evy processes.

The introduction of multifractal analysis was motivated by  classification and model selection issues in signal processing (the analysis of one-dimensional turbulence data), and its fast development was boosted  by its pertinence for an unexpectedly large number of applications. Recent developments in image processing have stimulated similar needs in 2D, and classification tools based on multifractal parameters have recently proven to be promising directions of research~\cite{Wendt:2009fk}. These results motivate the investigation of random fields  that could  be used in image modeling, and the study of their multifractal properties (and in particular for porous media, medical imaging, etc.). Since L\'evy processes were the first class of processes proven to be multifractal, and as they play an important role in both theoretical probability and modeling, it is very natural to wonder if their multivariate extensions display similar properties.

We shall prove that this is the case for the fields that we study in this paper, which are multivariate analogs of L\'evy processes that jump along random hyperplanes. Due to that very anisotropic situation, new ubiquity techniques will have to be employed in order to perform their multifractal analysis. This strongly differs from other situations where geometry does not play any specific role, and the results in the one-dimensional case may straightforwardly be extended to higher dimensions. Such is the case for several multifractal models of random wavelet series~\cite{Arneodo:1998vn,Aubry:2002up,Durand:2007fk}.

Two main extensions of L\'evy processes to the multivariate setting, which one could term as L\'evy fields, have been proposed. The first one is that of R.~Adler, D.~Monrad, R.~Scissors and R.~Wilson~\cite{Adler:1983kx}, which finds its origins in the work of M.L.~Straf~\cite{Straf:1972ve}, and whose regularity was studied by M.E.~Vares~\cite{Vares:1981kx} and S.~Lagaize~\cite{Lagaize:2001vn} in the two-dimensional case. The second extension is due to T.~Mori~\cite{Mori:1992qf}, and covers the following important particular cases. In the Gaussian case, the L\'evy Brownian motion~\cite{Levy:1965uq} is the first instance of L\'evy field (in the sense of Mori) which has been considered. This field is isotropic, but its geometric construction, due to Chentsov~\cite{Chentsov:1957uq}, may easily be extended to the anisotropic setting. Subsequently, Chentsov's construction motivated the introduction of a class of isotropic and stable L\'evy fields, which also satisfy a remarkable selfsimilarity property. These fields, usually termed as the L\'evy-Chentsov fields, are discussed in G.~Samorodnitsky and Taqqu's book~\cite{Samorodnitsky:1994fk}, and were studied by N.-R.~Shieh~\cite{Shieh:1996ul}, who established their local boundedness and the existence of local times. In order to compare the relevance of these models, we propose two natural criteria for selecting a multivariate extension of L\' evy processes:

\begin{itemize}
\item {\bf Stability  under linear transforms of coordinates:} If $M$ is an invertible (deterministic) linear transform, and $Y(t)$ is a L\'evy field, then $Y(Mt)$ should also be a L\'evy field. In particular, the coordinate axes do not play any specific role in the construction or the analysis of these fields.
\item {\bf Stability under trace:} The restriction of a $d$-dimensional L\'evy field to any $d'$-dimensional linear subspace is a $d'$-dimensional L\'evy field. Note that the first requirement implies that we do not need to specify a particular system of coordinates on that subspace. Furthermore, restrictions to arbitrary half lines are L\'evy processes.
\end{itemize}

Among the aforementioned extensions, Mori's is the only one that fulfills these two criteria. It turns out that Mori's definition of a L\'evy field is also a natural multivariate analog of  the definition of a L\'evy process. Indeed, recall that a stochastic process indexed by $[0,\infty)$ is a L\'evy process if it satisfies the following properties: it has stationary and independent increments, is stochastically continuous, and vanishes at zero almost surely, see e.g.~\cite{Bertoin:1996uq,Sato:1999fk}. This implies in particular that the one-dimensional marginals are infinitely divisible.

\begin{df}
A random field $Y=\{Y(t),\,t\in\R^{d}\}$ is a {\em L\'evy field} (in the sense of Mori) if the following conditions are satisfied:
\begin{enumerate}
\item it is stochastically continuous and vanishes at zero almost surely;
\item it has stationary increments, i.e.~$Y(a+\cdot\,)-Y(a)\eqlaw Y$ for any $a\in\R^{d}$;
\item its finite-dimensional marginals are infinitely divisible;
\item for any $a,b\in\R^{d}$, the increments of $\{Y(a+\lambda b),\,\lambda\in\R\}$ are independent.
\end{enumerate}
\end{df}

In the previous definition, $\eqlaw$ stands for equality of the finite dimensional distributions. It is easy to check that for any $a,b\in\R^{d}$, the process $\{Y(a+\lambda b)-Y(a),\,\lambda\geq 0\}$ has stationary and independent increments, is stochastically continuous, and vanishes at zero almost surely, thereby being a L\'evy process. More generally, it is clear that the two stability requirements listed above are satisfied.

%%%%%%%%%%%%%%%%%%%%%%%%%%%%%%%%%%%%%%%%%%%%%%%%%%%%%%%%%%%%%%%%%%%%%%%%%%%%%%
%%%%%%%%%%%%%%%%%%%%%%%%%%%%%%%%%%%%%%%%%%%%%%%%%%%%%%%%%%%%%%%%%%%%%%%%%%%%%%
%%%%%%%%%%%%%%%%%%%%%%%%%%%%%%%%%%%%%%%%%%%%%%%%%%%%%%%%%%%%%%%%%%%%%%%%%%%%%%
\section{Main results}
%%%%%%%%%%%%%%%%%%%%%%%%%%%%%%%%%%%%%%%%%%%%%%%%%%%%%%%%%%%%%%%%%%%%%%%%%%%%%%
%%%%%%%%%%%%%%%%%%%%%%%%%%%%%%%%%%%%%%%%%%%%%%%%%%%%%%%%%%%%%%%%%%%%%%%%%%%%%%
%%%%%%%%%%%%%%%%%%%%%%%%%%%%%%%%%%%%%%%%%%%%%%%%%%%%%%%%%%%%%%%%%%%%%%%%%%%%%%

%%%%%%%%%%%%%%%%%%%%%%%%%%%%%%%%%%%%%%%%%%%%%%%%%%%%%%%%%%%%%%%%%%%%%%%%%%%%%%
%%%%%%%%%%%%%%%%%%%%%%%%%%%%%%%%%%%%%%%%%%%%%%%%%%%%%%%%%%%%%%%%%%%%%%%%%%%%%%
\subsection{Representation of L\'evy fields}\label{subsec:represLevy}
%%%%%%%%%%%%%%%%%%%%%%%%%%%%%%%%%%%%%%%%%%%%%%%%%%%%%%%%%%%%%%%%%%%%%%%%%%%%%%
%%%%%%%%%%%%%%%%%%%%%%%%%%%%%%%%%%%%%%%%%%%%%%%%%%%%%%%%%%%%%%%%%%%%%%%%%%%%%%

Mori established a remarkable decomposition of L\'evy fields into three parts which is similar to the L\'evy-It\=o decomposition of a L\'evy process into a linear drift, a Brownian motion, and a jump component controlled by the L\' evy measure, see e.g.~\cite[Chapter~4]{Sato:1999fk}. In the $d$-dimensional case, the three components are the following.
\begin{enumerate}
\item A {\bf linear drift}, that is, a mapping of the form $t\mapsto\inpr{a}{t}$, where $a$ is a deterministic vector and $\inpr{\cdot}{\cdot}$ denotes the standard inner product in $\R^{d}$.
\item A {\bf Gaussian field} $B_{\mu}=\{B_{\mu}(t),\,t\in\R^{d}\}$ that depends on a finite nonnegative symmetric Borel measure $\mu$ defined on the unit sphere $\S^{d-1}$ of $\R^{d}$. Here, symmetric means invariant under the mapping $s\mapsto -s$. The construction of $B_{\mu}$ is detailed in Section~\ref{sec:defGaussianLevy}. Note that $B_{\mu}$ is almost surely constant equal to zero when $\mu=0$.
\item A {\bf jump field} $L_{\nu}=\{L_{\nu}(t),\,t\in\R^{d}\}$. Recall that the jump component of a L\'evy process is a sum of independent compensated compound Poisson processes; the structure is the same in the multivariate setting, except that Poisson processes, which jump at points, are replaced by random fields which jump along Poisson distributed hyperplanes. Their distribution is described by a $d$-dimensional analog of the L\'evy measure, namely, a nonnegative Borel measure $\nu$ defined on $\S^{d-1}\times\R^{*}$, where $\R^{*}$ means $\R\setminus\{0\}$. The measure  $\nu$ is symmetric, i.e.~invariant under $(s,x)\mapsto(-s,-x)$, and
\begin{equation}\label{eq:nuLevy}
\int_{s\in\S^{d-1}\atop x\in\R^{*}} (1\wedge x^{2})\,\nu(\dd s,\dd x)<\infty,
\end{equation}
where $\wedge$ stands for minimum. Intuitively, $\nu(\dd s,\dd x)$ describes the amount of hyperplanes orthogonal to $s$ where a jump of size $x$ occurs. In particular, the jump field $L_{\nu}$ is almost surely zero everywhere when $\nu=0$. The construction of $L_{\nu}$ is detailed in Section~\ref{sec:defnonGaussianLevy}.
\end{enumerate}

\begin{thm}[Mori]\label{thm:LevyIto}
Every L\'evy field $Y=\{Y(t),\,t\in\R^{d}\}$ may be represented in the following manner:
\[
Y\eqlaw\inpr{a}{\cdot}+B_{\mu}+L_{\nu},
\]
for some $a\in\R^{d}$, some finite symmetric measure $\mu$ on $\S^{d-1}$, and some symmetric measure $\nu$ on $\S^{d-1}\times\R^{*}$ satisfying~(\ref{eq:nuLevy}). Moreover, $B_{\mu}$ and $L_{\nu}$ are independent and the triple $(a,\mu,\nu)$ is uniquely determined by the field $Y$.
\end{thm}

Due to Theorem~\ref{thm:LevyIto}, studying the regularity of a L\'evy field reduces to analyzing each of these three components, and then to understanding what happens when combining them. To begin with, note that linear drifts are $C^{\infty}$ and thus play no role in our analysis.

%%%%%%%%%%%%%%%%%%%%%%%%%%%%%%%%%%%%%%%%%%%%%%%%%%%%%%%%%%%%%%%%%%%%%%%%%%%%%%
%%%%%%%%%%%%%%%%%%%%%%%%%%%%%%%%%%%%%%%%%%%%%%%%%%%%%%%%%%%%%%%%%%%%%%%%%%%%%%
\subsection{Regularity of the Gaussian component}
%%%%%%%%%%%%%%%%%%%%%%%%%%%%%%%%%%%%%%%%%%%%%%%%%%%%%%%%%%%%%%%%%%%%%%%%%%%%%%
%%%%%%%%%%%%%%%%%%%%%%%%%%%%%%%%%%%%%%%%%%%%%%%%%%%%%%%%%%%%%%%%%%%%%%%%%%%%%%

Our main result on the regularity of the Gaussian field $B_{\mu}$ is the next statement.

\begin{thm}\label{thm:HolderBmu}
Let $\mu$ be a finite symmetric measure on $\S^{d-1}$ such that $\mu\neq 0$. Then, $B_{\mu}$ is a homogeneous mono-H\"older field and
\[
\as \quad \forall t\in\R^{d} \qquad \alpha_{B_{\mu}}(t)=1/2.
\]
\end{thm}

That theorem is in fact a consequence of two slightly more precise results on the pointwise modulus of continuity of $B_{\mu}$, namely, Propositions~\ref{prp:modcontBmu} and~\ref{prp:irregBmu} below. Their proofs are given in Section~\ref{sec:regGaussian} and rely on standard techniques for studying the oscillations of Gaussian random fields.

%%%%%%%%%%%%%%%%%%%%%%%%%%%%%%%%%%%%%%%%%%%%%%%%%%%%%%%%%%%%%%%%%%%%%%%%%%%%%%
%%%%%%%%%%%%%%%%%%%%%%%%%%%%%%%%%%%%%%%%%%%%%%%%%%%%%%%%%%%%%%%%%%%%%%%%%%%%%%
\subsection{Regularity of the jump component}
%%%%%%%%%%%%%%%%%%%%%%%%%%%%%%%%%%%%%%%%%%%%%%%%%%%%%%%%%%%%%%%%%%%%%%%%%%%%%%
%%%%%%%%%%%%%%%%%%%%%%%%%%%%%%%%%%%%%%%%%%%%%%%%%%%%%%%%%%%%%%%%%%%%%%%%%%%%%%

Very precise results on the size of the iso-H\"older sets of $L_{\nu}$ are detailed in Section~\ref{sec:regnonGaussian}. In order to give a first insight into them, let us single out a representative consequence: Theorem~\ref{thm:specLnu}, which yields the spectrum of singularities of $L_{\nu}$. Its statement involves an index $\beta_{\nu}$ associated with $\nu$ and holds under an admissibility condition bearing on that measure.

\begin{df}\label{df:index}
Let $\nu$ be a symmetric measure on $\S^{d-1}\times\R^{*}$ satisfying~(\ref{eq:nuLevy}). 
The {\em index} of $\nu$ is
\begin{equation}\label{eq:defbetanu}
\beta_{\nu}=\inf\Biggl\{ \gamma\geq 0 \:\Biggl|\: \int_{s\in\S^{d-1}\atop x\in(0,1]} x^{\gamma} \,\nu(\dd s,\dd x)<\infty \Biggr\}.
\end{equation}
For any integer $j\geq 1$, let $\nu_{j}=\nu(\S^{d-1}\times(2^{-j},2^{-j+1}])$; $\nu$ is {\em admissible} if
\begin{equation}\label{eq:defchinu}
\chi_{\nu}=\sum_{j=1}^{\infty} 2^{-j}(j\nu_{j})^{1/2}<\infty.
\end{equation}
\end{df}

The index $\beta_{\nu}$ is the $d$-dimensional analog of the index that Blumenthal and Getoor associated to a L\'evy process through~\cite[Definition~2.1]{Blumenthal:1961kx}; note that~(\ref{eq:nuLevy}) implies that $\beta_{\nu}\in [0,2]$. We shall almost always assume below that the measure $\nu$ is admissible. A similar assumption had been made in~\cite{Jaffard:1999fg}. This condition is slightly stronger than~(\ref{eq:nuLevy}), which ensures the existence of $L_{\nu}$ and amounts to the finiteness of $\sum 2^{-2j}\nu_{j}$. Actually, assuming~(\ref{eq:defchinu}) is a mild restriction since, for instance, every measure $\nu$ with index less than two is admissible.

\begin{thm}\label{thm:specLnu}
Let $\nu$ be an admissible measure with $\beta_{\nu}>0$. Then, $L_{\nu}$ is a homogeneous multifractal field and with probability one,
\[
\forall h\in[0,\infty] \quad \forall W\neq\emptyset \mbox{ open} \qquad
d_{L_{\nu}}(h,W)=
\begin{cases}
d-1+\beta_{\nu}h & \mbox{if } h\leq 1/\beta_{\nu}\\
-\infty & \mbox{if } h>1/\beta_{\nu}.
\end{cases}
\]
\end{thm}

This theorem is a straightforward consequence of Proposition~\ref{prp:sizeEnuhsimple} and Corollary~\ref{cor:sizepropEnuh} below. As well as the results of Section~\ref{sec:regnonGaussian}, it covers the isotropic and stable case of the L\'evy-Chentsov random fields~\cite{Samorodnitsky:1994fk}, for which $\nu(\dd s,\dd x)=\sigma(\dd s)\dd x/|x|^{\alpha+1}$ where $\alpha\in (0,2)$ and $\sigma$ denotes the uniform measure on $\S^{d-1}$. In that situation, $\beta_{\nu}=\alpha$.

Contrary to what the hypothesis of Theorem~\ref{thm:specLnu} may suggest, the analysis developed in Section~\ref{sec:regnonGaussian} also includes the case in which $\beta_{\nu}=0$. In particular, we shall discuss the regularity of compound Poisson fields (multivariate analogs of compound Poisson processes), for which $\nu$ has finite total mass. In that situation, the field is {\em not} homogeneous. Indeed, locally, the field jumps on at most finitely many hyperplanes, thereby being $C^{\infty}$ except on a set of dimension $d-1$ where its H\"older exponent vanishes, see Proposition~\ref{prp:Enuhnufinite} for a precise statement.

%%%%%%%%%%%%%%%%%%%%%%%%%%%%%%%%%%%%%%%%%%%%%%%%%%%%%%%%%%%%%%%%%%%%%%%%%%%%%%
%%%%%%%%%%%%%%%%%%%%%%%%%%%%%%%%%%%%%%%%%%%%%%%%%%%%%%%%%%%%%%%%%%%%%%%%%%%%%%
\subsection{Spectrum of singularities of a general L\'evy field}
%%%%%%%%%%%%%%%%%%%%%%%%%%%%%%%%%%%%%%%%%%%%%%%%%%%%%%%%%%%%%%%%%%%%%%%%%%%%%%
%%%%%%%%%%%%%%%%%%%%%%%%%%%%%%%%%%%%%%%%%%%%%%%%%%%%%%%%%%%%%%%%%%%%%%%%%%%%%%

We call a L\'evy field {\em canonical} if it is of the form $Y_{a,\mu,\nu}=\inpr{a}{\cdot}+B_{\mu}+L_{\nu}$ for some $a\in\R^{d}$, some finite symmetric measure $\mu$ on $\S^{d-1}$, and some symmetric measure $\nu$ on $\S^{d-1}\times\R^{*}$ satisfying~(\ref{eq:nuLevy}). By virtue of Theorem~\ref{thm:LevyIto}, every L\'evy field has the same finite-dimensional distributions as a canonical one. In view of Definition~\ref{df:index}, we call $\beta_{\nu}$ the index of $Y_{a,\mu,\nu}$, and we call that field admissible if $\nu$ is admissible. Moreover, $Y_{a,\mu,\nu}$ is said to have a Gaussian component if and only if $\mu\neq 0$.

We show in Section~\ref{sec:Gaussianplusjumps} that the H\"older exponent of $Y_{a,\mu,\nu}$ is everywhere the minimum of that of its components $B_{\mu}$ and $L_{\nu}$. In view of Theorems~\ref{thm:HolderBmu} and~\ref{thm:specLnu}, this leads immediately to the next statement. This approach also easily enables one to deduce comparable results for the cases that are treated in Section~\ref{sec:regnonGaussian} but not covered by Theorem~\ref{thm:specLnu}. For instance, one could easily infer the spectrum of singularities of a canonical L\'evy field whose jump component is a compound Poisson field.

\begin{cor}\label{cor:specgeneral}
Let $Y$ be an admissible canonical L\'evy field with index $\beta>0$. Then, $Y$ is homogeneous. Moreover,
\begin{itemize}
\item if $Y$ has no Gaussian component, then with probability one,
\[
\forall h\in[0,\infty] \quad \forall W\neq\emptyset \mbox{ open} \qquad
d_{Y}(h,W)=
\begin{cases}
d-1+\beta h & \mbox{if } h\leq 1/\beta\\
-\infty & \mbox{if } h>1/\beta;
\end{cases}
\]
\item if $Y$ has a Gaussian component, then with probability one,
\[
\forall h\in[0,\infty] \quad \forall W\neq\emptyset \mbox{ open} \qquad
d_{Y}(h,W)=
\begin{cases}
d-1+\beta h & \mbox{if } h<1/2\\
d & \mbox{if } h=1/2\\
-\infty & \mbox{if } h>1/2.
\end{cases}
\]
\end{itemize}
\end{cor}

As expected, when $d=1$, the previous result boils down to~\cite[Theorem~1]{Jaffard:1999fg} which gives the spectrum of singularities of a L\'evy process.

%%%%%%%%%%%%%%%%%%%%%%%%%%%%%%%%%%%%%%%%%%%%%%%%%%%%%%%%%%%%%%%%%%%%%%%%%%%%%%
%%%%%%%%%%%%%%%%%%%%%%%%%%%%%%%%%%%%%%%%%%%%%%%%%%%%%%%%%%%%%%%%%%%%%%%%%%%%%%
\subsection{Large intersection properties and ubiquity}
%%%%%%%%%%%%%%%%%%%%%%%%%%%%%%%%%%%%%%%%%%%%%%%%%%%%%%%%%%%%%%%%%%%%%%%%%%%%%%
%%%%%%%%%%%%%%%%%%%%%%%%%%%%%%%%%%%%%%%%%%%%%%%%%%%%%%%%%%%%%%%%%%%%%%%%%%%%%%

We will not restrict our study to the determination of the Hausdorff dimension of the random sets of points related with  the definition of H\"older regularity, but we will also investigate some of their geometric properties. It turns out that, rather than the iso-H\"older sets, the singularity sets mentioned in Section~\ref{sec:intro} are those which display the most striking features. We shall show in Section~\ref{subsec:largeintpropEnuh} that the singularity sets of the jump component $L_{\nu}$ satisfy a remarkable counterintuitive property introduced by K.~Falconer~\cite{Falconer:1994hx}: they fall in the category of {\em sets with large intersection}. Recall that the intersection of two subsets of $\R^{d}$ with dimension $d_{1}$ and $d_{2}$ respectively is usually expected to be equal to $d_{1}+d_{2}-d$ (the codimensions add up), as in the case of affine subspaces, see~\cite[Chapter~8]{Falconer:2003oj} for precise statements. Sets with large intersection disprove this heuristic in a striking way: their size properties are not altered by taking countable intersection. As a matter of fact, the Hausdorff dimension of the intersection of countably many sets with large intersection is equal to the infimum of their Hausdorff dimensions.  Sets with large intersection have been shown to play a prominent role in metric number theory (Diophantine approximation) and dynamical systems, see e.g.~\cite{Durand:2007uq,Durand:2008jk,Falconer:1994hx} and the references therein. In the context of probability, they are relevant to the multifractal analysis of random wavelet series~\cite{Aubry:2002up,Durand:2007kx}, as well as the study of random coverings of the circle~\cite{Durand:2008kx}. Durand~\cite{Durand:2007fk} also proved that the singularity sets of L\'evy processes are sets with large intersection. In all these situations, large intersection properties arise because the sets under study are derived from an underlying {\em ubiquitous system}. This is also the case here, and our proofs make use of new extensions of ubiquity that we will develop (see Theorem~\ref{thm:sliKnualpha} below and its proof). Besides the aforementioned papers, we also refer to~\cite{Barral:2004ae,Barral:2006fk,Beresnevich:2006ve} for important results on ubiquity and its applications.

%%%%%%%%%%%%%%%%%%%%%%%%%%%%%%%%%%%%%%%%%%%%%%%%%%%%%%%%%%%%%%%%%%%%%%%%%%%%%%
%%%%%%%%%%%%%%%%%%%%%%%%%%%%%%%%%%%%%%%%%%%%%%%%%%%%%%%%%%%%%%%%%%%%%%%%%%%%%%
\subsection{Behavior of traces}
%%%%%%%%%%%%%%%%%%%%%%%%%%%%%%%%%%%%%%%%%%%%%%%%%%%%%%%%%%%%%%%%%%%%%%%%%%%%%%
%%%%%%%%%%%%%%%%%%%%%%%%%%%%%%%%%%%%%%%%%%%%%%%%%%%%%%%%%%%%%%%%%%%%%%%%%%%%%%

An important and difficult subject of investigation is to understand how the multifractal properties of a field and its traces on linear subspaces are related, see e.g.~\cite{Aubry:2010fk}. This question initially came up in the context of the analysis of turbulence: the only high precision experimental data available are one-dimensional cuts and the challenge is to infer from these cuts information about the multifractal properties of the whole field. In general, it is expected that the spectrum of the trace is the initial spectrum lowered by the codimension of the subspace (see~(\ref{eq:shiftspec}) below) and the parts which become negative are set to $-\infty$. L\'evy fields provide a case study of that effect, because the characteristic parameters of their traces are easy to obtain.

To be more specific, let $e=(e_{1},\ldots,e_{d'})$ be an arbitrary orthonormal system of $\R^{d}$ with $1\leq d'\leq d$, and let $Y_{a,\nu,\mu}$ denote a canonical L\'evy field. Then, the random field $Y^{e}_{a,\nu,\mu}$ defined by
\[
\forall t_{1},\ldots,t_{d'}\in\R^{d'} \qquad Y^{e}_{a,\nu,\mu}(t_{1},\ldots,t_{d'})=Y_{a,\nu,\mu}(t_{1}e_{1}+\ldots+t_{d'}e_{d'})
\]
is a canonical L\'evy field indexed by $\R^{d'}$ whose characteristic triple $(a_{e},\mu_{e},\nu_{e})$ may be deduced from $(a,\mu,\nu)$ with the help of the mapping $\proj_{e}:\R^{d}\to\R^{d'}$ defined by $\proj_{e}(t)=(\inpr{t}{e_{1}},\ldots,\inpr{t}{e_{d'}})$ for any $t\in\R^{d}$. To be specific, $a_{e}=\proj_{e}(a)$, the measure $\mu_{e}$ is the image under $s\mapsto\proj_{e}(s)/\|\proj_{e}(s)\|$ of $\|\proj_{e}(s)\|\,\mu(\dd s)$, and the measure $\nu_{e}$ is the image under $(s,x)\mapsto(\proj_{e}(s)/\|\proj_{e}(s)\|,x)$ of $\|\proj_{e}(s)\|\,\nu(\dd s,\dd x)$. In particular, the field $Y^{e}_{a,\nu,\mu}$ has no Gaussian component if and only if $\mu$ is supported in the orthogonal complement of the linear span of $e$. Moreover, $\nu_{e}$ is admissible whenever $\nu$ is, and the index of $\nu_{e}$ is at most that of $\nu$.

In the isotropic case, for which $\nu$ is the product of the uniform measure on $\S^{d-1}$ and a given measure on $\R^{*}$, the index of $\nu_{e}$ coincides with that of $\nu$, regardless of the choice of $e$. It follows from Corollary~\ref{cor:specgeneral} that, as expected,
\begin{equation}\label{eq:shiftspec}
d_{Y^{e}_{a,\nu,\mu}}=d_{Y_{a,\nu,\mu}}-(d-d').
\end{equation}
However, for appropriate anisotropic choices of the measure $\nu$, one may obtain a whole range of values for the index $\beta_{\nu_{e}}$. Our results then lead to a whole variety of spectra for $Y^{e}_{a,\nu,\mu}$, depending on the choice of the directions for the trace.

%%%%%%%%%%%%%%%%%%%%%%%%%%%%%%%%%%%%%%%%%%%%%%%%%%%%%%%%%%%%%%%%%%%%%%%%%%%%%%
%%%%%%%%%%%%%%%%%%%%%%%%%%%%%%%%%%%%%%%%%%%%%%%%%%%%%%%%%%%%%%%%%%%%%%%%%%%%%%
\subsection{Directional regularity}
%%%%%%%%%%%%%%%%%%%%%%%%%%%%%%%%%%%%%%%%%%%%%%%%%%%%%%%%%%%%%%%%%%%%%%%%%%%%%%
%%%%%%%%%%%%%%%%%%%%%%%%%%%%%%%%%%%%%%%%%%%%%%%%%%%%%%%%%%%%%%%%%%%%%%%%%%%%%%

Note that the notion of pointwise regularity given in Definition~\ref{df:point} does not take into account directional regularity but yields the worst possible regularity in all directions. Therefore, the results obtained in this paper do not take into account possible directional irregularity phenomena. However, such phenomena are to be expected in the case of L\'evy fields. Indeed, given that they display jumps along hyperplanes, they have, by construction, a very anisotropic nature. Therefore, it would be of great interest to perform a multifractal analysis of these fields using a more flexible notion of pointwise smoothness, which can take directionality into account (see e.g.~\cite{Ben-Slimane:1998uq,Jaffard:2010ys} for appropriate definitions).

%%%%%%%%%%%%%%%%%%%%%%%%%%%%%%%%%%%%%%%%%%%%%%%%%%%%%%%%%%%%%%%%%%%%%%%%%%%%%%
%%%%%%%%%%%%%%%%%%%%%%%%%%%%%%%%%%%%%%%%%%%%%%%%%%%%%%%%%%%%%%%%%%%%%%%%%%%%%%
\subsection{Roadmap}
%%%%%%%%%%%%%%%%%%%%%%%%%%%%%%%%%%%%%%%%%%%%%%%%%%%%%%%%%%%%%%%%%%%%%%%%%%%%%%
%%%%%%%%%%%%%%%%%%%%%%%%%%%%%%%%%%%%%%%%%%%%%%%%%%%%%%%%%%%%%%%%%%%%%%%%%%%%%%

The paper is organized as follows. In Sections~\ref{sec:defGaussianLevy} and~\ref{sec:defnonGaussianLevy}, we detail the construction of the two main components appearing in Mori's decompositon of L\'evy fields: the Gaussian part and the jump part. We will also derive some basic properties which will be useful for their multifractal analysis.

Section~\ref{sec:regGaussian} is devoted to the proof of Theorem~\ref{thm:HolderBmu}, according to which the Gaussian part has everywhere the H\"older exponent $1/2$. Precise results on the size (in terms of Hausdorff measures, Hausdorff dimension and packing dimension) and large intersection properties of the iso-H\"older and the singularity sets of the jump part are stated in Section~\ref{sec:regnonGaussian}. These results lead to Theorem~\ref{thm:specLnu} above. In Section~\ref{sec:Gaussianplusjumps}, we explain how the results on the Gaussian and the jump components may be combined to obtain the pointwise regularity of a general canonical L\'evy field.

The rest of the paper is devoted to establishing the results of Section~\ref{sec:regnonGaussian}. The structure of the proof is described in Section~\ref{sec:architecture}; we present there the required tools, specifically, a precise knowledge of the location of the singularities of the jump part, and a description of the size and large intersection properties of the set of points approximated at a certain rate by random hyperplanes that are distributed in a Poissonian way. The first tool is detailed in Section~\ref{sec:proofprpcharacLnuKnualpha}. The second one is presented in Sections~\ref{sec:proofprpKnualphaRd} and~\ref{sec:proofthmsliKnualpha}, and relies heavily on ubiquity. Last, Section~\ref{sec:proofsecregnonGaussian} details the proofs of the results of Section~\ref{sec:regnonGaussian}, and the paper ends with the proof of a lemma which is called upon by the first ingredient, see Section~\ref{sec:prooflemoscLnuj}.

%%%%%%%%%%%%%%%%%%%%%%%%%%%%%%%%%%%%%%%%%%%%%%%%%%%%%%%%%%%%%%%%%%%%%%%%%%%%%%
%%%%%%%%%%%%%%%%%%%%%%%%%%%%%%%%%%%%%%%%%%%%%%%%%%%%%%%%%%%%%%%%%%%%%%%%%%%%%%
%%%%%%%%%%%%%%%%%%%%%%%%%%%%%%%%%%%%%%%%%%%%%%%%%%%%%%%%%%%%%%%%%%%%%%%%%%%%%%
\section{The Gaussian component}\label{sec:defGaussianLevy} 
%%%%%%%%%%%%%%%%%%%%%%%%%%%%%%%%%%%%%%%%%%%%%%%%%%%%%%%%%%%%%%%%%%%%%%%%%%%%%%
%%%%%%%%%%%%%%%%%%%%%%%%%%%%%%%%%%%%%%%%%%%%%%%%%%%%%%%%%%%%%%%%%%%%%%%%%%%%%%
%%%%%%%%%%%%%%%%%%%%%%%%%%%%%%%%%%%%%%%%%%%%%%%%%%%%%%%%%%%%%%%%%%%%%%%%%%%%%%

The Gaussian component of a L\'evy field is a Gaussian random field $B_{\mu}=\{B_{\mu}(t),\,t\in\R^{d}\}$ depending on a finite nonnegative symmetric Borel measure $\mu$ defined on $\S^{d-1}$. We now recall the construction of such fields $B_{\mu}$, which is essentially due to Chentsov~\cite{Chentsov:1957uq}. Note that, in the isotropic case (where $\mu$ is the uniform measure on the sphere), one basically ends up with a geometric representation of the L\'evy Brownian motion~\cite{Levy:1965uq}.

%%%%%%%%%%%%%%%%%%%%%%%%%%%%%%%%%%%%%%%%%%%%%%%%%%%%%%%%%%%%%%%%%%%%%%%%%%%%%%
%%%%%%%%%%%%%%%%%%%%%%%%%%%%%%%%%%%%%%%%%%%%%%%%%%%%%%%%%%%%%%%%%%%%%%%%%%%%%%
\subsection{Definition of $B_{\mu}$}
%%%%%%%%%%%%%%%%%%%%%%%%%%%%%%%%%%%%%%%%%%%%%%%%%%%%%%%%%%%%%%%%%%%%%%%%%%%%%%
%%%%%%%%%%%%%%%%%%%%%%%%%%%%%%%%%%%%%%%%%%%%%%%%%%%%%%%%%%%%%%%%%%%%%%%%%%%%%%

Let us consider the collection $\Bcal_{0}(\H_{d})$ of all relatively compact Borel subsets of the product space $\H_{d}=(0,\infty)\times\S^{d-1}$, and the centered Gaussian process $\Bfrak_{\mu}=\{\Bfrak_{\mu}(V),\,V\in\Bcal_{0}(\H_{d})\}$ with covariance function given by
\[
\forall V,V'\in\Bcal_{0}(\H_{d}) \qquad \esp\left[\Bfrak_{\mu}(V)\Bfrak_{\mu}(V')\right]=\int_{(\rho,s)\in\H_{d}}\ind_{\{(\rho,s)\in V\cap V'\}}\,\dd\rho\,\mu(\dd s).
\]
Such a process $\Bfrak_{\mu}$ is often referred to as a white noise, and may roughly be regarded as a random signed measure on $\H_{d}$, although strictly speaking it is not. The reason is that for any disjoint sets $V,V'\in\Bcal_{0}(\H_{d})$, the random variables $\Bfrak_{\mu}(V\cup V')$ and $\Bfrak_{\mu}(V)+\Bfrak_{\mu}(V')$ coincide almost surely.

At this point, it is useful to mention that $\Bcal_{0}(\H_{d})$ is a $\delta$-ring, i.e.~is closed under symmetric difference and countable intersections. Moreover, there is a one-to-one correspondence between the set $\H_{d}$ and the collection of all $(d-1)$-dimensional hyperplanes in $\R^{d}$ that do not contain the origin. Indeed, any such hyperplane $h$ may be represented in a unique manner by a pair $(\rho,s)\in\H_{d}$ since it coincides with the set of all $t\in\R^{d}$ satisfying $\rho=\inpr{s}{t}$. In that correspondence, the hyperplanes that separate a given point $t\in\R^{d}$ and the origin are those which are represented by a pair that belongs to the set $V_{t}\in\Bcal_{0}(\H_{d})$ given by
\begin{equation}\label{eq:defVt}
V_{t}=\{(\rho,s)\in\H_{d} \:|\: \rho<\inpr{s}{t}\}.
\end{equation}

The geometric intuition behind the definition of the random field $B_{\mu}$ is that its value at a given point $t$ is determined by the mass that $\Bfrak_{\mu}$ assigns to the hyperplanes separating $t$ and the origin. Specifically,
\[
\forall t\in\R^{d} \qquad B_{\mu}(t)=\Bfrak_{\mu}(V_{t}).
\]
More generally, the increment of the field $B_{\mu}$ between two points in $\R^{d}$ is determined by the mass of the hyperplanes that separate them.

%%%%%%%%%%%%%%%%%%%%%%%%%%%%%%%%%%%%%%%%%%%%%%%%%%%%%%%%%%%%%%%%%%%%%%%%%%%%%%
%%%%%%%%%%%%%%%%%%%%%%%%%%%%%%%%%%%%%%%%%%%%%%%%%%%%%%%%%%%%%%%%%%%%%%%%%%%%%%
\subsection{Basic properties}
%%%%%%%%%%%%%%%%%%%%%%%%%%%%%%%%%%%%%%%%%%%%%%%%%%%%%%%%%%%%%%%%%%%%%%%%%%%%%%
%%%%%%%%%%%%%%%%%%%%%%%%%%%%%%%%%%%%%%%%%%%%%%%%%%%%%%%%%%%%%%%%%%%%%%%%%%%%%%

It is quite straightforward to establish that $B_{\mu}$ falls in the category of L\'evy fields. In particular, checking the stationary increments property for $B_{\mu}$ calls upon the fact that for any $t,t'\in\R^{d}$, the random variable $B_{\mu}(t')-B_{\mu}(t)$ is normally distributed with mean zero and variance
\begin{equation}\label{eq:varBmu}
\dfrak_{\mu}(t,t')^{2}=\esp[(B_{\mu}(t')-B_{\mu}(t))^{2}]=\frac{1}{2}\int_{s\in\S^{d-1}}|\inpr{s}{t'-t}|\,\mu(\dd s),
\end{equation}
which depends on $t$ and $t'$ only through $t'-t$. We refer to~\cite[Equation~(5.11)]{Mori:1992qf} for details on how to derive the last equality in~(\ref{eq:varBmu}). In addition, it is customary to observe that $\dfrak_{\mu}$ defines a pseudometric on $\R^{d}$ which satisfies
\begin{equation}\label{eq:boundmetricBmu}
\forall t,t'\in\R^{d} \qquad \dfrak_{\mu}(t,t')\leq c_{\mu}\|t-t'\|^{1/2},
\end{equation}
with $c_{\mu}=(\mu(\S^{d-1})/2)^{1/2}$, in view of~(\ref{eq:varBmu}) and the Cauchy-Schwarz inequality. Here and below, $\|\cdot\|$ denotes the Euclidean norm. This observation implies that the field $B_{\mu}$ is stochastically continuous; and also admits a separable modification, see e.g.~\cite[Section~4]{Lifshits:1995kx}. For technical reasons that come into play in Section~\ref{sec:regGaussian} (see the proof of Proposition~\ref{prp:modcontBmu}), we need to work with such a modification. Therefore, we assume throughout that the field $B_{\mu}$ is separable.

%%%%%%%%%%%%%%%%%%%%%%%%%%%%%%%%%%%%%%%%%%%%%%%%%%%%%%%%%%%%%%%%%%%%%%%%%%%%%%
%%%%%%%%%%%%%%%%%%%%%%%%%%%%%%%%%%%%%%%%%%%%%%%%%%%%%%%%%%%%%%%%%%%%%%%%%%%%%%
%%%%%%%%%%%%%%%%%%%%%%%%%%%%%%%%%%%%%%%%%%%%%%%%%%%%%%%%%%%%%%%%%%%%%%%%%%%%%%
\section{The jump component}\label{sec:defnonGaussianLevy}
%%%%%%%%%%%%%%%%%%%%%%%%%%%%%%%%%%%%%%%%%%%%%%%%%%%%%%%%%%%%%%%%%%%%%%%%%%%%%%
%%%%%%%%%%%%%%%%%%%%%%%%%%%%%%%%%%%%%%%%%%%%%%%%%%%%%%%%%%%%%%%%%%%%%%%%%%%%%%
%%%%%%%%%%%%%%%%%%%%%%%%%%%%%%%%%%%%%%%%%%%%%%%%%%%%%%%%%%%%%%%%%%%%%%%%%%%%%%

As seen in Section~\ref{subsec:represLevy}, the jump component of a L\'evy field is a random field $L_{\nu}=\{L_{\nu}(t),\,t\in\R^{d}\}$ that depends on a nonnegative Borel measure $\nu$ on $\S^{d-1}\times\R^{*}$. Recall that $\nu$ is symmetric and satisfies~(\ref{eq:nuLevy}). Moreover, let $\leb^{1}_{+}$ be the Lebesgue measure on $(0,\infty)$ and let $\Nu$ denote a Poisson random measure on $(0,\infty)\times\S^{d-1}\times\R^{*}$ with intensity $\leb^{1}_{+}\otimes\nu$, see e.g.~\cite{Neveu:1977mz}.

%%%%%%%%%%%%%%%%%%%%%%%%%%%%%%%%%%%%%%%%%%%%%%%%%%%%%%%%%%%%%%%%%%%%%%%%%%%%%%
%%%%%%%%%%%%%%%%%%%%%%%%%%%%%%%%%%%%%%%%%%%%%%%%%%%%%%%%%%%%%%%%%%%%%%%%%%%%%%
\subsection{Multivariate compound Poisson processes}
%%%%%%%%%%%%%%%%%%%%%%%%%%%%%%%%%%%%%%%%%%%%%%%%%%%%%%%%%%%%%%%%%%%%%%%%%%%%%%
%%%%%%%%%%%%%%%%%%%%%%%%%%%%%%%%%%%%%%%%%%%%%%%%%%%%%%%%%%%%%%%%%%%%%%%%%%%%%%

For any set $V\in\Bcal_{0}(\H_{d})$, let
\[
\Lfrak_{\nu,0}(V)=\int_{(\rho,s)\in V\atop|x|>1} x\,\Nu(\dd\rho,\dd s,\dd x).
\]
What plays the role of a compound Poisson process with jumps of magnitude larger than one in our situation is the random field $L_{\nu,0}=\{L_{\nu,0}(t),\,t\in\R^{d}\}$ defined by $L_{\nu,0}(t)=\Lfrak_{\nu,0}(V_{t})$ for any $t\in\R^{d}$, where $V_{t}$ is given by~(\ref{eq:defVt}). In fact, letting $(\Rho_{n},S_{n},X_{n})$, for $n\geq 1$, denote the atoms of the Poisson measure $\Nu$, we have
\begin{equation}\label{eq:Lnu0}
L_{\nu,0}(t)=\sum_{n=1}^{\infty} X_{n}\,\ind_{\{\Rho_{n}<\inpr{S_{n}}{t},\,|X_{n}|>1\}}
\end{equation}
for all $t\in\R^{d}$. For any integer $A\geq 1$, when $t$ ranges over the closed ball with radius $A$ centered at the origin, the sum in~(\ref{eq:Lnu0}) may actually be restricted to the almost surely finitely many integers $n\geq 1$ satisfying both $\Rho_{n}<A$ and $|X_{n}|>1$. Therefore, $L_{\nu,0}$ is almost surely piecewise constant, with jumps of magnitude $|X_{n}|>1$ located on the hyperplanes $H_{n}$ parametrized by $(\Rho_{n},S_{n})\in\H_{d}$, which are given by
\begin{equation}\label{eq:defHn}
H_{n}=\{ t\in\R^{d} \:|\: \Rho_{n}=\inpr{S_{n}}{t} \}.
\end{equation}

%%%%%%%%%%%%%%%%%%%%%%%%%%%%%%%%%%%%%%%%%%%%%%%%%%%%%%%%%%%%%%%%%%%%%%%%%%%%%%
%%%%%%%%%%%%%%%%%%%%%%%%%%%%%%%%%%%%%%%%%%%%%%%%%%%%%%%%%%%%%%%%%%%%%%%%%%%%%%
\subsection{Multivariate compensated sums of jumps}
%%%%%%%%%%%%%%%%%%%%%%%%%%%%%%%%%%%%%%%%%%%%%%%%%%%%%%%%%%%%%%%%%%%%%%%%%%%%%%
%%%%%%%%%%%%%%%%%%%%%%%%%%%%%%%%%%%%%%%%%%%%%%%%%%%%%%%%%%%%%%%%%%%%%%%%%%%%%%

For any integer $j\geq 1$, the compensated sum corresponding to the jumps of magnitude in $\Ical_{j}=(2^{-j},2^{-j+1}]$ is the field $L_{\nu,j}=\{L_{\nu,j}(t),\,t\in\R^{d}\}$ given by $L_{\nu,j}(t)=\Lfrak_{\nu,j}(V_{t})$ for any $t\in\R^{d}$, where
\begin{equation}\label{eq:defLfraknuj}
\Lfrak_{\nu,j}(V)=\int_{(\rho,s)\in V\atop|x|\in\Ical_{j}} x\,\Nu(\dd\rho,\dd s,\dd x)-\int_{(\rho,s)\in V\atop|x|\in\Ical_{j}} x\,\dd\rho\,\nu(\dd s,\dd x)
\end{equation}
for any Borel set $V\in\Bcal_{0}(\H_{d})$. Note that $\Lfrak_{\nu,j}(V)$ may be regarded as an integral with respect to the compensated Poisson measure $\Nu^{*}=\Nu-\leb^{1}_{+}\otimes\nu$ associated with $\Nu$. In addition, due to the symmetry of $\nu$, we have
\begin{equation}\label{eq:Lnujalt}
L_{\nu,j}(t)=\sum_{n=1}^{\infty} X_{n}\ind_{\{\Rho_{n}<\inpr{S_{n}}{t},\,|X_{n}|\in\Ical_{j}\}}-\int_{s\in\S^{d-1}\atop x\in\Ical_{j}} x\inpr{s}{t}\,\nu(\dd s,\dd x)
\end{equation}
for all $t\in\R^{d}$. The sum in~(\ref{eq:Lnujalt}) is almost surely piecewise constant with jumps of magnitude $|X_{n}|\in\Ical_{j}$ located on the hyperplanes $H_{n}$ given by~(\ref{eq:defHn}), while the compensating integral depends linearly on $t$.

%%%%%%%%%%%%%%%%%%%%%%%%%%%%%%%%%%%%%%%%%%%%%%%%%%%%%%%%%%%%%%%%%%%%%%%%%%%%%%
%%%%%%%%%%%%%%%%%%%%%%%%%%%%%%%%%%%%%%%%%%%%%%%%%%%%%%%%%%%%%%%%%%%%%%%%%%%%%%
\subsection{Definition and basic properties of $L_{\nu}$}
%%%%%%%%%%%%%%%%%%%%%%%%%%%%%%%%%%%%%%%%%%%%%%%%%%%%%%%%%%%%%%%%%%%%%%%%%%%%%%
%%%%%%%%%%%%%%%%%%%%%%%%%%%%%%%%%%%%%%%%%%%%%%%%%%%%%%%%%%%%%%%%%%%%%%%%%%%%%%

The series formed by the compensated sums~(\ref{eq:defLfraknuj}) for $j\geq 1$ converges and yields a L\'evy field with jumps of magnitude at most one. As a matter of fact, for any $V\in\Bcal_{0}(\H_{d})$ and $j\geq 1$, Campbell's theorem~\cite{Kingman:1993gf} ensures that the random variable $\Lfrak_{\nu,j}(V)$ is centered with variance
\[
\int_{(\rho,s)\in V\atop|x|\in\Ical_{j}} x^{2}\,\dd\rho\,\nu(\dd s,\dd x).
\]
In view of~(\ref{eq:nuLevy}), it follows that the series $\sum_{j\geq 1}\Lfrak_{\nu,j}(V)$ converges in $L^{2}$. Thus, it converges almost surely in view of the L\'evy-It\=o-Nisio theorem for sums of independent random variables, see~\cite[p.~151]{Ledoux:1991qy}. This enables us to define
\[
\Lfrak_{\nu}(V)=\Lfrak_{\nu,0}(V)+\sum_{j=1}^{\infty}\Lfrak_{\nu,j}(V).
\]
Then, $\Lfrak_{\nu}=\{\Lfrak_{\nu}(V),\,V\in\Bcal_{0}(\H_{d})\}$ is an infinitely divisible random measure on $\H_{d}$ with control measure $\nu$, in the sense that it satisfies the next properties:
\begin{itemize}
\item it is independently scattered, i.e.~for any disjoint sets $V_{1},\ldots,V_{n}\in\Bcal_{0}(\H_{d})$, the random variables $\Lfrak_{\nu}(V_{1}),\ldots,\Lfrak_{\nu}(V_{n})$ are independent;
\item it is $\sigma$-additive, i.e.~for any sequence $(V_{n})_{n\geq 1}$ of disjoint sets in $\Bcal_{0}(\H_{d})$ whose union belongs to $\Bcal_{0}(\H_{d})$, the series $\sum_{n\geq 1}\Lfrak_{\nu}(V_{n})$ converges almost surely and its sum is equal to $\Lfrak_{\nu}(\bigcup_{n\geq 1} V_{n})$;
\item for every $V\in\Bcal_{0}(\H_{d})$, the characteristic function of $\Lfrak_{\nu}(V)$ is given by
\begin{equation}\label{eq:chfLfraknuV}
\esp[\ee^{i\theta \Lfrak_{\nu}(V)}]=\exp\int_{(\rho,s)\in V\atop x\in\R^{*}}\left(\ee^{i\theta x}-1-i\theta x\,\ind_{\{|x|\leq 1\}} \right)\,\dd\rho\,\nu(\dd s,\dd x).
\end{equation}
\end{itemize}
The first two properties directly follow from standard results on Poisson random measures, while the third one is a consequence of Campbell's theorem, see~\cite{Kingman:1993gf,Neveu:1977mz}.

Making use of these properties, it is straightforward to check that the random field $L_{\nu}=\{L_{\nu}(t),\,t\in\R^{d}\}$ defined by
\[
\forall t\in\R^{d} \qquad L_{\nu}(t)=\Lfrak_{\nu}(V_{t})
\]
is a L\'evy field, see~\cite{Mori:1992qf} for details. Replacing the measure $\nu$ by appropriate restrictions, this implies that the fields $L_{\nu,j}$ defined above are of L\'evy type as well.

%%%%%%%%%%%%%%%%%%%%%%%%%%%%%%%%%%%%%%%%%%%%%%%%%%%%%%%%%%%%%%%%%%%%%%%%%%%%%%
%%%%%%%%%%%%%%%%%%%%%%%%%%%%%%%%%%%%%%%%%%%%%%%%%%%%%%%%%%%%%%%%%%%%%%%%%%%%%%
\subsection{Comments}
%%%%%%%%%%%%%%%%%%%%%%%%%%%%%%%%%%%%%%%%%%%%%%%%%%%%%%%%%%%%%%%%%%%%%%%%%%%%%%
%%%%%%%%%%%%%%%%%%%%%%%%%%%%%%%%%%%%%%%%%%%%%%%%%%%%%%%%%%%%%%%%%%%%%%%%%%%%%%

Mori opted for the cut-off function $1/(1+x^{2})$ instead of $\ind_{\{|x|\leq 1\}}$ in~(\ref{eq:chfLfraknuV}). This leads to a slightly different expression of the field $\Lfrak_{\nu}$ and of the drift coefficient $a$. Yet, our choice clearly does not compromise the validity of Theorem~\ref{thm:LevyIto} and does alter the value of neither $\mu$ nor $\nu$, given a L\'evy field $Y$. Furthermore, Mori did not detail the construction of the infinitely divisible random measure $\Lfrak_{\nu}$ on which the field $L_{\nu}$ is based. However, the proofs below call upon a precise knowledge of the jump structure of $L_{\nu}$, and this explains why we chose above to present the construction of $\Lfrak_{\nu}$ in terms of compensated Poisson integrals.

Furthermore, in order to study the regularity of $L_{\nu}$, we first need to make sure that its H\"older  exponent is a well-defined quantity. This boils down to verifying that the sample functions of that field exist almost surely. In fact, as yet, $L_{\nu}$ has been defined in a pointwise manner only: we merely proved the almost sure convergence of the series $\sum_{j\geq 0}L_{\nu,j}(t)$ defining $L_{\nu}(t)$, for every fixed $t\in\R^{d}$. This enables us to consider the finite-dimensional marginals of the field $L_{\nu}$, which is sufficient to state Theorem~\ref{thm:LevyIto}. However, in order to determine the value of the H\"older exponent $\alpha_{L_{\nu}}(t)$ at a given point $t\in\R^{d}$, we need to consider $L_{\nu}$ {\em everywhere} near $t$. The next result indicates that this is possible.

\begin{prp}\label{prp:existpathsLnu}
Let $\nu$ be an admissible measure. Then, with probability one,
\[
\forall t\in\R^{d} \qquad L_{\nu}(t)=\sum_{j=0}^{\infty} L_{\nu,j}(t)\mbox{ exists.}
\]
\end{prp}

The proof of Proposition~\ref{prp:existpathsLnu} is postponed to Section~\ref{subsec:proofprpexistpathsLnu}. It relies on precise estimates of the increments of the fields $L_{\nu,j}$ that are given by Lemma~\ref{lem:oscLnuj} below. It can also easily be adapted to show that, when the measure $\nu$ is admissible, the sample paths of the field $L_{\nu}$ are almost surely locally bounded; this may be seen as an extension of~\cite[Theorem~2.3]{Shieh:1996ul}, which concerns the L\'evy-Chentsov fields only.

In what follows, whenever $\nu$ is supposed to be admissible, we implicitly work on the almost sure event where the sample paths of the field $L_{\nu}$ exist.

%%%%%%%%%%%%%%%%%%%%%%%%%%%%%%%%%%%%%%%%%%%%%%%%%%%%%%%%%%%%%%%%%%%%%%%%%%%%%%
%%%%%%%%%%%%%%%%%%%%%%%%%%%%%%%%%%%%%%%%%%%%%%%%%%%%%%%%%%%%%%%%%%%%%%%%%%%%%%
%%%%%%%%%%%%%%%%%%%%%%%%%%%%%%%%%%%%%%%%%%%%%%%%%%%%%%%%%%%%%%%%%%%%%%%%%%%%%%
\section{Regularity of the Gaussian component}\label{sec:regGaussian}
%%%%%%%%%%%%%%%%%%%%%%%%%%%%%%%%%%%%%%%%%%%%%%%%%%%%%%%%%%%%%%%%%%%%%%%%%%%%%%
%%%%%%%%%%%%%%%%%%%%%%%%%%%%%%%%%%%%%%%%%%%%%%%%%%%%%%%%%%%%%%%%%%%%%%%%%%%%%%
%%%%%%%%%%%%%%%%%%%%%%%%%%%%%%%%%%%%%%%%%%%%%%%%%%%%%%%%%%%%%%%%%%%%%%%%%%%%%%

This section is devoted to the proof of Theorem~\ref{thm:HolderBmu}; we shall in fact establish two slightly more precise results. First, the H\"older exponent of the random field $B_{\mu}$ is almost surely at least $1/2$ everywhere, as a direct consequence of the following result on its modulus of continuity. Recall that $c_{\mu}=(\mu(\S^{d-1})/2)^{1/2}$.

\begin{prp}\label{prp:modcontBmu}
There exists a universal constant $\Kappa>0$ such that for any finite symmetric measure $\mu$ on $\S^{d-1}$ satisfying $c_{\mu}>0$ and for any integer $A\geq 1$,
\begin{equation}\label{eq:modcontBmu}
\as \qquad \limsup_{\delta\to 0}\frac{1}{(\delta\log(1/\delta))^{1/2}}\sup_{t,t'\in[-A,A]^{d} \atop \|t-t'\|\leq\delta}|B_{\mu}(t')-B_{\mu}(t)|\leq \Kappa c_{\mu} d^{1/2}.
\end{equation}
\end{prp}

\begin{proof}
As mentioned in Section~\ref{sec:defGaussianLevy}, the field $B_{\mu}$ is assumed to be separable. We may therefore apply~\cite[Theorem~1.3.5]{Adler:2007rt}. Accordingly, there exists a universal constant $\Kappa>0$ such that with probability one, for $\eta>0$ small enough,
\begin{equation}\label{eq:entropyboundBmu}
\sup_{t,t'\in[-A,A]^{d} \atop \dfrak_{\mu}(t,t')\leq\eta}|B_{\mu}(t')-B_{\mu}(t)|\leq\Kappa\int_{0}^{\eta} (\log N([-A,A]^{d},\dfrak_{\mu},\eps))^{1/2} \,\dd\eps.
\end{equation}
Here, $N([-A,A]^{d},\dfrak_{\mu},\eps)$ denotes the minimal number of balls with radius $\eps$ that cover the cube $[-A,A]^{d}$, the balls being closed, centered in that cube and taken in the sense of the pseudometric $\dfrak_{\mu}$ defined by~(\ref{eq:varBmu}). Letting $\lfloor\cdot\rfloor$ stand for the floor function, it is easy to check that the aforementioned cube is covered by $(1+2\lfloor A(c_{\mu}/\eps)^{2} d^{1/2}\rfloor)^{d}$ closed Euclidean balls centered in it with radius $(\eps/c_{\mu})^{2}$, and~(\ref{eq:boundmetricBmu}) implies that each of these balls is included in a closed ball with radius $\eps$ for the pseudometric $\dfrak_{\mu}$. Hence, the right-hand side of~(\ref{eq:entropyboundBmu}) is bounded above by
\[
\Kappa\int_{0}^{\eta}\biggl(d\log\biggl(1+\frac{2Ac_{\mu}^{2}d^{1/2}}{\eps^{2}}\biggr)\biggr)^{1/2}\,\dd\eps \ \sim\ \Kappa\eta\left(2d\log\frac{c_{\mu}}{\eta}\right)^{1/2} \quad\mbox{as}\quad \eta\to 0.
\]
To conclude, it now suffices to let $\eta=c_{\mu}\sqrt{\delta}$ and observe that the supremum in left-hand side of~(\ref{eq:modcontBmu}) is bounded above by the left-hand side of~(\ref{eq:entropyboundBmu}), thanks to~(\ref{eq:boundmetricBmu}).
\end{proof}

Second, the fact that the H\"older exponent of $B_{\mu}$ is almost surely at most $1/2$ everywhere follows directly from the next proposition; $\clball{t}{\delta}$ will denote the closed Euclidean ball centered at $t$ with radius $\delta$.

\begin{prp}\label{prp:irregBmu}
For any finite symmetric measure $\mu$ on $\S^{d-1}$ satisfying $c_{\mu}>0$, there exists a real $\kappa_{d,\mu}>0$ such that
\[
\as \quad \forall t\in\R^{d} \quad \forall \delta>0 \quad \exists t'\in\clball{t}{\delta} \qquad \left|B_{\mu}(t')-B_{\mu}(t)\right|>\kappa_{d,\mu}\|t'-t\|^{1/2}.
\]
\end{prp}

\begin{proof}
We shall adapt some ideas that Dvoretzky~\cite{Dvoretzky:1963vn} employed in the case of Brownian motion. To begin with, let $(e_{1},\ldots,e_{d})$ denote the canonical basis of $\R^{d}$. We necessarily have $\dfrak_{\mu}(0,e_{i})>0$ for some $i\in\{1,\ldots,d\}$, because
\[
\sum_{i=1}^{d} \dfrak_{\mu}(0,e_{i})^{2}\geq\frac{1}{2}\sum_{i=1}^{d}\int_{s\in\S^{d-1}}|\inpr{s}{e_{i}}|^{2}\,\mu(\dd s)=\int_{s\in\S^{d-1}}\|s\|^{2}\,\mu(\dd s)=2c_{\mu}^{2}>0,
\]
due to~(\ref{eq:varBmu}). Then, let us consider an integer $A\geq 1$, a real number $\kappa>0$ and let us assume that for any real $\delta>0$, there exists a point $t\in[-A,A]^{d}$ such that
\begin{equation}\label{eq:holderBmu}
\forall t'\in\clball{t}{\delta} \qquad \left|B_{\mu}(t')-B_{\mu}(t)\right|\leq \kappa\|t'-t\|^{1/2}.
\end{equation}
Hence, for every integer $n\geq 1$, there exists a $k\in\{-An,\ldots,An-1\}^{d}$ such that~(\ref{eq:holderBmu}) holds for some $t\in\R^{d}$ with $t-k/n\in[0,1/n]^{d}$. Letting $\log_{2}$ stand for base two logarithm and $j_{n}=\lfloor\log_{2}(n\delta/d^{1/2})\rfloor-1$, and assuming that $n$ is large enough to ensure that $j_{n}\geq 1$, we now see that for any $j\in\{0,\ldots,j_{n}\}$, the point $(k+2^{j}e_{i})/n$ belongs to the ball $\clball{t}{\delta}$, so that
\[
\left|B_{\mu}\biggl(\frac{k+2^{j}e_{i}}{n}\biggr)-B_{\mu}(t)\right|\leq \kappa\left\|\frac{k+2^{j}e_{i}}{n}-t\right\|^{1/2}\leq \kappa\biggl(\frac{2^{j+1} d^{1/2}}{n}\biggr)^{1/2},
\]
thanks to~(\ref{eq:holderBmu}). Due to the triangle inequality, this implies that
\begin{equation}\label{eq:incdiscBmu}
\forall j\in\{1,\ldots,j_{n}\} \ \ \left|B_{\mu}\biggl(\frac{k+2^{j}e_{i}}{n}\biggr)-B_{\mu}\biggl(\frac{k+2^{j-1}e_{i}}{n}\biggr)\right|\leq \kappa\frac{d^{1/4}}{\sqrt{n}}2^{(j+3)/2}.
\end{equation}
The increments in the left-hand side, for $j\in\{1,\ldots,j_{n}\}$, are independent and normally distributed with mean zero and variance $2^{j-1}\dfrak_{\mu}(0,e_{i})^{2}/n$ by virtue of~(\ref{eq:varBmu}). Thus,~(\ref{eq:incdiscBmu}) occurs with probability at most
\[
\prod_{j=1}^{j_{n}}\prob\left(\biggl|\frac{2^{j/2}}{\sqrt{2n}}\dfrak_{\mu}(0,e_{i})\zeta\biggr|\leq \kappa\frac{d^{1/4}}{\sqrt{n}}2^{(j+3)/2}\right)=q(\kappa)^{j_{n}},
\]
where $q(\kappa)=\prob(|\zeta|\leq 4\kappa d^{1/4}/\dfrak_{\mu}(0,e_{i}))$ and $\zeta$ denotes a standard normal random variable. As a result, the probability that~(\ref{eq:holderBmu}) holds for some point $t\in[-A,A]^{d}$ is at most $(2An)^{d}q(\kappa)^{j_{n}}$ for every integer $n$ such that $j_{n}\geq 1$. Clearly, we have $q(\kappa_{d,\mu})<2^{-d}$ for some $\kappa_{d,\mu}>0$, which implies that $(2An)^{d}q(\kappa_{d,\mu})^{j_{n}}$ tends to zero as $n\to\infty$. We deduce that for any integer $A\geq 1$ and any real $\delta>0$,
\[
\as \quad \forall t\in[-A,A]^{d}\quad\exists t'\in\clball{t}{\delta} \qquad \left|B_{\mu}(t')-B_{\mu}(t)\right|>\kappa_{d,\mu}\|t'-t\|^{1/2}.
\]
The desired result clearly follows.
\end{proof}

%%%%%%%%%%%%%%%%%%%%%%%%%%%%%%%%%%%%%%%%%%%%%%%%%%%%%%%%%%%%%%%%%%%%%%%%%%%%%%
%%%%%%%%%%%%%%%%%%%%%%%%%%%%%%%%%%%%%%%%%%%%%%%%%%%%%%%%%%%%%%%%%%%%%%%%%%%%%%
%%%%%%%%%%%%%%%%%%%%%%%%%%%%%%%%%%%%%%%%%%%%%%%%%%%%%%%%%%%%%%%%%%%%%%%%%%%%%%
\section{Regularity of the jump component}\label{sec:regnonGaussian}
%%%%%%%%%%%%%%%%%%%%%%%%%%%%%%%%%%%%%%%%%%%%%%%%%%%%%%%%%%%%%%%%%%%%%%%%%%%%%%
%%%%%%%%%%%%%%%%%%%%%%%%%%%%%%%%%%%%%%%%%%%%%%%%%%%%%%%%%%%%%%%%%%%%%%%%%%%%%%
%%%%%%%%%%%%%%%%%%%%%%%%%%%%%%%%%%%%%%%%%%%%%%%%%%%%%%%%%%%%%%%%%%%%%%%%%%%%%%

Let us detail our results on the regularity of the random field $L_{\nu}$ defined in Section~\ref{sec:defnonGaussianLevy}. The corresponding proofs are given in Section~\ref{sec:proofsecregnonGaussian}, and their architecture is presented in Section~\ref{sec:architecture}. Throughout, we assume that $\nu$ is admissible.

In Section~\ref{subsec:sizepropEnuh}, we describe the size properties of the iso-H\"older sets $E_{L_{\nu}}(h)$ in a very precise manner: we give the value of their Hausdorff $g$-measure in every open subset of $\R^{d}$, for every gauge function $g$. In what follows, the iso-H\"older sets are denoted by $E_{\nu}(h)$ instead of $E_{L_{\nu}}(h)$, for the sake of brevity. Specifically,
\[
\forall h\in[0,\infty] \qquad E_{\nu}(h)=\{ t\in\R^{d} \:|\: \alpha_{L_{\nu}}(t)=h \}.
\]
We get similar results for the singularity sets raised in Section~\ref{sec:intro} and defined by
\begin{equation}\label{eq:defEnuhp}
E'_{\nu}(h)=\{ t\in\R^{d}\setminus J_{\nu} \:|\: \alpha_{L_{\nu}}(t)\leq h \}.
\end{equation}
We also establish that the sets $E'_{\nu}(h)$ fall in the category of sets with large intersection, see Section~\ref{subsec:largeintpropEnuh}. As explained in Section~\ref{subsec:packdimEnuh}, this property has straightforward consequences on the value of their packing dimension.

The set $J_{\nu}$ in~(\ref{eq:defEnuhp}) consists of the points at which $L_{\nu}$ jumps. Specifically, $J_{\nu}$ is the union over $n\geq 1$ of the hyperplanes $H_{n}$ defined by~(\ref{eq:defHn}), in terms of the atoms of the Poisson measure $\Nu$ that arises in the construction of $L_{\nu}$. Equivalently,
\begin{equation}\label{eq:defJnu}
J_{\nu}=\left\{ t\in\R^{d} \:\Biggl|\: \int_{(\rho,s)\in\H_{d}\atop x\in\R^{*}}\ind_{\{\rho=\inpr{s}{t}\}} \,\Nu(\dd\rho,\dd s,\dd x)\geq 1 \right\}.
\end{equation}

%%%%%%%%%%%%%%%%%%%%%%%%%%%%%%%%%%%%%%%%%%%%%%%%%%%%%%%%%%%%%%%%%%%%%%%%%%%%%%
%%%%%%%%%%%%%%%%%%%%%%%%%%%%%%%%%%%%%%%%%%%%%%%%%%%%%%%%%%%%%%%%%%%%%%%%%%%%%%
\subsection{Preliminary remark}
%%%%%%%%%%%%%%%%%%%%%%%%%%%%%%%%%%%%%%%%%%%%%%%%%%%%%%%%%%%%%%%%%%%%%%%%%%%%%%
%%%%%%%%%%%%%%%%%%%%%%%%%%%%%%%%%%%%%%%%%%%%%%%%%%%%%%%%%%%%%%%%%%%%%%%%%%%%%%

We adopt the convention $1/\beta_{\nu}=\infty$ if the index $\beta_{\nu}$ defined by~(\ref{eq:defbetanu}) vanishes. The size and large intersection properties of $E_{\nu}(h)$ and $E'_{\nu}(h)$ are nontrivial only when $h\in[0,1/\beta_{\nu})$. In fact, when $h\geq 1/\beta_{\nu}$, this description follows essentially from the next result, which is proven in Section~\ref{subsec:proofprpsizeEnuhsimple}.

\begin{prp}\label{prp:sizeEnuhsimple}
Let $\nu$ be an admissible measure. Then,
\[
\as \quad \forall h\in[1/\beta_{\nu},\infty] \qquad E'_{\nu}(h)=\R^{d}\setminus J_{\nu}.
\]
Moreover, $E_{\nu}(1/\beta_{\nu})$ has full Lebesgue measure in $\R^{d}$ with probability one, and
\[
\as \quad \forall h\in(1/\beta_{\nu},\infty] \qquad E_{\nu}(h)=\emptyset.
\]
\end{prp}

Given that $J_{\nu}$ has Lebesgue measure zero, Proposition~\ref{prp:sizeEnuhsimple} ensures that $E_{\nu}(1/\beta_{\nu})$ and $E'_{\nu}(1/\beta_{\nu})$ have full Lebesgue measure in $\R^{d}$ almost surely. This result has direct implications in terms of Hausdorff measures and large intersection properties, which we shall detail when appropriate (see the comments following Theorems~\ref{thm:sizepropEnuh} and~\ref{thm:lipropEnuh}).

%%%%%%%%%%%%%%%%%%%%%%%%%%%%%%%%%%%%%%%%%%%%%%%%%%%%%%%%%%%%%%%%%%%%%%%%%%%%%%
%%%%%%%%%%%%%%%%%%%%%%%%%%%%%%%%%%%%%%%%%%%%%%%%%%%%%%%%%%%%%%%%%%%%%%%%%%%%%%
\subsection{Size properties of the sets $E_{\nu}(h)$ and $E'_{\nu}(h)$: Hausdorff measures and dimension}\label{subsec:sizepropEnuh}
%%%%%%%%%%%%%%%%%%%%%%%%%%%%%%%%%%%%%%%%%%%%%%%%%%%%%%%%%%%%%%%%%%%%%%%%%%%%%%
%%%%%%%%%%%%%%%%%%%%%%%%%%%%%%%%%%%%%%%%%%%%%%%%%%%%%%%%%%%%%%%%%%%%%%%%%%%%%%

We call a gauge function any nondecreasing function $g$ defined on $[0,\infty)$ such that $\lim_{0+}g=g(0)=0$ and $r\mapsto g(r)/r^{d}$ is positive and nonincreasing near zero (this last assumption is not particularly restrictive and may be removed using~\cite[Proposition~2]{Durand:2007uq}), and we let $\gauge_{d}$ denote the set of gauge functions. For any $g\in\gauge_{d}$, the Hausdorff $g$-measure of a subset $F$ of $\R^{d}$ is defined by
\[
\hau^{g}(F)=\lim_{\delta\downarrow 0}\uparrow \hau^{g}_{\delta}(F) \qquad\mbox{with}\qquad \hau^{g}_{\delta}(F)=\inf_{F\subseteq\bigcup_{n} U_{n} \atop \diam{U_{n}}<\delta} \sum_{n=1}^{\infty} g(\diam{U_{n}}).
\]
Here, the infimum is taken over all sequences $(U_{n})_{n\geq 1}$ of subsets of $\R^{d}$ satisfying $F\subseteq\bigcup_{n} U_{n}$ and $\diam{U_{n}}<\delta$ for all $n\geq 1$, where $\diam{\,\cdot\,}$ denotes diameter. As mentioned in~\cite{Rogers:1970wb}, $\hau^{g}$ is a Borel measure on $\R^{d}$. For simplicity, the Hausdorff measure corresponding to a gauge function of the form $r\mapsto r^{d-1}g(r)$ with $g\in\gauge_{1}$ (obtained by letting $d=1$ in the definition of $\gauge_{d}$) is denoted by $\hau^{d-1,g}$.

When $g$ is of the form $r\mapsto r^{s}$ for some $s\in(0,d]$, the Hausdorff $g$-measure is referred to as the $s$-dimensional Hausdorff measure and denoted by $\hau^{s}$. These particular measures lead to the notion of Hausdorff dimension. Specifically, the Hausdorff dimension of a nonempty set $F\subseteq\R^{d}$ is defined by
\[
\Hdim F=\sup\{ s\in (0,d) \:|\: \hau^{s}(F)=\infty \}=\inf\{ s\in (0,d) \:|\: \hau^{s}(F)=0 \},
\]
thus giving an abridged description of its size properties, see~\cite{Falconer:2003oj}. Here, we adopt the convention that $\sup\emptyset=0$ and $\inf\emptyset=d$. We also assume that $\Hdim\emptyset=-\infty$.

%%%%%%%%%%%%%%%%%%%%%%%%%%%%%%%%%%%%%%%%%%%%%%%%%%%%%%%%%%%%%%%%%%%%%%%%%%%%%%
%%%%%%%%%%%%%%%%%%%%%%%%%%%%%%%%%%%%%%%%%%%%%%%%%%%%%%%%%%%%%%%%%%%%%%%%%%%%%%
\subsubsection{General results}
%%%%%%%%%%%%%%%%%%%%%%%%%%%%%%%%%%%%%%%%%%%%%%%%%%%%%%%%%%%%%%%%%%%%%%%%%%%%%%
%%%%%%%%%%%%%%%%%%%%%%%%%%%%%%%%%%%%%%%%%%%%%%%%%%%%%%%%%%%%%%%%%%%%%%%%%%%%%%

When $h\in[0,1/\beta_{\nu})$, the size properties of $E_{\nu}(h)$ and $E'_{\nu}(h)$ are described by the next result. In its statement, $h_{\nu}(g)$ is given by
\begin{equation}\label{eq:defhnug}
h_{\nu}(g)=\inf\Biggl\{ h>0 \:\Biggl|\: \int_{s\in\S^{d-1}\atop x\in(0,1]} g(x^{1/h}) \,\nu(\dd s,\dd x)=\infty \Biggr\},
\end{equation}
for $g\in\gauge_{1}$, with the usual convention that $\inf\emptyset=\infty$. Clearly, $h_{\nu}(g)\leq 1/\beta_{\nu}$.

\begin{thm}\label{thm:sizepropEnuh}
Let $\nu$ be an admissible measure. Then, with probability one, for any real $h\in [0,1/\beta_{\nu})$, any gauge function $g\in\gauge_{1}$ and any nonempty open $W\subseteq\R^{d}$,
\[
\hau^{d-1,g}(E_{\nu}(h)\cap W)=\left\{\begin{array}{ll}
0 & \mbox{if }h<h_{\nu}(g) \\
\infty & \mbox{if }h=h_{\nu}(g)
\end{array}\right.
\]
and
\[
\hau^{d-1,g}(E'_{\nu}(h)\cap W)=\left\{\begin{array}{ll}
0 & \mbox{if }h<h_{\nu}(g) \\
\infty & \mbox{if }h\geq h_{\nu}(g).
\end{array}\right.
\]
\end{thm}

We refer to Section~\ref{subsec:proofslipropEnuh} for a proof of Theorem~\ref{thm:sizepropEnuh}. What is remarkable in this result is that the almost sure event on which its statement holds does not depend on the real $h$, the gauge function $g$ or the open set $W$. In other words, the previous description of the size properties of the sets $E_{\nu}(h)$ and $E'_{\nu}(h)$ holds for almost every sample function of the random field $L_{\nu}$. On top of that, the description is both precise and local, in the sense that we do not restrict our attention to the mere gauge functions of the form $r\mapsto r^{s}$ that lead to the Hausdorff dimension, or to the case in which the open set $W$ is equal to the whole space $\R^{d}$.

With the help of Proposition~\ref{prp:sizeEnuhsimple} above, it is easy to obtain an analog of Theorem~\ref{thm:sizepropEnuh} for the case where $h\geq 1/\beta_{\nu}$. Specifically, with probability one, the sets $E_{\nu}(1/\beta_{\nu})$ and $E'_{\nu}(1/\beta_{\nu})$ both have full Lebesgue measure in $\R^{d}$, so that
\[
\forall g\in\gauge_{d} \quad \forall W\mbox{ open} \qquad
\left.\begin{array}{r}
\hau^{g}(E_{\nu}(1/\beta_{\nu})\cap W)\\
\hau^{g}(E'_{\nu}(1/\beta_{\nu})\cap W)
\end{array}\right\}=\hau^{g}(W).
\]
The same result holds for the sets $E'_{\nu}(h)$, with $h\in(1/\beta_{\nu},\infty]$, because they all coincide with $E'_{\nu}(1/\beta_{\nu})$. For those values of $h$, the study of the size of the sets $E_{\nu}(h)$ is pointless because they are empty almost surely.

In the remainder of this section, we make additional assumptions on the positivity of $\beta_{\nu}$ or the finiteness of $\nu$ with a view to obtaining further results. Due to~(\ref{eq:defbetanu}), it is clear that $\beta_{\nu}=0$ when $\nu$ has finite mass. This corresponds to the case of a compound Poisson field.

%%%%%%%%%%%%%%%%%%%%%%%%%%%%%%%%%%%%%%%%%%%%%%%%%%%%%%%%%%%%%%%%%%%%%%%%%%%%%%
%%%%%%%%%%%%%%%%%%%%%%%%%%%%%%%%%%%%%%%%%%%%%%%%%%%%%%%%%%%%%%%%%%%%%%%%%%%%%%
\subsubsection{Case where $\beta_{\nu}>0$}
%%%%%%%%%%%%%%%%%%%%%%%%%%%%%%%%%%%%%%%%%%%%%%%%%%%%%%%%%%%%%%%%%%%%%%%%%%%%%%
%%%%%%%%%%%%%%%%%%%%%%%%%%%%%%%%%%%%%%%%%%%%%%%%%%%%%%%%%%%%%%%%%%%%%%%%%%%%%%

Here, Theorem~\ref{thm:sizepropEnuh} leads to the following more compact result, which is established in Section~\ref{subsec:proofcorsizepropEnuh}. In its statement, $\gauge_{1}^{*}$ is the set of gauge functions $g\in\gauge_{1}$ for which $\lim_{r\to 0}(\log g(r))/\log r$ exists.

\begin{cor}\label{cor:sizepropEnuh}
Let $\nu$ be an admissible measure satisfying $\beta_{\nu}>0$. Then, with probability one, for any real $h\in [0,1/\beta_{\nu})$, any gauge function $g\in\gauge_{1}^{*}$ and any nonempty open subset $W$ of $\R^{d}$,
\[
\hau^{d-1,g}(E_{\nu}(h)\cap W)=\hau^{d-1,g}(E'_{\nu}(h)\cap W)
=\left\{\begin{array}{ll}
0 & \mbox{if }h<h_{\nu}(g) \\
\infty & \mbox{if }h\geq h_{\nu}(g),
\end{array}\right.
\]
so that
\[
\Hdim(E_{\nu}(h)\cap W)=\Hdim(E'_{\nu}(h)\cap W)=d-1+\beta_{\nu}h.
\]
\end{cor}

In particular, the sets $E_{\nu}(h)\cap W$ have Hausdorff dimension $d-1+\beta_{\nu}h$, with an infinite $(d-1+\beta_{\nu}h)$-dimensional Hausdorff measure.

%%%%%%%%%%%%%%%%%%%%%%%%%%%%%%%%%%%%%%%%%%%%%%%%%%%%%%%%%%%%%%%%%%%%%%%%%%%%%%
%%%%%%%%%%%%%%%%%%%%%%%%%%%%%%%%%%%%%%%%%%%%%%%%%%%%%%%%%%%%%%%%%%%%%%%%%%%%%%
\subsubsection{Case where $\beta_{\nu}=0$ and $\nu$ has infinite total mass}
%%%%%%%%%%%%%%%%%%%%%%%%%%%%%%%%%%%%%%%%%%%%%%%%%%%%%%%%%%%%%%%%%%%%%%%%%%%%%%
%%%%%%%%%%%%%%%%%%%%%%%%%%%%%%%%%%%%%%%%%%%%%%%%%%%%%%%%%%%%%%%%%%%%%%%%%%%%%%

In that situation, we have $h_{\nu}(r\mapsto r^{s})=\infty$ for all $s\in(0,1]$ and, borrowing ideas from the proof of~\cite[Proposition~5]{Durand:2007fk}, we may build a gauge function $g\in\gauge_{1}$ such that $h_{\nu}(g)=0$. Hence, Theorem~\ref{thm:sizepropEnuh} leads to the following statement.

\begin{cor}
Let $\nu$ be a measure with $\beta_{\nu}=0$ and infinite total mass. Then, with probability one, for any real $h\in [0,\infty)$ and any nonempty open subset $W$ of $\R^{d}$,
\[
\Hdim(E'_{\nu}(h)\cap W)=d-1
\]
and the $(d-1)$-dimensional Hausdorff measure of $E'_{\nu}(h)\cap W$ is infinite.
\end{cor}

Let us now focus on the iso-H\"older sets $E_{\nu}(h)$. Theorem~\ref{thm:sizepropEnuh}, coupled with the preceding remarks, also implies that with probability one, for any real $h\in [0,\infty)$ and any nonempty open subset $W$ of $\R^{d}$,
\begin{equation}\label{eq:dimEnuhbetanu0}
\Hdim(E_{\nu}(h)\cap W)\leq d-1,
\end{equation}
However, we cannot conclude that~(\ref{eq:dimEnuhbetanu0}) is an equality, except when $h=0$. According to our approach, inferring that equality holds in~(\ref{eq:dimEnuhbetanu0}) demands that we build, for any given $h\in [0,\infty)$, a gauge function $g\in\gauge_{1}$ with $h_{\nu}(g)=h$, and apply Theorem~\ref{thm:sizepropEnuh} with this gauge function. As mentioned above, such a construction is feasible when $h=0$, but is not always possible otherwise. Indeed, some measures $\nu$ yield $h_{\nu}(g)\in\{0,\infty\}$ for all $g\in\gauge_{1}$. This is the case of the product of the uniform measure on $\S^{d-1}$ and the sum over $j\in\Z\setminus\{0\}$ of the atoms $\delta_{1/j}/|j|$ on $\R^{*}$.

%%%%%%%%%%%%%%%%%%%%%%%%%%%%%%%%%%%%%%%%%%%%%%%%%%%%%%%%%%%%%%%%%%%%%%%%%%%%%%
%%%%%%%%%%%%%%%%%%%%%%%%%%%%%%%%%%%%%%%%%%%%%%%%%%%%%%%%%%%%%%%%%%%%%%%%%%%%%%
\subsubsection{Case where $\nu$ has finite total mass}
%%%%%%%%%%%%%%%%%%%%%%%%%%%%%%%%%%%%%%%%%%%%%%%%%%%%%%%%%%%%%%%%%%%%%%%%%%%%%%
%%%%%%%%%%%%%%%%%%%%%%%%%%%%%%%%%%%%%%%%%%%%%%%%%%%%%%%%%%%%%%%%%%%%%%%%%%%%%%

Making use of~(\ref{eq:Lnu0}) and~(\ref{eq:Lnujalt}), we get
\begin{equation}\label{eq:Lnufinite}
\as \quad \forall t\in\R^{d} \qquad L_{\nu}(t)=\sum_{n=1}^{\infty} X_{n}\ind_{\{\Rho_{n}<\inpr{S_{n}}{t}\}}-\int_{s\in\S^{d-1}\atop x\in (0,1]} x\inpr{s}{t}\,\nu(\dd s,\dd x).
\end{equation}
With probability one, the above sum is piecewise constant, with jumps located on the set $J_{\nu}$ given by~(\ref{eq:defJnu}). Moreover, the integral in~(\ref{eq:Lnufinite}) depends linearly on $t$. Thus, it is natural that the previous results lead to the next statement.

\begin{prp}\label{prp:Enuhnufinite}  Let $\nu$ be a measure with finite total mass. Then, with probability one, for any real $h\in [0,\infty)$,
\[
E'_{\nu}(h)=\emptyset \qquad\mbox{and}\qquad
E_{\nu}(h)=
\begin{cases}
J_{\nu} & \mbox{if } h=0\\
\emptyset & \mbox{if }h>0.
\end{cases}
\]
\end{prp}

Although Proposition~\ref{prp:Enuhnufinite} is elementary, a formal proof is given in Section~\ref{subsec:proofprpEnuhnufinite} for the sake of completeness. By virtue of Proposition~\ref{prp:sizeEnuhsimple}, we also have $E_{\nu}(\infty)=E'_{\nu}(\infty)=\R^{d}\setminus J_{\nu}$ almost surely. Hence,
\[
\as \quad \forall h\in[0,\infty] \qquad \Hdim E_{\nu}(h)=
\begin{cases}
d-1 & \mbox{if } h=0\\
-\infty & \mbox{if } 0<h<\infty\\
d & \mbox{if } h=\infty.
\end{cases}
\]
In the present case, the sample paths of $L_{\nu}$ are not homogeneous. Indeed, with probability one, for any $A\geq 1$, the ball $\clball{0}{A}$ intersects only finitely many hyperplanes $H_{n}$. Thus, one may find a nonempty open set $W$ such that $J_{\nu}\cap W=\emptyset$, which implies that $d_{L_{\nu}}(0,W)=-\infty\neq d-1=d_{L_{\nu}}(0,\R^{d})$.

%%%%%%%%%%%%%%%%%%%%%%%%%%%%%%%%%%%%%%%%%%%%%%%%%%%%%%%%%%%%%%%%%%%%%%%%%%%%%%
%%%%%%%%%%%%%%%%%%%%%%%%%%%%%%%%%%%%%%%%%%%%%%%%%%%%%%%%%%%%%%%%%%%%%%%%%%%%%%
\subsection{Large intersection properties of the sets $E'_{\nu}(h)$}\label{subsec:largeintpropEnuh}
%%%%%%%%%%%%%%%%%%%%%%%%%%%%%%%%%%%%%%%%%%%%%%%%%%%%%%%%%%%%%%%%%%%%%%%%%%%%%%
%%%%%%%%%%%%%%%%%%%%%%%%%%%%%%%%%%%%%%%%%%%%%%%%%%%%%%%%%%%%%%%%%%%%%%%%%%%%%%

As shown below, the sets $E'_{\nu}(h)$ defined by~(\ref{eq:defEnuhp}) are sets with large intersection. Specifically, they belong to certain classes $\lic{g}{W}$ of subsets of $\R^{d}$, which were introduced in~\cite{Durand:2007uq} in order to generalize the original classes of sets with large intersection defined by Falconer~\cite{Falconer:1994hx}. Given a gauge function $g\in\gauge_{d}$ and a nonempty open set $W\subseteq\R^{d}$, the class $\lic{g}{W}$ may be defined in the following manner. To begin with, recall that the function $r\mapsto g(r)/r^d$ is positive and nonincreasing in a neighborhood of zero. Let $\eps_{g}$ denote the supremum of all $\eps\in(0,1]$ such that this property holds on the interval $(0,\eps]$. Moreover, let $\Lambda_{g}$ be the set of all dyadic cubes of diameter less than $\eps_{g}$, that is, sets of the form $\lambda=2^{-j}(k+[0,1)^d)$, where $j$ is an integer larger than $\log_{2}(\sqrt{d}/\eps_{g})$ and $k\in\Z^d$. The outer net measure associated with $g\in\gauge_{d}$ is defined by
\[
\forall F\subseteq\R^{d} \qquad \netm^{g}_{\infty}(F) = \inf_{(\lambda_{n})_{n\geq 1}} \sum_{n=1}^{\infty} g(\diam{\lambda_{n}}),
\]
where the infimum is over all sequences $(\lambda_{n})_{n\geq 1}$ in $\Lambda_{g}\cup\{\emptyset\}$ satisfying $F\subseteq\bigcup_{n}\lambda_{n}$. The outer measure $\netm^{g}_{\infty}$ is intimately related with the corresponding Hausdorff measure  $\hau^{g}$, so that in particular $\hau^{g}(F)>0$ for any set $F\subseteq\R^{d}$ with $\netm^{g}_{\infty}(F)>0$, see~\cite[Theorem~49]{Rogers:1970wb}. In addition, for $\overline{g},g\in\gauge_{d}$, let us write $\overline{g}\prec g$ if $\overline{g}/g$ tends monotonically to infinity at zero. This means essentially that $\overline{g}$ grows faster than $g$ near the origin. We can now define the class $\lic{g}{W}$. Recall that a $G_{\delta}$-set is one that may be expressed as a countable intersection of open sets.

\begin{df}
Let $g\in\gauge_{d}$ and let $W$ be a nonempty open subset of $\R^d$. The class $\lic{g}{W}$ of {\em subsets of $\R^d$ with large intersection in $W$ with respect to $g$} is the collection of all $G_{\delta}$-subsets $F$ of $\R^d$ such that
\[
\forall U\subseteq W\mbox{ open} \quad \forall\overline{g}\in\gauge_{d} \qquad \overline{g}\prec g \quad\Longrightarrow\quad \netm^{\overline{g}}_{\infty}(F\cap U)=\netm^{\overline{g}}_{\infty}(U).
\]
\end{df}

The class $\lic{g}{W}$ satisfies several remarkable properties which are detailed in~\cite{Durand:2007uq}. We collect the most significant ones in the following statement.

\begin{thm}\label{thm:sli}
Let $g\in\gauge_{d}$ and let $W$ be a nonempty open subset of $\R^d$. Then,
\begin{enumerate}
\item\label{item:sliclosedint} the class $\lic{g}{W}$ is closed under countable intersections;
\item the set $f^{-1}(F)$ belongs to $\lic{g}{W}$ for any bi-Lipschitz mapping $f:W\to\R^{d}$ and any set $F\in\lic{g}{f(W)}$;
\item\label{item:slisize} any set $F\in\lic{g}{W}$ satisfies $\hau^{\overline{g}}(F\cap W)=\infty$ for any $\overline{g}\in\gauge_{d}$ with $\overline{g}\prec g$;
\item\label{item:slifull} any $G_{\delta}$-subset of $\R^{d}$ with full Lebesgue measure in $W$ belongs to $\lic{g}{W}$.
\end{enumerate}
\end{thm}

Making use of Theorem~\ref{thm:sli}(\ref{item:slisize}), it is easy to check that any set that belongs to the class $\lic{g}{W}$ is of Hausdorff dimension at least
\begin{equation}\label{eq:defsigmag}
\sigma_{g}=\sup\{ s\in (0,d) \:|\: (r\mapsto r^s)\prec g \},
\end{equation}
with the convention that $\sup\emptyset=0$. In other words, the fact that a set satisfies a large intersection property leads to a lower bound on its Hausdorff dimension, an information which is usually difficult to derive.

More generally, Theorem~\ref{thm:sli}(\ref{item:slisize}) embodies the connection between size and large intersection properties, thereby suggesting a link between the following result and Theorem~\ref{thm:sizepropEnuh}. We shall make this link more apparent when proving these two theorems concurrently in Section~\ref{subsec:proofslipropEnuh}. In the next statement, $\lic{d-1,g}{W}$ denotes the class of sets with large intersection in $W$ with respect to $r\mapsto r^{d-1}g(r)$, where $g\in\gauge_{1}$.

\begin{thm}\label{thm:lipropEnuh}
Let $\nu$ be an admissible measure. With probability one, for any real $h\in [0,1/\beta_{\nu})$, any gauge function $g\in\gauge_{1}$ and any nonempty open subset $W$ of $\R^{d}$,
\[
E'_{\nu}(h)\in\lic{d-1,g}{W} \qquad\Longleftrightarrow\qquad h\geq h_{\nu}(g).
\]
\end{thm}

The previous result only concerns the case where $h<1/\beta_{\nu}$. As regards the opposite case, recall that $E'_{\nu}(1/\beta_{\nu})$ has full Lebesgue measure in $\R^{d}$ with probability one, by Proposition~\ref{prp:sizeEnuhsimple}. Hence, owing to Theorem~\ref{thm:sli}(\ref{item:slifull}),
\[
\as \quad \forall g\in\gauge_{d} \quad \forall W\neq\emptyset\mbox{ open} \qquad E'_{\nu}(1/\beta_{\nu})\in\lic{g}{W}.
\]
The sets $E'_{\nu}(h)$, for $h\in(1/\beta_{\nu},\infty]$, satisfy the same property because they are all identical to $E'_{\nu}(1/\beta_{\nu})$.

%%%%%%%%%%%%%%%%%%%%%%%%%%%%%%%%%%%%%%%%%%%%%%%%%%%%%%%%%%%%%%%%%%%%%%%%%%%%%%
%%%%%%%%%%%%%%%%%%%%%%%%%%%%%%%%%%%%%%%%%%%%%%%%%%%%%%%%%%%%%%%%%%%%%%%%%%%%%%
\subsection{Size properties of the sets $E'_{\nu}(h)$: packing dimension}\label{subsec:packdimEnuh}
%%%%%%%%%%%%%%%%%%%%%%%%%%%%%%%%%%%%%%%%%%%%%%%%%%%%%%%%%%%%%%%%%%%%%%%%%%%%%%
%%%%%%%%%%%%%%%%%%%%%%%%%%%%%%%%%%%%%%%%%%%%%%%%%%%%%%%%%%%%%%%%%%%%%%%%%%%%%%

We do not recall the definition of the packing dimension $\Pdim F$ of a subset $F$ of $\R^{d}$ here, and refer to~\cite[Chapter~3]{Falconer:2003oj} or~\cite[Chapter~5]{Mattila:1995fk} for a full exposition. The fact that the sets $E'_{\nu}(h)$ satisfy a large intersection property has a straightforward consequence on the value of their packing dimension, as we now explain.

First, in view of Proposition~\ref{prp:sizeEnuhsimple}, it is clear that with probability one, the sets $E'_{\nu}(h)$, for $h\geq 1/\beta_{\nu}$, all have packing dimension equal to $d$. We may therefore restrict our attention to the case in which $h<1/\beta_{\nu}$ in what follows.

Then, let us assume that $\nu$ has infinite total mass and that $d>1$. As mentioned above, there exists a gauge function $g\in\gauge_{1}$ such that $h_{\nu}(g)=0$. Thanks to Theorem~\ref{thm:lipropEnuh}, almost surely, for every $h\in[0,1/\beta_{\nu})$, the set $E'_{\nu}(h)$ belongs to the class $\lic{\widetilde{g}}{\R^{d}}$, where $\widetilde{g}:r\mapsto r^{d-1}g(r)$. It is shown in~\cite{Durand:2007uq} that if the parameter $\sigma_{\widetilde{g}}$ defined by~(\ref{eq:defsigmag}) is positive, then this class is included in the class $\grintfalc^{\sigma_{\widetilde{g}}}$ introduced by Falconer~\cite{Falconer:1994hx}. This is clearly the case here, since $\sigma_{\widetilde{g}}\geq d-1>0$. By virtue of~\cite[Theorem~D]{Falconer:1994hx}, every set of the latter class has packing dimension $d$ in every nonempty open subset of $\R^{d}$. Hence,  we end up with the next statement.

\begin{cor}
If $d>1$, then for any admissible measure $\nu$ with infinite total mass,
\[
\as \quad \forall h\in[0,1/\beta_{\nu}) \quad \forall W\neq\emptyset\mbox{ open} \qquad \Pdim(E'_{\nu}(h)\cap W)=d.
\]
\end{cor}

The preceding result remains valid for $d=1$ under the additional assumption that $\beta_{\nu}>0$, but for $h\in(0,1/\beta_{\nu})$ only. This follows from applying the previous method with the gauge function $\widetilde{g}:r\mapsto r^{\beta_{\nu}h}$, which satisfies $\sigma_{\widetilde{g}}=\beta_{\nu}h>0$.

%%%%%%%%%%%%%%%%%%%%%%%%%%%%%%%%%%%%%%%%%%%%%%%%%%%%%%%%%%%%%%%%%%%%%%%%%%%%%%
%%%%%%%%%%%%%%%%%%%%%%%%%%%%%%%%%%%%%%%%%%%%%%%%%%%%%%%%%%%%%%%%%%%%%%%%%%%%%%
%%%%%%%%%%%%%%%%%%%%%%%%%%%%%%%%%%%%%%%%%%%%%%%%%%%%%%%%%%%%%%%%%%%%%%%%%%%%%%
\section{Superposition of the Gaussian and the jump components}\label{sec:Gaussianplusjumps}
%%%%%%%%%%%%%%%%%%%%%%%%%%%%%%%%%%%%%%%%%%%%%%%%%%%%%%%%%%%%%%%%%%%%%%%%%%%%%%
%%%%%%%%%%%%%%%%%%%%%%%%%%%%%%%%%%%%%%%%%%%%%%%%%%%%%%%%%%%%%%%%%%%%%%%%%%%%%%
%%%%%%%%%%%%%%%%%%%%%%%%%%%%%%%%%%%%%%%%%%%%%%%%%%%%%%%%%%%%%%%%%%%%%%%%%%%%%%

The next lemma gives an expression of the H\"older exponent of an admissible canonical L\'evy field $Y_{a,\mu,\nu}=\inpr{a}{\cdot}+B_{\mu}+L_{\nu}$, in terms of that of its Gaussian component $B_{\mu}$ and that of its jump component $L_{\nu}$. In particular, it allows one to easily infer the spectrum of singularities of $Y_{a,\mu,\nu}$ from those of $B_{\mu}$ and $L_{\nu}$, that is, to deduce Corollary~\ref{cor:specgeneral} from Theorems~\ref{thm:HolderBmu} and~\ref{thm:specLnu}.

\begin{lem}\label{lem:superpos}
Let $a\in\R^{d}$, let $\mu$ be a finite symmetric measure on $\S^{d-1}$, and let $\nu$ be an admissible measure. Then,
\begin{equation}\label{eq:lowbouholexp}
\as \quad \forall t \in \R^{d} \qquad \alpha_{Y_{a,\mu,\nu}}(t)=\alpha_{B_{\mu}}(t)\wedge \alpha_{L_{\nu}}(t).
\end{equation}
\end{lem}

\begin{proof}
As linear drifts are $C^{\infty}$ everywhere, we may clearly assume that $a$ is zero. Moreover, note that the H\"older exponent of the sum of two functions is the minimum of the two exponents, except maybe when these exponents coincide, in which case the exponent of the sum may be larger. In view of Theorem~\ref{thm:HolderBmu}, it remains to show that with probability one, $\alpha_{Y_{0,\mu,\nu}}(t)\leq 1/2$ for any $t\in\R^{d}$ with $\alpha_{L_{\nu}}(t)=1/2$.

To this end, it suffices to observe that Proposition~\ref{prp:uppbouhold} below (which corresponds to the upper bound in Proposition~\ref{prp:characLnuKnualpha} below) still holds when replacing $L_{\nu}$ by $B_{\mu}+L_{\nu}$, because $B_{\mu}$ is continuous everywhere. Therefore, with probability one, $\alpha_{Y_{0,\mu,\nu}}(t)\leq\Alpha_{\nu}(t)$ for any $t\in\R^{d}$ with $\alpha_{L_{\nu}}(t)=1/2$. We conclude by Proposition~\ref{prp:characLnuKnualpha}.
\end{proof}

%%%%%%%%%%%%%%%%%%%%%%%%%%%%%%%%%%%%%%%%%%%%%%%%%%%%%%%%%%%%%%%%%%%%%%%%%%%%%%
%%%%%%%%%%%%%%%%%%%%%%%%%%%%%%%%%%%%%%%%%%%%%%%%%%%%%%%%%%%%%%%%%%%%%%%%%%%%%%
%%%%%%%%%%%%%%%%%%%%%%%%%%%%%%%%%%%%%%%%%%%%%%%%%%%%%%%%%%%%%%%%%%%%%%%%%%%%%%
\section{Architecture of the proofs concerning the jump component}\label{sec:architecture}
%%%%%%%%%%%%%%%%%%%%%%%%%%%%%%%%%%%%%%%%%%%%%%%%%%%%%%%%%%%%%%%%%%%%%%%%%%%%%%
%%%%%%%%%%%%%%%%%%%%%%%%%%%%%%%%%%%%%%%%%%%%%%%%%%%%%%%%%%%%%%%%%%%%%%%%%%%%%%
%%%%%%%%%%%%%%%%%%%%%%%%%%%%%%%%%%%%%%%%%%%%%%%%%%%%%%%%%%%%%%%%%%%%%%%%%%%%%%

Let us now present the key ideas involved in the proofs of the main results of Section~\ref{sec:regnonGaussian}, which describe the size and large intersection properties of the iso-H\"older sets $E_{\nu}(h)$ and the singularity sets $E'_{\nu}(h)$.

%%%%%%%%%%%%%%%%%%%%%%%%%%%%%%%%%%%%%%%%%%%%%%%%%%%%%%%%%%%%%%%%%%%%%%%%%%%%%%
%%%%%%%%%%%%%%%%%%%%%%%%%%%%%%%%%%%%%%%%%%%%%%%%%%%%%%%%%%%%%%%%%%%%%%%%%%%%%%
\subsection{Location of the singularities}\label{subsec:locsing}
%%%%%%%%%%%%%%%%%%%%%%%%%%%%%%%%%%%%%%%%%%%%%%%%%%%%%%%%%%%%%%%%%%%%%%%%%%%%%%
%%%%%%%%%%%%%%%%%%%%%%%%%%%%%%%%%%%%%%%%%%%%%%%%%%%%%%%%%%%%%%%%%%%%%%%%%%%%%%

The first ingredient in our proofs is a precise determination of the location of the singularities of $L_{\nu}$. This follows from a characterization of the H\"older exponent of its sample paths in terms of the atoms of the Poisson measure $\Nu$ arising in the construction described in Section~\ref{sec:defnonGaussianLevy}. In fact, the value of the exponent $\alpha_{L_{\nu}}(t)$ at a point $t\in\R^{d}$ depends on how well $t$ is approximated by the hyperplanes on which $L_{\nu}$ has a jump of size at most one.

To be more specific, for any real $\alpha>0$, let us consider
\[
K_{\nu}(\alpha)=\bigl\{ t\in\R^{d} \:\bigl|\: \dist(t,H_{n})<|X_{n}|^{1/\alpha} \mbox{ for i.m. } n\geq 1 \mbox{ with } |X_{n}|\leq 1 \bigr\},
\]
where i.m.~stands for infinitely many and $\dist(t,H_{n})$ denotes the distance, equal to $|\Rho_{n}-\inpr{S_{n}}{t}|$, between the point $t$ and the hyperplane $H_{n}$ defined by~(\ref{eq:defHn}). In the terminology of Diophantine approximation, $K_{\nu}(\alpha)$ is the set of points that are approximated by the hyperplanes $H_{n}$ with a precision given by $|X_{n}|^{1/\alpha}$. Equivalently,
\[
K_{\nu}(\alpha)=\left\{ t\in\R^{d} \:\Biggl|\: \int_{(\rho,s)\in\H_{d}\atop |x|\in(0,1]}\ind_{\{|\rho-\inpr{s}{t}|<|x|^{1/\alpha}\}} \,\Nu(\dd\rho,\dd s,\dd x)=\infty \right\}.
\]
Given that the mapping $\alpha\mapsto K_{\nu}(\alpha)$ is nondecreasing, it is possible to define
\[
\Alpha_{\nu}(t)=\inf\{\alpha>0 \:|\: t\in K_{\nu}(\alpha) \}
\]
for each $t\in\R^{d}$. Proposition~\ref{prp:KnualphaRd} below implies that with probability one, $K_{\nu}(\alpha)$ covers the whole space $\R^{d}$ when $\alpha>1/\beta_{\nu}$. As a consequence,
\begin{equation}\label{eq:Alphanutleqbetanu}
\as \quad \forall t\in\R^{d} \qquad 0\leq\Alpha_{\nu}(t)\leq 1/\beta_{\nu}.
\end{equation}

The next result, which is proven in Section~\ref{sec:proofprpcharacLnuKnualpha}, gives a simple connection between the value of the H\"older exponent of $L_{\nu}$ at a given point $t$ and that of $\Alpha_{\nu}(t)$.

\begin{prp}\label{prp:characLnuKnualpha}
If $\nu$ is admissible, then
\[
\as \quad \forall t\in\R^{d} \qquad
\alpha_{L_{\nu}}(t)=
\begin{cases}
0 & \mbox{if } t\in J_{\nu}\\
\Alpha_{\nu}(t) & \mbox{else.}
\end{cases}
\]
\end{prp}

Recall that $J_{\nu}$ is the set given by~(\ref{eq:defJnu}) and composed of the points at which the field $L_{\nu}$ jumps. This is why the H\"older exponent of this field vanishes everywhere in $J_{\nu}$. Furthermore, it follows from Proposition~\ref{prp:characLnuKnualpha} that with probability one, the H\"older exponent of $L_{\nu}$ is at most $1/\beta_{\nu}$ everywhere. So, for any $h>1/\beta_{\nu}$,
\begin{equation}\label{eq:sizeEnuhsimple}
E'_{\nu}(h)=\R^{d}\setminus J_{\nu} \qquad\mbox{and}\qquad E_{\nu}(h)=\emptyset,
\end{equation}
thus confirming some results announced in Proposition~\ref{prp:sizeEnuhsimple}. When $h\leq 1/\beta_{\nu}$, Proposition~\ref{prp:characLnuKnualpha} enables us to relate $E_{\nu}(h)$ and $E'_{\nu}(h)$ with the sets $K_{\nu}(\alpha)$ and $J_{\nu}$ as follows.

\begin{cor}\label{cor:characLnuKnualpha}
If $\nu$ is admissible, then with probability one,
\begin{enumerate}
\item for any $h\in[0,1/\beta_{\nu}]$,
\[
E'_{\nu}(h)=\Biggl(\bigcap_{h<\alpha\leq 1/\beta_{\nu}}K_{\nu}(\alpha)\Biggr)\setminus J_{\nu}\,;
\]
\item $E_{\nu}(0)=J_{\nu}\cup E'_{\nu}(0)$ and for any $h\in(0,1/\beta_{\nu}]$,
\[
E_{\nu}(h)=E'_{\nu}(h)\setminus\bigcup_{0<\alpha<h} K_{\nu}(\alpha).
\]
\end{enumerate}
\end{cor}

We adopt here the standard convention that an intersection and a union indexed by the empty set are equal to $\R^{d}$ and the empty set, respectively. In particular, $E'_{\nu}(1/\beta_{\nu})=\R^{d}\setminus J_{\nu}$ almost surely. The proof of Corollary~\ref{cor:characLnuKnualpha}, assuming that Proposition~\ref{prp:characLnuKnualpha} holds, is straightforward and therefore omitted.

%%%%%%%%%%%%%%%%%%%%%%%%%%%%%%%%%%%%%%%%%%%%%%%%%%%%%%%%%%%%%%%%%%%%%%%%%%%%%%
%%%%%%%%%%%%%%%%%%%%%%%%%%%%%%%%%%%%%%%%%%%%%%%%%%%%%%%%%%%%%%%%%%%%%%%%%%%%%%
\subsection{Diophantine approximation by Poisson hyperplanes}\label{subsec:Diophapphyp}
%%%%%%%%%%%%%%%%%%%%%%%%%%%%%%%%%%%%%%%%%%%%%%%%%%%%%%%%%%%%%%%%%%%%%%%%%%%%%%
%%%%%%%%%%%%%%%%%%%%%%%%%%%%%%%%%%%%%%%%%%%%%%%%%%%%%%%%%%%%%%%%%%%%%%%%%%%%%%

Corollary~\ref{cor:characLnuKnualpha} suggests that the proofs of the results of Section~\ref{sec:regnonGaussian} follow from a precise understanding of the size and large intersection properties of the sets $K_{\nu}(\alpha)$. In fact, these proofs make a crucial use of the next two results, see Section~\ref{subsec:proofslipropEnuh}. The first one shows that $K_{\nu}(\alpha)=\R^{d}$ almost surely whenever $\alpha>1/\beta_{\nu}$, thereby leading to~(\ref{eq:Alphanutleqbetanu}). We refer to Section~\ref{sec:proofprpKnualphaRd} for its proof.

\begin{prp}\label{prp:KnualphaRd}
For any real $\alpha>1/\beta_{\nu}$, with probability one, $K_{\nu}(\alpha)=\R^{d}$.
\end{prp}

The second one deals with the size and large intersection properties of $K_{\nu}(\alpha)$ and is proven in Section~\ref{sec:proofthmsliKnualpha}. Recall that $h_{\nu}(g)$ is defined by~(\ref{eq:defhnug}).

\begin{thm}\label{thm:sliKnualpha}
Let $\alpha>0$. Then, with probability one, for any gauge function $g\in\gauge_{1}$,
\[\left\{\begin{array}{l}
h_{\nu}(g)>\alpha \qquad\Longrightarrow\qquad \hau^{d-1,g}(K_{\nu}(\alpha))=0 \\[2mm]
h_{\nu}(g)<\alpha \qquad\Longrightarrow\qquad K_{\nu}(\alpha)\in\lic{d-1,g}{\R^{d}}.
\end{array}\right.
\]
\end{thm}

The fact that we make use of that result in the proof of Theorems~\ref{thm:sizepropEnuh} and~\ref{thm:lipropEnuh} hints at why they describe the size and large intersection properties of the set $E_{\nu}(h)$ and $E'_{\nu}(h)$ only in terms of the gauge functions of the form $r\mapsto r^{d-1}g(r)$ with $g\in\gauge_{1}$. This comes from the ubiquity techniques that we use in Section~~\ref{sec:proofthmsliKnualpha} (see Lemma~\ref{lem:ubiqFAilg} in particular), and is understandable because $K_{\nu}(\alpha)$ may be regarded as an enlargement of the random hyperplanes $H_{n}$, which are of dimension $d-1$.

%%%%%%%%%%%%%%%%%%%%%%%%%%%%%%%%%%%%%%%%%%%%%%%%%%%%%%%%%%%%%%%%%%%%%%%%%%%%%%
%%%%%%%%%%%%%%%%%%%%%%%%%%%%%%%%%%%%%%%%%%%%%%%%%%%%%%%%%%%%%%%%%%%%%%%%%%%%%%
%%%%%%%%%%%%%%%%%%%%%%%%%%%%%%%%%%%%%%%%%%%%%%%%%%%%%%%%%%%%%%%%%%%%%%%%%%%%%%
\section{Location of the singularities of the jump component}\label{sec:proofprpcharacLnuKnualpha}
%%%%%%%%%%%%%%%%%%%%%%%%%%%%%%%%%%%%%%%%%%%%%%%%%%%%%%%%%%%%%%%%%%%%%%%%%%%%%%
%%%%%%%%%%%%%%%%%%%%%%%%%%%%%%%%%%%%%%%%%%%%%%%%%%%%%%%%%%%%%%%%%%%%%%%%%%%%%%
%%%%%%%%%%%%%%%%%%%%%%%%%%%%%%%%%%%%%%%%%%%%%%%%%%%%%%%%%%%%%%%%%%%%%%%%%%%%%%

%%%%%%%%%%%%%%%%%%%%%%%%%%%%%%%%%%%%%%%%%%%%%%%%%%%%%%%%%%%%%%%%%%%%%%%%%%%%%%
%%%%%%%%%%%%%%%%%%%%%%%%%%%%%%%%%%%%%%%%%%%%%%%%%%%%%%%%%%%%%%%%%%%%%%%%%%%%%%
\subsection{Preliminaries}
%%%%%%%%%%%%%%%%%%%%%%%%%%%%%%%%%%%%%%%%%%%%%%%%%%%%%%%%%%%%%%%%%%%%%%%%%%%%%%
%%%%%%%%%%%%%%%%%%%%%%%%%%%%%%%%%%%%%%%%%%%%%%%%%%%%%%%%%%%%%%%%%%%%%%%%%%%%%%

The purpose of this section is to establish Proposition~\ref{prp:characLnuKnualpha}. The proof relies on suitable estimates of the increments of the random fields $L_{\nu,j}$ that come into play in the construction of $L_{\nu}$, as detailed in Section~\ref{sec:defnonGaussianLevy}. To be specific, for all integers $A,j,k\geq 1$, let
\begin{equation}\label{eq:defzetanuAjk}
\zeta_{\nu}(A,j,k)=\sup_{\|t\|\leq A \atop \|\tau\|\leq 2^{-k}} |L_{\nu,j}(t+\tau)-L_{\nu,j}(t)|.
\end{equation}
Even though the supremum is taken over an uncountable set of parameters, there is no measurability issue here, in the sense that $\zeta_{\nu}(A,j,k)$ is a random variable. In fact, it is easy to check that the field $L_{\nu,j}$ is separable and that $\Q^{d}$ may be taken as the {\em separant} dense countable subset of $\R^{d}$ involved in this property. Accordingly, in addition to~(\ref{eq:defzetanuAjk}), we have
\begin{equation}\label{eq:defzetanuAjkalt}
\zeta_{\nu}(A,j,k)=\sup_{t\in\Q^{d}\cap\clball{0}{A} \atop \tau\in\Q^{d}\cap\clball{0}{2^{-k}}} |L_{\nu,j}(t+\tau)-L_{\nu,j}(t)|.
\end{equation}

The next lemma yields an upper bound on $\zeta_{\nu}(A,j,k)$, that is, a control of the increments of $L_{\nu,j}$. It is proven in Section~\ref{sec:prooflemoscLnuj}, and also comes into play in the proof of Proposition~\ref{prp:existpathsLnu}, see Section~\ref{subsec:proofprpexistpathsLnu}. Recall that $(\nu_{j})_{j\geq 1}$ is the sequence arising in~(\ref{eq:defchinu}).

\begin{lem}\label{lem:oscLnuj}
For every integer $A\geq 1$, with probability one,
\[
\Zeta_{\nu}(A)=\sup_{(j,k)\in\N^{2}}\frac{\zeta_{\nu}(A,j,k)}{2^{-j}k(j+2^{-k/2}(j\nu_{j})^{1/2})}<\infty.
\]
\end{lem}

In view of Lemma~\ref{lem:oscLnuj}, we shall work below on the almost sure event consisting in the finiteness of $\Zeta_{\nu}(A)$, for all $A\geq 1$. Furthermore, recall that $(\Rho_{n},S_{n},X_{n})$, for $n\geq 1$, are the atoms of the Poisson random measure $\Nu$ arising in the construction of $L_{\nu}$. We may assume that the reals $\Rho_{n}$ are distinct. Indeed, for any $j\geq 1$, the image under $(\rho,s,x)\mapsto\rho$ of the restriction of $\Nu$ to $(0,\infty)\times\S^{d-1}\times (\R\setminus[-2^{-j},2^{-j}])$ is a Poisson measure on $(0,\infty)$ with intensity proportional to $\leb^{1}_{+}$, thereby being almost surely simple, see~\cite[p.~299]{Neveu:1977mz}. In addition, we define $L^{j_{1},j_{2}}_{\nu}=\sum_{j=j_{1}}^{j_{2}}L_{\nu,j}$ for $0\leq j_{1}\leq j_{2}\leq\infty$.

We now split the proof of Proposition~\ref{prp:characLnuKnualpha} into three parts. First, we show that $\alpha_{L_{\nu}}(t)$ vanishes at every jump point $t\in J_{\nu}$. Second, we show that $\alpha_{L_{\nu}}(t)\leq\Alpha_{\nu}(t)$ at every $t\not\in J_{\nu}$. Third, we show that $\alpha_{L_{\nu}}(t)\geq\Alpha_{\nu}(t)$ at any such $t$. Throughout, we assume that the measure $\nu$ is admissible.

%%%%%%%%%%%%%%%%%%%%%%%%%%%%%%%%%%%%%%%%%%%%%%%%%%%%%%%%%%%%%%%%%%%%%%%%%%%%%%
%%%%%%%%%%%%%%%%%%%%%%%%%%%%%%%%%%%%%%%%%%%%%%%%%%%%%%%%%%%%%%%%%%%%%%%%%%%%%%
\subsection{Value of the H\"older exponent at the jump points}
%%%%%%%%%%%%%%%%%%%%%%%%%%%%%%%%%%%%%%%%%%%%%%%%%%%%%%%%%%%%%%%%%%%%%%%%%%%%%%
%%%%%%%%%%%%%%%%%%%%%%%%%%%%%%%%%%%%%%%%%%%%%%%%%%%%%%%%%%%%%%%%%%%%%%%%%%%%%%

With regard to the next statement, recall that the hyperplanes $H_{n}$ are defined by~(\ref{eq:defHn}).

\begin{lem}\label{lem:onejump}
Almost surely, for any $t\in\R^{d}$ such that $t\in H_{n_{0}}$ for a unique $n_{0}\geq 1$,
\[
\lim_{\ell\to\infty} L_{\nu}\left(t+\frac{S_{n_{0}}}{\ell}\right)=L_{\nu}(t)+X_{n_{0}}.
\]
\end{lem}

\begin{proof}
Given an integer $A\geq 1$, let $t$ be a point in $\clball{0}{A}$ such that $t\in H_{n_{0}}$ for a unique $n_{0}\geq 1$. Then, for any integer $\ell$ larger than some $\ell_{0}\geq 1$, the only $n\geq 1$ satisfying $(\Rho_{n},S_{n})\in V_{t+S_{n_{0}}/\ell}\Delta V_{t}$ is $n_{0}$, and in fact $(\Rho_{n_{0}},S_{n_{0}})\in V_{t+S_{n_{0}}/\ell}\setminus V_{t}$. Here, $\Delta$ stands for symmetric difference of sets, and both $V_{t}$ and $V_{t+S_{n_{0}}/\ell}$ are given by~(\ref{eq:defVt}). Hence, for all natural numbers $\ell>\ell_{0}$ and $j_{0}>-\log_{2}|X_{n_{0}}|$, we have
\[
L^{0,j_{0}}_{\nu}\left(t+\frac{S_{n_{0}}}{\ell}\right)-L^{0,j_{0}}_{\nu}(t)=X_{n_{0}}-\frac{1}{\ell}\int_{s\in\S^{d-1}\atop x\in(2^{-j_{0}},1]}x\inpr{s}{S_{n_{0}}}\,\nu(\dd s,\dd x),
\]
thanks to~(\ref{eq:Lnu0}) and~(\ref{eq:Lnujalt}). In addition, given that $\|t\|\leq A$ and $\|S_{n_{0}}/\ell\|\leq 1/2$, we get
\[
\left|L^{j_{0}+1,\infty}_{\nu}\left(t+\frac{S_{n_{0}}}{\ell}\right)-L^{j_{0}+1,\infty}_{\nu}(t)\right|\leq\sum_{j=j_{0}+1}^{\infty}\zeta_{\nu}(A,j,1),
\]
by virtue of~(\ref{eq:defzetanuAjk}). As a consequence, making use of Lemma~\ref{lem:oscLnuj}, we deduce that
\[
\left|L_{\nu}\left(t+\frac{S_{n_{0}}}{\ell}\right)-L_{\nu}(t)-X_{n_{0}}\right|\leq\frac{1}{\ell}\sum_{j=1}^{j_{0}}\nu_{j}+\Zeta_{\nu}(A)\sum_{j=j_{0}+1}^{\infty} 2^{-j}(j+(j\nu_{j})^{1/2}),
\]
and conclude by letting $\ell\to\infty$, and then by letting $j_{0}\to\infty$ while using of the fact that the sum $\chi_{\nu}$ defined by~(\ref{eq:defchinu}) is finite, because $\nu$ is admissible.
\end{proof}

Thanks to Lemma~\ref{lem:onejump}, we may now prove that the H\"older exponent of $L_{\nu}$ vanishes at every jump point.

\begin{prp}\label{prp:Holderjumpoint}
With probability one, $\alpha_{L_{\nu}}(t)=0$ for every $t\in J_{\nu}$.
\end{prp}

\begin{proof}
Let $t\in J_{\nu}$ and suppose that $\alpha_{L_{\nu}}(t)>0$. This means that there are three reals $\eps,\delta,C>0$ such that for any $\tau\in\R^{d}$,
\begin{equation}\label{eq:Holderpos}
\|\tau\|\leq\delta \qquad\Longrightarrow\qquad |L_{\nu}(t+\tau)-L_{\nu}(t)|\leq C\|\tau\|^{\eps}.
\end{equation}

Since $t\in J_{\nu}$, there is an $n_{0}\geq 1$ such that $t\in H_{n_{0}}$. However, $n_{0}$ need not be unique, and we cannot apply Lemma~\ref{lem:onejump} directly. To cope with that problem, recall that the reals $\Rho_{n}$ are distinct, so that the hyperplanes $H_{n}$ are distinct too. Hence, for any integer $m\geq 1$, the set $\opball{t}{1/m}\cap(H_{n_{0}}\setminus\bigcup_{n\neq n_{0}} H_{n})$ contains a point $t_{m}$. Here, $\opball{t}{1/m}$ is the open Euclidean ball centered at $t$ with radius $1/m$. Then, $n_{0}$ is the only integer such that $t_{m}\in H_{n_{0}}$. Applying Lemma~\ref{lem:onejump} with $t_{m}$, we get
\begin{equation}\label{eq:onejumptm}
\forall m\geq 1 \qquad \lim_{\ell\to\infty} L_{\nu}\left(t_{m}+\frac{S_{n_{0}}}{\ell}\right)=L_{\nu}(t_{m})+X_{n_{0}}.
\end{equation}

Now, for all integers $m,\ell\geq 1$, we have $\|t_{m}-t\|<1/m$ and $\|t_{m}+S_{n_{0}}/\ell-t\|<1/m+1/\ell$. So, assuming that $1/m+1/\ell<\delta$ and applying~(\ref{eq:Holderpos}), we obtain
\[
\left|L_{\nu}\left(t_{m}+\frac{S_{n_{0}}}{\ell}\right)-L_{\nu}(t_{m})\right|\leq C\left(\frac{1}{m^{\eps}}+\left(\frac{1}{m}+\frac{1}{\ell}\right)^{\eps}\right).
\]
Letting $\ell\to\infty$ and using~(\ref{eq:onejumptm}), we infer that $|X_{n_{0}}|\leq 2C/m^{\eps}$ for any $m>1/\delta$. Then, letting $m\to\infty$, we get $X_{n_{0}}=0$, which contradicts the fact that $X_{n_{0}}\in\R^{*}$.
\end{proof}

%%%%%%%%%%%%%%%%%%%%%%%%%%%%%%%%%%%%%%%%%%%%%%%%%%%%%%%%%%%%%%%%%%%%%%%%%%%%%%
%%%%%%%%%%%%%%%%%%%%%%%%%%%%%%%%%%%%%%%%%%%%%%%%%%%%%%%%%%%%%%%%%%%%%%%%%%%%%%
\subsection{Upper bound on the H\"older exponent}
%%%%%%%%%%%%%%%%%%%%%%%%%%%%%%%%%%%%%%%%%%%%%%%%%%%%%%%%%%%%%%%%%%%%%%%%%%%%%%
%%%%%%%%%%%%%%%%%%%%%%%%%%%%%%%%%%%%%%%%%%%%%%%%%%%%%%%%%%%%%%%%%%%%%%%%%%%%%%

We now consider the points at which the field $L_{\nu}$ does not jump.

\begin{prp}\label{prp:uppbouhold} 
With probability one, $\alpha_{L_{\nu}}(t)\leq\Alpha_{\nu}(t)$ for every $t\in\R^{d}\setminus J_{\nu}$.
\end{prp}

\begin{proof}
Let $t\in\R^{d}\setminus J_{\nu}$ and $\alpha>\Alpha_{\nu}(t)$, and suppose that $\alpha_{L_{\nu}}(t)\geq\alpha+\eps$ for some $\eps>0$. First, there are $\delta,C>0$ and a polynomial $Q_{t}$ such that for any $\tau\in\R^{d}$,
\begin{equation}\label{eq:Holderalphaeps}
\|\tau\|\leq\delta \qquad\Longrightarrow\qquad |L_{\nu}(t+\tau)-Q_{t}(\tau)|\leq C\|\tau\|^{\alpha+\eps}.
\end{equation}
Second, $t\in K_{\nu}(\alpha)\setminus J_{\nu}$, so that for any $B\in(0,1]$ with $B^{-\alpha/\eps}>3\,C\,2^{\alpha+\eps}$ and $2B^{1/\alpha}<\delta$, there is an $n_{0}\geq 1$ with $|X_{n_{0}}|\leq B$ and $0<\dist(t,H_{n_{0}})\leq|X_{n_{0}}|^{1/\alpha}$.

Let us suppose that for any $\tau\in\R^{d}$,
\begin{equation}\label{eq:contradHolderup}
\|\tau\|<2\,\dist(t,H_{n_{0}}) \qquad\Longrightarrow\qquad 3\,|L_{\nu}(t+\tau)-Q_{t}(\tau)|<|X_{n_{0}}|.
\end{equation}
As in the proof of Proposition~\ref{prp:Holderjumpoint}, the set $\opball{t}{2\,\dist(t,H_{n_{0}})}\cap(H_{n_{0}}\setminus\bigcup_{n\neq n_{0}} H_{n})$ contains a point $t'$. For $\ell$ large enough, $\|t'-t\|$ and $\|t'+S_{n_{0}}/\ell-t\|$ are both smaller than $2\,\dist(t,H_{n_{0}})$, so that~(\ref{eq:contradHolderup}) leads to
\[
\left| L_{\nu}\left(t'+\frac{S_{n_{0}}}{\ell}\right)-L_{\nu}(t') \right| \leq \frac{2}{3}|X_{n_{0}}|+\left|Q_{t}\left(t'+\frac{S_{n_{0}}}{\ell}-t\right)-Q_{t}(t'-t)\right|.
\]
Given that $n_{0}$ is the only integer such that $t'\in H_{n_{0}}$, and due to Lemma~\ref{lem:onejump}, the left-hand side tends to $|X_{n_{0}}|$ as $\ell\to\infty$. This contradicts the fact that the right-hand side goes to $2\,|X_{n_{0}}|/3$. As a result, there is a $\tau\in\R^{d}$ for which~(\ref{eq:contradHolderup}) does not hold.

Therefore, we have $\|\tau\|<2\,\dist(t,H_{n_{0}})\leq 2|X_{n_{0}}|^{1/\alpha}\leq 2B^{1/\alpha}<\delta$ and
\[
|X_{n_{0}}|\leq 3\,|L_{\nu}(t+\tau)-Q_{t}(\tau)|\leq 3C\|\tau\|^{\alpha+\eps}\leq 3\,C\,2^{\alpha+\eps} |X_{n_{0}}|^{1+\eps/\alpha},
\]
thanks to~(\ref{eq:Holderalphaeps}). This implies that $B^{-\alpha/\eps}\leq 3\,C\,2^{\alpha+\eps}$, which contradicts the choice of $B$. Finally, $\alpha_{L_{\nu}}(t)\leq\alpha$, and we conclude by letting $\alpha\downarrow\Alpha_{\nu}(t)$.
\end{proof}

%%%%%%%%%%%%%%%%%%%%%%%%%%%%%%%%%%%%%%%%%%%%%%%%%%%%%%%%%%%%%%%%%%%%%%%%%%%%%%
%%%%%%%%%%%%%%%%%%%%%%%%%%%%%%%%%%%%%%%%%%%%%%%%%%%%%%%%%%%%%%%%%%%%%%%%%%%%%%
\subsection{Lower bound on the H\"older exponent}
%%%%%%%%%%%%%%%%%%%%%%%%%%%%%%%%%%%%%%%%%%%%%%%%%%%%%%%%%%%%%%%%%%%%%%%%%%%%%%
%%%%%%%%%%%%%%%%%%%%%%%%%%%%%%%%%%%%%%%%%%%%%%%%%%%%%%%%%%%%%%%%%%%%%%%%%%%%%%

It remains to establish the next result. Its proof is split into several parts for the sake of clarity.

\begin{prp}
With probability one, $\alpha_{L_{\nu}}(t)\geq\Alpha_{\nu}(t)$ for every $t\in\R^{d}\setminus J_{\nu}$.
\end{prp}

In view of~(\ref{eq:nuLevy}), with probability one, for any integer $A\geq 1$ and any real $\eps>0$, there are only finitely many $n\geq 1$ such that $\Rho_{n}<A$ and $|X_{n}|>\eps$ simultaneously. We may therefore suppose, in addition to the assumptions made at the beginning of this section, that the corresponding almost sure event occurs. We may also assume that the almost sure event given by~(\ref{eq:Alphanutleqbetanu}) occurs too.

Now, given an integer $A\geq 1$, let $t\in\opball{0}{A}\setminus J_{\nu}$. To show that $\alpha_{L_{\nu}}(t)\geq\Alpha_{\nu}(t)$, we may obviously assume that $\Alpha_{\nu}(t)>0$. Then, let $\alpha\in(0,\Alpha_{\nu}(t))$. As $t\not\in K_{\nu}(\alpha)$, there are only finitely many $n\geq 1$ such that $\dist(t,H_{n})<|X_{n}|^{1/\alpha}\leq 1$. Hence, there is an integer $k_{0}\geq 1$ such that $\dist(t,H_{n})\geq |X_{n}|^{1/\alpha}$ for any $n\geq 1$ with $|X_{n}|\leq 2^{-\lfloor\alpha k_{0}\rfloor}$. Besides, note that $\alpha<1/\beta_{\nu}$, owing to~(\ref{eq:Alphanutleqbetanu}).

%%%%%%%%%%%%%%%%%%%%%%%%%%%%%%%%%%%%%%%%%%%%%%%%%%%%%%%%%%%%%%%%%%%%%%%%%%%%%%
%%%%%%%%%%%%%%%%%%%%%%%%%%%%%%%%%%%%%%%%%%%%%%%%%%%%%%%%%%%%%%%%%%%%%%%%%%%%%%
\subsubsection{Reduction to the study of the component with small jumps}
%%%%%%%%%%%%%%%%%%%%%%%%%%%%%%%%%%%%%%%%%%%%%%%%%%%%%%%%%%%%%%%%%%%%%%%%%%%%%%
%%%%%%%%%%%%%%%%%%%%%%%%%%%%%%%%%%%%%%%%%%%%%%%%%%%%%%%%%%%%%%%%%%%%%%%%%%%%%%

Thanks to~(\ref{eq:Lnu0}) and~(\ref{eq:Lnujalt}), the value at $t$ of the component of the field $L_{\nu}$ that corresponds to the jumps of size larger than $2^{-\lfloor\alpha k_{0}\rfloor}$ may be written as
\[
L^{0,\lfloor\alpha k_{0}\rfloor}_{\nu}(t)=\sum_{n=1}^{\infty}X_{n}\ind_{\{(\Rho_{n},S_{n})\in V_{t},\,|X_{n}|>2^{-\lfloor\alpha k_{0}\rfloor}\}}-\int_{s\in\S^{d-1}\atop x\in(2^{-\lfloor\alpha k_{0}\rfloor},1]} x\inpr{s}{t} \,\nu(\dd s,\dd x),
\]
where $V_{t}$ is defined by~(\ref{eq:defVt}). Moreover, for any $\delta\in(0,A-\|t\|)$, let $V_{t,\delta}$ denote the complement of $\bigcap_{\|\tau\|\leq\delta} V_{t+\tau}$ in $\bigcup_{\|\tau\|\leq\delta} V_{t+\tau}$. If $n\geq 1$ satisfies $(\Rho_{n},S_{n})\in V_{t,\delta}$, then $\Rho_{n}<A$. Thus, the set $\Ncal_{\delta}$ of all $n\geq 1$ such that $(\Rho_{n},S_{n})\in V_{t,\delta}$ and $|X_{n}|>2^{-\lfloor\alpha k_{0}\rfloor}$ is finite. Moreover, given that $t\not\in J_{\nu}$, it is clear that $\bigcap_{\delta>0}\downarrow\Ncal_{\delta}=\emptyset$, so that $\Ncal_{\delta}=\emptyset$ for $\delta$ small enough. For such a $\delta$ and for $\|\tau\|\leq\delta$, we have $V_{t+\tau}\Delta V_{t}\subseteq V_{t,\delta}$. Hence, no integer $n\geq 1$ can satisfy both $(\Rho_{n},S_{n})\in V_{t}\Delta V_{t+\tau}$ and $|X_{n}|>2^{-\lfloor\alpha k_{0}\rfloor}$, so
\[
L^{0,\lfloor\alpha k_{0}\rfloor}_{\nu}(t+\tau)-L^{0,\lfloor\alpha k_{0}\rfloor}_{\nu}(t)=-\int_{s\in\S^{d-1}\atop x\in(2^{-\lfloor\alpha k_{0}\rfloor},1]} x\inpr{s}{\tau} \,\nu(\dd s,\dd x).
\]
Hence, $L^{0,\lfloor\alpha k_{0}\rfloor}_{\nu}$ coincides with an affine form near $t$, and the H\"older exponent at $t$ of $L_{\nu}$ is equal to that of the component with jumps of size at most $2^{-\lfloor\alpha k_{0}\rfloor}$.

%%%%%%%%%%%%%%%%%%%%%%%%%%%%%%%%%%%%%%%%%%%%%%%%%%%%%%%%%%%%%%%%%%%%%%%%%%%%%%
%%%%%%%%%%%%%%%%%%%%%%%%%%%%%%%%%%%%%%%%%%%%%%%%%%%%%%%%%%%%%%%%%%%%%%%%%%%%%%
\subsubsection{Study of the component with small jumps}
%%%%%%%%%%%%%%%%%%%%%%%%%%%%%%%%%%%%%%%%%%%%%%%%%%%%%%%%%%%%%%%%%%%%%%%%%%%%%%
%%%%%%%%%%%%%%%%%%%%%%%%%%%%%%%%%%%%%%%%%%%%%%%%%%%%%%%%%%%%%%%%%%%%%%%%%%%%%%

In order to study the H\"older exponent of the component with jumps of size at most $2^{-\lfloor\alpha k_{0}\rfloor}$, let us consider a vector $\tau\in\R^{d}$ such that $2^{-(k+1)}<\|\tau\|\leq 2^{-k}$ for some integer $k\geq k_{0}$. First,
\[
|L^{\lfloor\alpha k\rfloor+1,\infty}_{\nu}(t+\tau)-L^{\lfloor\alpha k\rfloor+1,\infty}_{\nu}(t)|\leq \Zeta_{\nu}(A)\sum_{j=\lfloor\alpha k\rfloor+1}^{\infty}2^{-j}k(j+2^{-k/2}(j\nu_{j})^{1/2}),
\]
in view of Lemma~\ref{lem:oscLnuj}. Furthermore, if $\beta_{\nu}<2$, then there is a real $\gamma\in(\beta_{\nu},(1/\alpha)\wedge 2)$ and~(\ref{eq:defbetanu}) implies that $c_{\nu,\gamma}=\sum_{j\geq 1} 2^{-\gamma j}j\nu_{j}$ is finite. Hence, we have
\[
k\,2^{-k/2}\sum_{j=\lfloor\alpha k\rfloor+1}^{\infty} 2^{-j}(j\nu_{j})^{1/2} \leq k\,2^{-k/2}c_{\nu,\gamma}^{1/2}\sum_{j=\lfloor\alpha k\rfloor+1}^{\infty} 2^{(\gamma/2-1)j}\leq\frac{c_{\nu,\gamma}^{1/2}}{1-2^{\gamma/2-1}}k\,2^{-\alpha k}.
\]
If $\beta_{\nu}=2$, we observe that $\alpha<1/2$, so that the left-hand side above is at most $\chi_{\nu}k\,2^{-\alpha k}$, where $\chi_{\nu}$ is finite and defined by~(\ref{eq:defchinu}). The upshot is that there exists a deterministic real $D_{\nu,\alpha}>0$ that depends on $\nu$ and $\alpha$ only such that
\begin{equation}\label{eq:smalljumps}
|L^{\lfloor\alpha k\rfloor+1,\infty}_{\nu}(t+\tau)-L^{\lfloor\alpha k\rfloor+1,\infty}_{\nu}(t)|\leq \Zeta_{\nu}(A) D_{\nu,\alpha} k^{2}2^{-\alpha k}.
\end{equation}
Second, no integer $n\geq 1$ can verify at the same time $2^{-\lfloor\alpha k\rfloor}<|X_{n}|\leq 2^{-\lfloor\alpha k_{0}\rfloor}$ and $(\Rho_{n},S_{n})\in V_{t}\Delta V_{t+\tau}$ (otherwise, $|\Rho_{n}-\inpr{S_{n}}{t}|$ would be at most $2^{-k}$ and at least $|X_{n}|^{1/\alpha}>2^{-k}$ simultaneously, which is impossible). Along with~(\ref{eq:Lnujalt}), this yields
\[
|L^{\lfloor\alpha k_{0}\rfloor+1,\lfloor\alpha k\rfloor}_{\nu}(t+\tau)-L^{\lfloor\alpha k_{0}\rfloor+1,\lfloor\alpha k\rfloor}_{\nu}(t)|\leq\|\tau\|\int_{s\in\S^{d-1}\atop x\in(2^{-\lfloor\alpha k\rfloor},2^{-\lfloor\alpha k_{0}\rfloor}]} x\,\nu(\dd s,\dd x).
\]

%%%%%%%%%%%%%%%%%%%%%%%%%%%%%%%%%%%%%%%%%%%%%%%%%%%%%%%%%%%%%%%%%%%%%%%%%%%%%%
%%%%%%%%%%%%%%%%%%%%%%%%%%%%%%%%%%%%%%%%%%%%%%%%%%%%%%%%%%%%%%%%%%%%%%%%%%%%%%
\subsubsection{End of the proof for $\beta_{\nu}\geq 1$}
%%%%%%%%%%%%%%%%%%%%%%%%%%%%%%%%%%%%%%%%%%%%%%%%%%%%%%%%%%%%%%%%%%%%%%%%%%%%%%
%%%%%%%%%%%%%%%%%%%%%%%%%%%%%%%%%%%%%%%%%%%%%%%%%%%%%%%%%%%%%%%%%%%%%%%%%%%%%%

There is a real $\gamma\in [1,2]$ such that $\gamma<1/\alpha$ and the integral $I_{\gamma}$ of $(s,x)\mapsto x^\gamma$ over $\S^{d-1}\times(0,1]$ with respect to $\nu$ is finite. Indeed, one may choose $\gamma=\beta_{\nu}$ if $\beta_{\nu}=2$ and $\gamma>\beta_{\nu}$ sufficiently small otherwise. Hence, using both~(\ref{eq:smalljumps}) and the above bound, we infer that
\begin{align*}
|L^{\lfloor\alpha k_{0}\rfloor+1,\infty}_{\nu}(t+\tau)-L^{\lfloor\alpha k_{0}\rfloor+1,\infty}_{\nu}(t)|&\leq \|\tau\|I_{\gamma}\,2^{(\gamma-1)\lfloor\alpha k\rfloor}+\Zeta_{\nu}(A)D_{\nu,\alpha}k^{2}2^{-\alpha k}\\
&\leq \|\tau\|^{\alpha}\left(I_{\gamma}+2^{\alpha}\Zeta_{\nu}(A)D_{\nu,\alpha}(\log_{2}\|\tau\|)^{2}\right).
\end{align*}
It follows that $\alpha_{L_{\nu}}(t)\geq\alpha$. To conclude, it remains to let $\alpha\uparrow\Alpha_{\nu}(t)$.

%%%%%%%%%%%%%%%%%%%%%%%%%%%%%%%%%%%%%%%%%%%%%%%%%%%%%%%%%%%%%%%%%%%%%%%%%%%%%%
%%%%%%%%%%%%%%%%%%%%%%%%%%%%%%%%%%%%%%%%%%%%%%%%%%%%%%%%%%%%%%%%%%%%%%%%%%%%%%
\subsubsection{End of the proof for $\beta_{\nu}<1$}
%%%%%%%%%%%%%%%%%%%%%%%%%%%%%%%%%%%%%%%%%%%%%%%%%%%%%%%%%%%%%%%%%%%%%%%%%%%%%%
%%%%%%%%%%%%%%%%%%%%%%%%%%%%%%%%%%%%%%%%%%%%%%%%%%%%%%%%%%%%%%%%%%%%%%%%%%%%%%

Here, thanks to~(\ref{eq:Lnujalt}), we have
\[
L^{\lfloor\alpha k_{0}\rfloor+1,\infty}_{\nu}(t)=\sum_{n=1}^{\infty}X_{n}\ind_{\{(\Rho_{n},S_{n})\in V_{t},\,|X_{n}|\leq 2^{-\lfloor\alpha k_{0}\rfloor}\}}-\int_{s\in\S^{d-1}\atop x\in(0,2^{-\lfloor\alpha k_{0}\rfloor}]} x\inpr{s}{t} \,\nu(\dd s,\dd x).
\]
The second term is a linear form, so we just need to study the increments of the first term. To this end, observe that for an arbitrary real $\gamma\in(\beta_{\nu},(1/\alpha)\wedge 1)$,
\begin{align*}
&\left| \sum_{n=1}^{\infty}X_{n}\ind_{\{(\Rho_{n},S_{n})\in V_{t+\tau},\,|X_{n}|\leq 2^{-\lfloor\alpha k_{0}\rfloor}\}} - \sum_{n=1}^{\infty}X_{n}\ind_{\{(\Rho_{n},S_{n})\in V_{t},\,|X_{n}|\leq 2^{-\lfloor\alpha k_{0}\rfloor}\}} \right| \\
\leq & | L^{\lfloor\alpha k\rfloor+1,\infty}_{\nu}(t+\tau) - L^{\lfloor\alpha k\rfloor+1,\infty}_{\nu}(t) |+\|\tau\|\int_{s\in\S^{d-1}\atop x\in(0,2^{-\lfloor\alpha k\rfloor}]} x \,\nu(\dd s,\dd x) \\
\leq & \Zeta_{\nu}(A)D_{\nu,\alpha}k^{2}2^{-\alpha k}+\|\tau\|I_{\gamma}\,2^{-(1-\gamma)\lfloor\alpha k\rfloor} \\
\leq & \|\tau\|^{\alpha}\bigl(2^{(1-\gamma)(1+\alpha)}I_{\gamma}+2^{\alpha}\Zeta_{\nu}(A)D_{\nu,\alpha}(\log_{2}\|\tau\|)^{2}\bigr),
\end{align*}
because of~(\ref{eq:smalljumps}) and the fact that no integer $n\geq 1$ can satisfy $(\Rho_{n},S_{n})\in V_{t}\Delta V_{t+\tau}$ and $2^{-\lfloor\alpha k\rfloor}<|X_{n}|\leq 2^{-\lfloor\alpha k_{0}\rfloor}$ simultaneously. We deduce that $\alpha_{L_{\nu}}(t)\geq\alpha$, and conclude by letting $\alpha\uparrow\Alpha_{\nu}(t)$.

%%%%%%%%%%%%%%%%%%%%%%%%%%%%%%%%%%%%%%%%%%%%%%%%%%%%%%%%%%%%%%%%%%%%%%%%%%%%%%
%%%%%%%%%%%%%%%%%%%%%%%%%%%%%%%%%%%%%%%%%%%%%%%%%%%%%%%%%%%%%%%%%%%%%%%%%%%%%%
%%%%%%%%%%%%%%%%%%%%%%%%%%%%%%%%%%%%%%%%%%%%%%%%%%%%%%%%%%%%%%%%%%%%%%%%%%%%%%
\section{Approximation by Poisson hyperplanes: covering the whole space}\label{sec:proofprpKnualphaRd}
%%%%%%%%%%%%%%%%%%%%%%%%%%%%%%%%%%%%%%%%%%%%%%%%%%%%%%%%%%%%%%%%%%%%%%%%%%%%%%
%%%%%%%%%%%%%%%%%%%%%%%%%%%%%%%%%%%%%%%%%%%%%%%%%%%%%%%%%%%%%%%%%%%%%%%%%%%%%%
%%%%%%%%%%%%%%%%%%%%%%%%%%%%%%%%%%%%%%%%%%%%%%%%%%%%%%%%%%%%%%%%%%%%%%%%%%%%%%

We now prove Proposition~\ref{prp:KnualphaRd}. To begin with, given $\alpha>1/\beta_{\nu}$, let $\nu_{\alpha}$ denote the image under $(s,x)\mapsto(s,|x|^{1/\alpha})$ of the restriction to $\S^{d-1}\times([-1,1]\setminus\{0\})$ of the measure $\nu$. Then, for an arbitrary orthonormal basis $(e_{1},\ldots,e_{d})$ of $\R^{d}$ and for any $s\in\S^{d-1}$, there necessarily exists an integer $i\in\{1,\ldots,d\}$ such that $\inpr{s}{e_{i}}\neq 0$. Together with~(\ref{eq:defbetanu}), this shows that for some $i$ and some $\eps>0$,
\begin{equation}\label{eq:condcoverall}
\int_{s\in\S^{d-1} \atop r\in (0,1]}\ind_{\{\inpr{s}{e_{i}}\neq 0\}}r^{1+\eps} \,\nu_{\alpha}(\dd s,\dd r)=\infty.
\end{equation}

Now, for any integers $A,j\geq 1$, let $U_{A,j}=\opball{0}{A}\cap(2^{-j}/\sqrt{d})\Z^{d}$. Moreover, for $j>j_{0}\geq 1$, let $\event_{A,j_{0},j}$ denote the event consisting in the existence of a point $u\in U_{A,j}$ satisfying $\dist(u,H_{n})\geq |X_{n}|^{1/\alpha}-2^{-j}$ for any integer $n\geq 1$ with $\inpr{S_{n}}{e_{i}}\neq 0$ and $2^{-j}<|X_{n}|^{1/\alpha}\leq 2^{-j_{0}}$. Then, we have
\begin{equation}\label{eq:KnualphaneqRdincl}
\{ K_{\nu}(\alpha)\neq\R^{d} \} \subseteq \bigcup_{A=1}^{\infty}\uparrow\bigcup_{j_{0}=1}^{\infty}\uparrow \bigcap_{j=j_{0}+1}^{\infty} \event_{A,j_{0},j}.
\end{equation}

The event $\event_{A,j_{0},j}$ occurs with probability at most $\sum_{u\in U_{A,j}} \ee^{-I_{j_{0},j,u}}$, where
\[
I_{j_{0},j,u}=\int_{(\rho,s)\in\H_{d} \atop r\in(0,1]} f_{j_{0},j,u}(\rho,s,r) \,\dd\rho\,\nu_{\alpha}(\dd s,\dd r)
\]
and $f_{j_{0},j,u}(\rho,s,r)$ is equal to one when $|\rho-\inpr{s}{u}|<r-2^{-j}$, $\inpr{s}{e_{i}}\neq 0$ and $2^{-j}<r\leq 2^{-j_{0}}$, and is equal to zero otherwise. Using the symmetry of $\nu$, we infer that $I_{j_{0},j,u}$ is equal to
\[
\int_{\inpr{s}{e_{i}}>0 \atop r\in(2^{-j},2^{-j_{0}}]} \int_{\rho\in\R} \ind_{\{|\rho-\inpr{s}{u}|<r-2^{-j}\}}\,\dd\rho\,\nu_{\alpha}(\dd s,\dd r)
=\int_{\inpr{s}{e_{i}}\neq 0 \atop r\in(2^{-j},2^{-j_{0}}]} (r-2^{-j})\,\nu_{\alpha}(\dd s,\dd r).
\]
Let $I'_{j_{0},j}=\nu_{\alpha}(\{s\in\S^{d-1}\:|\:\inpr{s}{e_{i}}\neq 0\}\times(2^{-j},2^{-j_{0}}])$. Then, due to Fubini's theorem,
\[
I_{j_{0},j,u}=\int_{w\in(2^{-j},2^{-j_{0}})}\int_{\inpr{s}{e_{i}}\neq 0 \atop r\in(w,2^{-j_{0}}]} \nu_{\alpha}(\dd s,\dd r)\,\dd w\geq 2^{-j}I'_{j_{0},j-1}.
\]
Therefore, given that $U_{A,j}$ has cardinality at most $(2^{j+2}A\sqrt{d})^{d}$, we get
\begin{equation}\label{eq:boundprobEAj0j}
\prob(\event_{A,j_{0},j}) \leq (4A\sqrt{d})^{d} \exp(j d \log 2-2^{-j}I'_{j_{0},j-1}).
\end{equation}

Finally, employing Fubini's theorem again, we also have
\begin{align*}
\int_{\inpr{s}{e_{i}}\neq 0 \atop r\in (0,2^{-j_{0}}]} r^{1+\eps} \,\nu_{\alpha}(\dd s,\dd r) &= (1+\eps)\int_{w\in(0,2^{-j_{0}})}w^{\eps}\int_{\inpr{s}{e_{i}}\neq 0 \atop r\in(w,2^{-j_{0}}]}\nu_{\alpha}(\dd s,\dd r)\,\dd w \\
&\leq (1+\eps)\sum_{j=j_{0}}^{\infty} I'_{j_{0},j+1}\int_{w\in(2^{-(j+1)},2^{-j})}w^{\eps} \,\dd w.
\end{align*}
Together with~(\ref{eq:condcoverall}), this ensures that $I'_{j_{0},j-1}>2^{(1+\eps)j}/j^{2}$ for infinitely many $j>j_{0}$. We conclude with the help of~(\ref{eq:KnualphaneqRdincl}) and~(\ref{eq:boundprobEAj0j}).

%%%%%%%%%%%%%%%%%%%%%%%%%%%%%%%%%%%%%%%%%%%%%%%%%%%%%%%%%%%%%%%%%%%%%%%%%%%%%%
%%%%%%%%%%%%%%%%%%%%%%%%%%%%%%%%%%%%%%%%%%%%%%%%%%%%%%%%%%%%%%%%%%%%%%%%%%%%%%
%%%%%%%%%%%%%%%%%%%%%%%%%%%%%%%%%%%%%%%%%%%%%%%%%%%%%%%%%%%%%%%%%%%%%%%%%%%%%%
\section{Approximation by Poisson hyperplanes: size and large intersection properties}\label{sec:proofthmsliKnualpha}
%%%%%%%%%%%%%%%%%%%%%%%%%%%%%%%%%%%%%%%%%%%%%%%%%%%%%%%%%%%%%%%%%%%%%%%%%%%%%%
%%%%%%%%%%%%%%%%%%%%%%%%%%%%%%%%%%%%%%%%%%%%%%%%%%%%%%%%%%%%%%%%%%%%%%%%%%%%%%
%%%%%%%%%%%%%%%%%%%%%%%%%%%%%%%%%%%%%%%%%%%%%%%%%%%%%%%%%%%%%%%%%%%%%%%%%%%%%%

This section is devoted to the proof of Theorem~\ref{thm:sliKnualpha}. We begin by establishing a series of preliminary lemmas. Then, we deal with the case where $h_{\nu}(g)>\alpha$. We finally end the proof with the case where $h_{\nu}(g)<\alpha$. Note that Theorem~\ref{thm:sliKnualpha} clearly holds when $\nu$ has finite total mass. Indeed, in this case, it is easy to check that $h_{\nu}(g)=\infty$ for any $g\in\gauge_{1}$, while $K_{\nu}(\alpha)=\emptyset$ with probability one for every $\alpha>0$ (given $A\geq 1$, there are almost surely finitely many $n\geq 1$ such that $H_{n}\cap\clball{0}{A}\neq\emptyset$ or, equivalently, such that $\Rho_{n}\leq A$). Therefore, we may assume throughout the section that $\nu$ has infinite total mass.

%%%%%%%%%%%%%%%%%%%%%%%%%%%%%%%%%%%%%%%%%%%%%%%%%%%%%%%%%%%%%%%%%%%%%%%%%%%%%%
%%%%%%%%%%%%%%%%%%%%%%%%%%%%%%%%%%%%%%%%%%%%%%%%%%%%%%%%%%%%%%%%%%%%%%%%%%%%%%
\subsection{A Bernstein-type inequality}
%%%%%%%%%%%%%%%%%%%%%%%%%%%%%%%%%%%%%%%%%%%%%%%%%%%%%%%%%%%%%%%%%%%%%%%%%%%%%%
%%%%%%%%%%%%%%%%%%%%%%%%%%%%%%%%%%%%%%%%%%%%%%%%%%%%%%%%%%%%%%%%%%%%%%%%%%%%%%

The first preliminary lemma yields an analog of Bernstein's inequality for integrals with respect to a compensated Poisson random measure, and is a direct consequence of~\cite[Corollary~5.1]{Houdre:2002kx}, see also~\cite[Proposition~7]{Reynaud-Bouret:2003kx}. It comes into play in the proof of Lemma~\ref{lem:oscLnuj} too, see Section~\ref{sec:prooflemoscLnuj}.

\begin{lem}\label{lem:Bernstein}
Let $(E,\mathcal{E})$ be a measurable space endowed with a finite nonnegative measure $\mu$, let $\Mu$ be a Poisson random measure on $E$ with intensity $\mu$ and let $\Mu^{*}=\Mu-\mu$. Then, for any real-valued measurable function $f$ defined on $E$ such that $S=\sup_{E}|f|$ and $V=\int_{E} f^{2} \, \dd\mu$ are both positive and finite, we have
\[
\forall \xi>0 \qquad \prob\left(\left|\int_{E} f \, \dd\Mu^{*}\right|\geq\xi\right)\leq 2\exp\left(-\frac{3\xi^2}{2S\xi+6V}\right).
\]
\end{lem}

%%%%%%%%%%%%%%%%%%%%%%%%%%%%%%%%%%%%%%%%%%%%%%%%%%%%%%%%%%%%%%%%%%%%%%%%%%%%%%
%%%%%%%%%%%%%%%%%%%%%%%%%%%%%%%%%%%%%%%%%%%%%%%%%%%%%%%%%%%%%%%%%%%%%%%%%%%%%%
\subsection{A law of large numbers for Poisson measures}
%%%%%%%%%%%%%%%%%%%%%%%%%%%%%%%%%%%%%%%%%%%%%%%%%%%%%%%%%%%%%%%%%%%%%%%%%%%%%%
%%%%%%%%%%%%%%%%%%%%%%%%%%%%%%%%%%%%%%%%%%%%%%%%%%%%%%%%%%%%%%%%%%%%%%%%%%%%%%

The second preliminary lemma concerns the behavior at zero of a Poisson random measure on $(0,1]$ and is reminiscent of the strong law of large numbers for the homogeneous Poisson process on $(0,\infty)$, see~\cite{Kingman:1993gf}. Let $\Pcal$ be the set of all nonnegative Borel measures $\pi$ on $(0,1]$ such that $\pi((0,1])=\infty$ and $\pi([\eps,1])<\infty$ for any $\eps>0$.

\begin{lem}\label{lem:asymptPi}
For any Poisson random measure $\Pi$ on $(0,1]$ with intensity $\pi\in\Pcal$,
\[
\as \qquad \Pi([w,1]) \underset{w\to 0}{\sim} \pi([w,1]).
\]
\end{lem}

\begin{proof}
Let $\Acal_{0}$ be the (countable) set of all $r\in (0,1]$ such that $\pi(\{r\})\geq 1$ and let $\Acal_{1}=(0,1]\setminus\Acal_{0}$. Then, for all $\ell\in\{0,1\}$ and $w\in(0,1]$, let $\Phi_{\ell}(w)=\Pi(\Acal_{\ell}\cap[w,1])$ and $\ph_{\ell}(w)=\pi(\Acal_{\ell}\cap[w,1])$. It is easy to see that the proof reduces to showing that for any $\ell\in\{0,1\}$ such that $\pi(\Acal_{\ell})=\infty$,
\begin{equation}\label{eq:asymptPhiell}
\as \qquad \Phi_{\ell}(w) \underset{w\to 0}{\sim} \ph_{\ell}(w).
\end{equation}

To this purpose, let us begin by observing that for any $\xi>0$ and any $w>0$ small enough to ensure that $\ph_{\ell}(w)>0$,
\begin{equation}\label{eq:BernsteinPhi}
\prob\left( \left|\frac{\Phi_{\ell}(w)}{\ph_{\ell}(w)}-1\right|\geq\xi \right) \leq 2\exp\left(-\frac{3\xi^2}{2\xi+6}\ph_{\ell}(w)\right),
\end{equation}
as a consequence of Lemma~\ref{lem:Bernstein}. Now, if $\pi(\Acal_{0})=\infty$, there exists a decreasing sequence $(a_{n})_{n\geq 1}$ of positive reals that converges to zero and whose terms form the set $\Acal_{0}$. The previous inequality then implies that for all integers $m,n\geq 1$,
\begin{equation}\label{eq:BernsteinPhi0}
\prob\left( \left|\frac{\Phi_{0}(a_{n})}{\ph_{0}(a_{n})}-1\right|\geq\frac{1}{m} \right)\leq 2\exp\left(-\frac{3n}{(6m+2)m}\right),
\end{equation}
because $\ph_{0}(a_{n})=\pi(\{a_{1},\ldots,a_{n}\})\geq n$. Summing over $n\geq 1$ and making use of the Borel-Cantelli lemma, we infer that for any integer $m\geq 1$, with probability one, $\Phi_{0}(a_{n})/\ph_{0}(a_{n})$ is between $1-1/m$ and $1+1/m$, for $n$ large enough. The same property holds for $\Phi_{0}(w)/\ph_{0}(w)$ with $w>0$ small enough, due to the fact that $\Phi_{0}(w)=\Phi_{0}(a_{n(w)})$ and $\ph_{0}(w)=\ph_{0}(a_{n(w)})$, where $n(w)$ is the number of integers $n\geq 1$ such that $a_{n}\geq w$. As a consequence,
\[
\as \quad \forall m\geq 1 \qquad \limsup_{w\to 0}\left|\frac{\Phi_{0}(w)}{\ph_{0}(w)}-1\right|\leq\frac{1}{m},
\]
and we get~(\ref{eq:asymptPhiell}) for $\ell=0$ by letting $m\to\infty$.

If $\pi(\Acal_{1})=\infty$, it is possible to consider, for each integer $n\geq 1$,
\[
w_{n}=\sup\left\{w>0 \:\bigl|\: \ph_{1}(w)\geq n\right\}.
\]
The reals $w_{n}$ are positive, satisfy $\ph_{1}(w_{n})\geq n$ and form a nonincreasing sequence that converges to zero. Moreover,~(\ref{eq:BernsteinPhi}) ensures that for any $m,n\geq 1$, the bound~(\ref{eq:BernsteinPhi0}) holds with $a_{n}$ replaced by $w_{n}$, and $\Phi_{0}$ and $\ph_{0}$ replaced by $\Phi_{1}$ and $\ph_{1}$, respectively. Summing over $n\geq 1$ and using the Borel-Cantelli lemma again, it follows that with probability one, for $n$ large enough, $\Phi_{1}(w_{n})/\ph_{1}(w_{n})$ is between $1-1/m$ and $1+1/m$. In addition, we have
\[
n\leq\ph_{1}(w_{n})=\pi(\Acal_{1}\cap\{w_{n}\})+\lim_{w\downarrow w_{n}}\uparrow\ph_{1}(w)\leq 1+n,
\]
by definition of $\Acal_{1}$ and $w_{n}$. Therefore, since $\Phi_{1}$ and $\ph_{1}$ are nonincreasing, we infer that with probability one, for $n$ large enough and for $w\in [w_{n+1},w_{n}]$,
\[
\frac{n}{n+2}\left(1-\frac{1}{m}\right)\leq\frac{\Phi_{1}(w)}{\ph_{1}(w)}\leq\frac{n+2}{n}\left(1+\frac{1}{m}\right),
\]
and~(\ref{eq:asymptPhiell}) with $\ell=1$ follows in a straightforward manner.
\end{proof}

%%%%%%%%%%%%%%%%%%%%%%%%%%%%%%%%%%%%%%%%%%%%%%%%%%%%%%%%%%%%%%%%%%%%%%%%%%%%%%
%%%%%%%%%%%%%%%%%%%%%%%%%%%%%%%%%%%%%%%%%%%%%%%%%%%%%%%%%%%%%%%%%%%%%%%%%%%%%%
\subsection{Integrability with respect to a Poisson measure}
%%%%%%%%%%%%%%%%%%%%%%%%%%%%%%%%%%%%%%%%%%%%%%%%%%%%%%%%%%%%%%%%%%%%%%%%%%%%%%
%%%%%%%%%%%%%%%%%%%%%%%%%%%%%%%%%%%%%%%%%%%%%%%%%%%%%%%%%%%%%%%%%%%%%%%%%%%%%%

Our third preliminary lemma is a direct consequence of the second one and deals with the integrability of a gauge function with respect to a Poisson random measure.

\begin{lem}\label{lem:convPi}
Let $\Pi$ be a Poisson random measure on $(0,1]$ with intensity $\pi\in\Pcal$. Then, with probability one, for any gauge function $g\in\gauge_{1}$,
\[
\int_{r\in (0,1]} g(r)\,\Pi(\dd r)=\infty \quad\Longleftrightarrow\quad \int_{r\in (0,1]} g(r)\,\pi(\dd r)=\infty.
\]
\end{lem}

\begin{proof}
Let us assume that the almost sure event on which the statement of Lemma~\ref{lem:asymptPi} holds occurs. Then, let $g\in\gauge_{1}$ and let $\gamma$ denote the Lebesgue-Stieltjes measure associated with $g$. Given that $\Pi([w,1])$ is equivalent to $\pi([w,1])$ as $w\to 0$, we have
\[
\int_{w\in(0,1]} \Pi([w,1])\,\gamma(\dd w)=\infty \qquad\Longleftrightarrow\qquad \int_{w\in(0,1]} \pi([w,1])\,\gamma(\dd w)=\infty.
\]
We conclude by remarking that, owing to the Fubini-Tonelli theorem, the above integrals are equal to those appearing in the statement of the lemma.
\end{proof}

%%%%%%%%%%%%%%%%%%%%%%%%%%%%%%%%%%%%%%%%%%%%%%%%%%%%%%%%%%%%%%%%%%%%%%%%%%%%%%
%%%%%%%%%%%%%%%%%%%%%%%%%%%%%%%%%%%%%%%%%%%%%%%%%%%%%%%%%%%%%%%%%%%%%%%%%%%%%%
\subsection{Approximation by homogeneously distributed hyperplanes}
%%%%%%%%%%%%%%%%%%%%%%%%%%%%%%%%%%%%%%%%%%%%%%%%%%%%%%%%%%%%%%%%%%%%%%%%%%%%%%
%%%%%%%%%%%%%%%%%%%%%%%%%%%%%%%%%%%%%%%%%%%%%%%%%%%%%%%%%%%%%%%%%%%%%%%%%%%%%%

The next lemma is a general result on Diophantine approximation by hyperplanes, under the assumption that the hyperplanes are homogeneously distributed, in a specific sense that we now introduce. Given $\sfrak\in\S^{d-1}$, let $\Hyp_{\sfrak}$ be the set of all hyperplanes $h$ represented by a pair $(\rho,s)\in\H_{d}$ with $\inpr{s}{\sfrak}\neq 0$. Such a hyperplane $h$ is not parallel to $\sfrak$. So, for any $t\in\R^{d}$, the line $t+\R\sfrak$ meets $h$ at a single point $t+\xi^{\sfrak}_{t}(h)\,\sfrak$, where
\[
\xi^{\sfrak}_{t}(h)=\frac{\rho-\inpr{s}{t}}{\inpr{s}{\sfrak}}.
\]

Now, let $\Hfrak=(\Hfrak_{n})_{n\geq 1}$ be a sequence in $\Hyp_{\sfrak}$ and let $W$ be a nonempty open subset of $\R^{d}$. Given $t\in W$, $\delta\in (0,1)$ and $j\geq 0$, let us consider the first $\lfloor 2^{j}/\delta\rfloor$ hyperplanes $\Hfrak_{n}$ and focus on those which intersect the line $t+\R\sfrak$ at a single point located at a distance less than $\delta$ from $t$ or, equivalently, those for which $|\xi^{\sfrak}_{t}(\Hfrak_{n})|<\delta$. The fact that the sequence $\Hfrak$ is homogeneously distributed in $W$ basically means that the resulting intersection points are dispersed in a regular manner around almost every point $t\in W$ in the direction $\sfrak$, in the sense that the set
\[
\Qrm^{\sfrak}_{t,\delta,j}(\Hfrak)=\left\{ q\in\{0,\ldots,2^{j}-1\} \:\biggl|\: q=\left\lfloor \frac{2^{j}}{\delta}|\xi^{\sfrak}_{t}(\Hfrak_{n})|\right\rfloor \mbox{ for some } n\leq \frac{2^j}{\delta} \right\}
\]
has cardinality of the order of $2^j$, asymptotically. Accordingly, we let $\Hom_{\sfrak}(W)$ be the set of sequences $\Hfrak$ in $\Hyp_{\sfrak}$ such that for $\leb^{d}$-almost every $t\in W$,
\[
\limsup_{\delta\to 0}\underline{\Qrm}^{\sfrak}_{t,\delta}(\Hfrak)>0 \qquad\mbox{with}\qquad \underline{\Qrm}^{\sfrak}_{t,\delta}(\Hfrak)=\liminf_{j\to\infty} \frac{1}{2^{j}}\#\Qrm^{\sfrak}_{t,\delta,j}(\Hfrak),
\]
where $\#$ stands for cardinality, and we say that such a sequence is homogeneously distributed in $W$. Here, $\leb^{d}$ denotes the Lebesgue measure in $\R^{d}$.

\begin{lem}\label{eq:homdistfull}
Let $\Hfrak=(\Hfrak_{n})_{n\geq 1}$ be a sequence in $\Hom_{\sfrak}(W)$ and let $\Rfrak=(\Rfrak_{n})_{n\geq 1}$ be a nonincreasing sequence of positive reals such that $\sum\nolimits_{n} \Rfrak_{n}=\infty$. Then, the set
\[
\Ffrak(\Hfrak,\Rfrak)=\{ t\in\R^{d} \:|\: \dist(t,\Hfrak_{n})<\Rfrak_{n} \mbox{ for i.m. } n\geq 1 \}
\]
has full Lebesgue measure in $W$.
\end{lem}

\begin{proof}
For any $t\in\R^{d}$, let $\Efrak_{t}(\Hfrak,\Rfrak)$ be the set of $\tau\in\R$ such that $|\tau-\xi^{\sfrak}_{t}(\Hfrak_{n})|<\Rfrak_{n}$ for infinitely many $n\geq 1$. It suffices to prove that
\begin{equation}\label{eq:minlebEtHR}
\forall t\in\R^d \quad \forall\delta\in (0,1) \qquad \leb^{1}(\Efrak_{t}(\Hfrak,\Rfrak)\cap (-\delta,\delta))\geq\frac{\delta}{18}\,\underline{\Qrm}^{\sfrak}_{t,\delta}(\Hfrak)^{2}.
\end{equation}
Indeed, letting $(\R\sfrak)^{\perp}$ be the orthogonal complement of $\R\sfrak$ in $\R^{d}$, we deduce from~(\ref{eq:minlebEtHR}) that for all $b\in (\R\sfrak)^{\perp}$ and $\tau\in\R$,
\[
\forall\delta\in (0,1) \qquad \frac{1}{2\delta}\int_{\tau-\delta}^{\tau+\delta} \ind_{\Efrak_{b}(\Hfrak,\Rfrak)}(v) \,\dd v\geq \left(\frac{\underline{\Qrm}^{\sfrak}_{b+\tau\sfrak,\delta}(\Hfrak)}{6}\right)^{2}.
\]
By virtue of Lebesgue's density theorem~\cite[Corollary~2.14]{Mattila:1995fk}, this implies that for every $b\in (\R\sfrak)^{\perp}$ and $\leb^{1}$-almost every $\tau\in\R$,
\[
\ind_{\Efrak_{b}(\Hfrak,\Rfrak)}(\tau)\geq\left(\frac{1}{6}\limsup_{\delta\to 0}\underline{\Qrm}^{\sfrak}_{b+\tau\sfrak,\delta}(\Hfrak)\right)^{2}.
\]
As the sequence $\Hfrak$ is in $\Hom_{\sfrak}(W)$, we infer that for $\hau^{d-1}$-almost every $b\in (\R\sfrak)^{\perp}$ and $\leb^{1}$-almost every $\tau\in\R$,
\[
b+\tau\sfrak\in W \qquad\Longrightarrow\qquad \ind_{\Efrak_{b}(\Hfrak,\Rfrak)}(\tau)>0 \qquad\Longrightarrow\qquad b+\tau\sfrak\in\Ffrak(\Hfrak,\Rfrak),
\]
and the result follows. It now remains to establish~(\ref{eq:minlebEtHR}).  To this end, we may clearly assume that $\underline{\Qrm}^{\sfrak}_{t,\delta}(\Hfrak)>0$, and then show that for every $\eta>1$,
\begin{equation}\label{eq:minlebEtHReta}
\leb^{1}(\Efrak_{t}(\Hfrak,\Rfrak)\cap (-\delta,\delta))\geq \frac{\delta}{18}\left(\frac{\underline{\Qrm}^{\sfrak}_{t,\delta}(\Hfrak)}{\eta}\right)^2.
\end{equation}
In order to prove~(\ref{eq:minlebEtHReta}), let us observe that, by definition of $\underline{\Qrm}^{\sfrak}_{t,\delta}(\Hfrak)$, there is an integer $\underline{j}\geq 0$ such that for any $j\geq\underline{j}$, there exists a set $\Ncal_{j}\subseteq\{1,\ldots,\lfloor 2^{j}/\delta \rfloor\}$ with:
\begin{itemize}
\item $2^{j-1}\underline{\Qrm}^{\sfrak}_{t,\delta}(\Hfrak)/\eta\leq\#\Ncal_{j}\leq 2^{j}$;
\item $|\xi^{\sfrak}_{t}(\Hfrak_{n})|<\delta$ for any $n\in\Ncal_{j}$;
\item $\bigl||\xi^{\sfrak}_{t}(\Hfrak_{n})|-|\xi^{\sfrak}_{t}(\Hfrak_{n'})|\bigr|\geq\delta 2^{-j}$ for any distinct $n,n'\in\Ncal_{j}$.
\end{itemize}
In addition, for each $n\geq 1$, let $\Rfrak'_{n}=\Rfrak_{n}\wedge(2(n+1))^{-1}$. Then, $(\Rfrak'_{n})_{n\geq 1}$ is nonincreasing and $\sum_{n}\Rfrak'_{n}=\infty$, due to Olivier's theorem. Indeed, owing to that result, the convergence of $\sum_{n}\Rfrak'_{n}$ would imply that $\Rfrak'_{n}={\rm o}(1/n)$ as $n\to\infty$, thus contradicting the divergence of $\sum_{n}\Rfrak_{n}$.

For every $j\geq\underline{j}$, let us consider the set
\[
U_{j}=\bigcup_{n\in\Ncal_{j}} \left(\xi^{\sfrak}_{t}(\Hfrak_{n})-\Rfrak'_{\lfloor 2^{j}/\delta\rfloor},\xi^{\sfrak}_{t}(\Hfrak_{n})+\Rfrak'_{\lfloor 2^{j}/\delta\rfloor}\right).
\]
All the points in the $\limsup$ of the sets $U_{j}$, except maybe those of the form $\xi^{\sfrak}_{t}(\Hfrak_{n})$, belong to both $[-\delta,\delta]$ and $\Efrak_{t}(\Hfrak,\Rfrak)$. Therefore,
\[
\leb^{1}\left(\limsup_{j\to\infty} U_{j}\right)\leq\leb^{1}(\Efrak_{t}(\Hfrak,\Rfrak)\cap (-\delta,\delta)).
\]
Lemma~5 in~\cite{Sprindzuk:1979eu} enables us to obtain an appropriate lower bound on the left-hand side above. To be more specific, this result ensures that
\[
\leb^{1}\left(\limsup_{j\to\infty} U_{j}\right)\geq\limsup_{J\to\infty} \frac{\biggl(\sum\limits_{\underline{j}\leq j\leq J}\leb^{1}(U_{j})\biggr)^{2}}{\sum\limits_{\underline{j}\leq j,j'\leq J} \leb^{1}(U_{j}\cap U_{j'})},
\]
with the proviso that $\sum_{j}\leb^{1}(U_{j})$ diverges. To check this last condition, it is crucial to observe that for each $j\geq\underline{j}$, the intervals forming the set $U_{j}$ are disjoint, so that
\begin{equation}\label{eq:minlebUj}
\leb^{1}(U_{j})=2\Rfrak'_{\lfloor 2^{j}/\delta\rfloor}\#\Ncal_{j}\geq \frac{\underline{\Qrm}^{\sfrak}_{t,\delta}(\Hfrak)}{\eta} 2^{j} \Rfrak'_{\lfloor 2^{j}/\delta\rfloor}.
\end{equation}
In view of the divergence of $\sum_{n}\Rfrak'_{n}$ and the Cauchy condensation test, $\sum_{j} 2^{j} \Rfrak'_{\lfloor 2^{j}/\delta\rfloor}$ diverges too. The divergence of $\sum_{j}\leb^{1}(U_{j})$ then follows from~(\ref{eq:minlebUj}). As a consequence, in order to obtain~(\ref{eq:minlebEtHReta}), it now suffices to show that for $J$ large enough,
\begin{equation}\label{eq:sumlebcapUj}
\sum_{\underline{j}\leq j,j'\leq J} \leb^{1}(U_{j}\cap U_{j'})\leq \frac{18}{\delta}\left(\frac{\eta}{\underline{\Qrm}^{\sfrak}_{t,\delta}(\Hfrak)}\sum_{\underline{j}\leq j\leq J}\leb^{1}(U_{j})\right)^{2}.
\end{equation}
To this end, let us derive an upper bound on the Lebesgue measure of $U_{j}\cap U_{j'}$, where $\underline{j}\leq j<j'$. This set is the union over $n\in\Ncal_{j}$ and $n'\in\Ncal_{j',n}$ of the sets
\[
(\xi^{\sfrak}_{t}(\Hfrak_{n})-\Rfrak'_{\lfloor 2^{j}/\delta\rfloor},\xi^{\sfrak}_{t}(\Hfrak_{n})+\Rfrak'_{\lfloor 2^{j}/\delta\rfloor})\cap(\xi^{\sfrak}_{t}(\Hfrak_{n'})-\Rfrak'_{\lfloor 2^{j'}/\delta\rfloor},\xi^{\sfrak}_{t}(\Hfrak_{n'})+\Rfrak'_{\lfloor 2^{j'}/\delta\rfloor}),
\]
where $\Ncal_{j',n}$ is the collection of all $n'\in\Ncal_{j'}$ such that this last intersection is nonempty. The cardinality of $\Ncal_{j',n}$ is clearly bounded above by the number of integers of the form $\lfloor 2^{j'}|\xi|/\delta\rfloor$ with $|\xi-\xi^{\sfrak}_{t}(\Hfrak_{n})|<2\Rfrak'_{\lfloor 2^{j}/\delta\rfloor}$, which is itself at most $2+2^{j'+2}\Rfrak'_{\lfloor 2^{j}/\delta\rfloor}/\delta$. Along with the fact that $\#\Ncal_{j}\leq 2^{j}$, this yields
\[
\leb^{1}(U_{j}\cap U_{j'})\leq 2^{j} \left(2+\frac{2^{j'+2}}{\delta}\Rfrak'_{\lfloor 2^{j}/\delta\rfloor}\right) \left(2\Rfrak'_{\lfloor 2^{j'}/\delta\rfloor}\right).
\]
As a consequence, for any integer $J>\underline{j}$, the left-hand side of~(\ref{eq:sumlebcapUj}) is at most
\[
2\sum_{\underline{j}\leq j\leq J} 2^{j} \Rfrak'_{\lfloor 2^{j}/\delta\rfloor} + 8\sum_{\underline{j}\leq j<j'\leq J} 2^{j} \Rfrak'_{\lfloor 2^{j'}/\delta\rfloor} + \frac{16}{\delta}\sum_{\underline{j}\leq j<j'\leq J} 2^{j+j'}\Rfrak'_{\lfloor 2^{j}/\delta\rfloor}\Rfrak'_{\lfloor 2^{j'}/\delta\rfloor}.
\]
The third sum is smaller than half the sum over all $j,j'\in\{\underline{j},\ldots,J\}$, and the second sum is smaller than the first one. Thus, the left-hand side of~(\ref{eq:sumlebcapUj}) is at most
\[
10\sum_{\underline{j}\leq j\leq J} 2^{j} \Rfrak'_{\lfloor 2^{j}/\delta\rfloor} + \frac{8}{\delta}\left(\sum_{\underline{j}\leq j\leq J} 2^{j}\Rfrak'_{\lfloor 2^{j}/\delta\rfloor}\right)^{2}\leq\frac{18}{\delta}\left(\sum_{\underline{j}\leq j\leq J} 2^{j}\Rfrak'_{\lfloor 2^{j}/\delta\rfloor}\right)^{2},
\]
where the last bound holds for $J$ large enough, due to the divergence of the series $\sum_{j} 2^{j} \Rfrak'_{\lfloor 2^{j}/\delta\rfloor}$. We conclude using~(\ref{eq:minlebUj}).
\end{proof}

%%%%%%%%%%%%%%%%%%%%%%%%%%%%%%%%%%%%%%%%%%%%%%%%%%%%%%%%%%%%%%%%%%%%%%%%%%%%%%
%%%%%%%%%%%%%%%%%%%%%%%%%%%%%%%%%%%%%%%%%%%%%%%%%%%%%%%%%%%%%%%%%%%%%%%%%%%%%%
\subsection{Poisson hyperplanes are homogeneously distributed}
%%%%%%%%%%%%%%%%%%%%%%%%%%%%%%%%%%%%%%%%%%%%%%%%%%%%%%%%%%%%%%%%%%%%%%%%%%%%%%
%%%%%%%%%%%%%%%%%%%%%%%%%%%%%%%%%%%%%%%%%%%%%%%%%%%%%%%%%%%%%%%%%%%%%%%%%%%%%%

The next lemma asserts essentially that the random hyperplanes arising in the definition of $K_{\nu}(\alpha)$ are homogeneously distributed. To be more specific, we need to introduce some additional notations. Given $\alpha>0$, let $\Nu_{\alpha}$ denote the image under the mapping $(\rho,s,x)\mapsto(\rho,s,|x|^{1/\alpha})$ of the restriction to $(0,\infty)\times\S^{d-1}\times([-1,1]\setminus\{0\})$ of the measure $\Nu$. Then, $\Nu_{\alpha}$ is a Poisson random measure with intensity $\leb^{1}_{+}\otimes\nu_{\alpha}$, where $\nu_{\alpha}$ is defined at the beginning of Section~\ref{sec:proofprpKnualphaRd}.

Given an arbitrary orthonormal basis $(e_{1},\ldots,e_{d})$ of $\R^{d}$, for each $i\in\{1,\ldots,d\}$, let $\Ccal_{i}$ be the set of all $s\in\S^{d-1}$ satisfying $|\inpr{s}{e_{i}}|\geq d^{-1/2}$. Furthermore, let $\Acal_{i,0}$ be the set of all $r\in(0,1]$ such that $\nu_{\alpha}(\Ccal_{i}\times\{r\})\geq 1$, and let $\Acal_{i,1}=(0,1]\setminus\Acal_{i,0}$. Then, for each $\ell\in\{0,1\}$, let $\nu_{\alpha,i,\ell}$ be the restriction of $\nu_{\alpha}$ to $\Ccal_{i}\times\Acal_{i,\ell}$.

Let us assume that $\nu_{\alpha,i,\ell}$ has infinite total mass. Thus, $\nu_{\alpha,i,\ell}$ belongs to $\Pcal$, in view of~(\ref{eq:nuLevy}). Given an integer $A\geq 1$, the restriction $\Nu^{A,i,\ell}_{\alpha}$ of $\Nu_{\alpha}$ to $(0,A)\times\Ccal_{i}\times\Acal_{i,\ell}$ may be written almost surely as
\begin{equation}\label{eq:NuAalphaiellRhoSR}
\Nu^{A,i,\ell}_{\alpha}=\sum_{n=1}^{\infty} \delta_{(\Rho^{A,i,\ell}_{n},S^{A,i,\ell}_{n},R^{A,i,\ell}_{n})},
\end{equation}
for some sequence $(\Rho^{A,i,\ell}_{n},S^{A,i,\ell}_{n},R^{A,i,\ell}_{n})_{n\geq 1}$ in $(0,A)\times\Ccal_{i}\times\Acal_{i,\ell}$. Since $\nu_{\alpha,i,\ell}\in\Pcal$, we see that almost surely, for all $\eps>0$, only finitely many $n\geq 1$ satisfy $R^{A,i,\ell}_{n}\geq\eps$. So, up to a reordering, we may assume that the sequence $(R^{A,i,\ell}_{n})_{n\geq 1}$ is nonincreasing and converges to zero. Last, for any $n\geq 1$, let $H^{A,i,\ell}_{n}$ denote the hyperplane defined in terms of $\Rho^{A,i,\ell}_{n}$ and $S^{A,i,\ell}_{n}$ as in~(\ref{eq:defHn}).

\begin{lem}\label{lem:HAilhomdist}
Almost surely, for any $A\geq 1$, $i\in\{1,\ldots,d\}$ and $\ell\in\{0,1\}$ such that $\nu_{\alpha,i,\ell}$ has infinite total mass, $H^{A,i,\ell}=(H^{A,i,\ell}_{n})_{n\geq 1}$ is in $\Hom_{e_{i}}(\opball{0}{A})$.
\end{lem}

\begin{proof}
Given $A\geq 1$, $i\in\{1,\ldots,d\}$ and $\ell\in\{0,1\}$, it is clear that $H^{A,i,\ell}_{n}\in\Hyp_{e_{i}}$ for all $n\geq 1$. Then, let us assume that $\nu_{\alpha,i,\ell}$ has infinite total mass. In view of Fubini's theorem, it suffices to let $t\in\opball{0}{A}$ and show that with probability one, $\limsup_{\delta\to 0}\underline{\Qrm}^{e_{i}}_{t,\delta}(H^{A,i,\ell})>0$. To this end, given $\delta\in (0,1)$ and $j\geq 1$, suppose that $\#\Qrm^{e_{i}}_{t,\delta,j}(H^{A,i,\ell})\leq\theta 2^{j}$, for some $\theta\in(0,1)$. So, there is a subset $\Qrm$ of $\{0,\ldots,2^{j}-1\}$ with cardinality $\lfloor\theta 2^{j}\rfloor$ which satisfies $\Qrm\supseteq\Qrm^{e_{i}}_{t,\delta,j}(H^{A,i,\ell})$. We have
\begin{equation}\label{eq:xieiHail}
\forall n\leq\left\lfloor\frac{2^{j}}{\delta}\right\rfloor \qquad |\xi^{e_{i}}_{t}(H^{A,i,\ell}_{n})|=\left|\frac{\Rho^{A,i,\ell}_{n}-\inpr{S^{A,i,\ell}_{n}}{t}}{\inpr{S^{A,i,\ell}_{n}}{e_{i}}}\right|\in[\delta,\infty]\cup\bigcup_{q\in\Qrm} \lambda_{q},
\end{equation}
where $\lambda_{q}$ denotes the interval $[q\delta 2^{-j},(q+1)\delta 2^{-j})$.

Let us first assume that $\ell=0$ and derive an upper bound on the probability that~(\ref{eq:xieiHail}) happens. There is a decreasing sequence $(a^{i}_{p})_{p\geq 1}$ of positive reals that converges to zero and whose terms form the set $\Acal_{i,0}$. Then, basic properties of Poisson random measures enable us to write that almost surely,
\[
\Nu^{A,i,0}_{\alpha}=\sum_{p=1}^{\infty}\sum_{m=1}^{N_{p}}\delta_{(\widetilde\Rho^{A,i,0}_{p,m},\widetilde S^{A,i,0}_{p,m},a^{i}_{p})},
\]
where each $N_{p}$ is Poisson distributed with mean $A\,\nu_{\alpha,i,0}(\Ccal_{i}\times\{a^{i}_{p}\})$, each $\widetilde\Rho^{A,i,0}_{p,m}$ is uniformly distributed on $(0,A)$ and each $\widetilde S^{A,i,0}_{p,m}$ is distributed according to the probability measure $\mu_{i,a^{i}_{p}}=\nu_{\alpha,i,0}(\,\cdot\times\{a^{i}_{p}\})/\nu_{\alpha,i,0}(\Ccal_{i}\times\{a^{i}_{p}\})$ on $\Ccal_{i}$, and all these variables are independent, see e.g.~\cite{Neveu:1977mz}. This means that, in~(\ref{eq:NuAalphaiellRhoSR}) above, we have the property that, conditional on the $\sigma$-algebra generated by the sequence $(R^{A,i,0}_{n})_{n\geq 1}$, the variables $\Rho^{A,i,0}_{n}$ and $S^{A,i,0}_{n}$, for $n\geq 1$, are independent and distributed according to the uniform law on $(0,A)$ and the law $\mu_{i,R^{A,i,0}_{n}}$, respectively. Hence, the conditional probability that~(\ref{eq:xieiHail}) holds given that $\sigma$-algebra is
\[
\prod_{n=1}^{\lfloor 2^{j}/\delta\rfloor} \int_{\rho\in(0,A) \atop s\in\Ccal_{i}} \ind_{\left\{\left|\frac{\rho-\inpr{s}{t}}{\inpr{s}{e_{i}}}\right|\in[\delta,\infty]\cup\bigcup\limits_{q\in\Qrm} \lambda_{q}\right\}} \,\frac{\dd\rho}{A}\,\mu_{i,R^{A,i,0}_{n}}(\dd s).
\]
Exploiting the symmetry of $\nu$, it is easy to check that this is equal to
\[
\prod_{n=1}^{\lfloor 2^{j}/\delta\rfloor} \left(1-\frac{\delta-\lfloor\theta 2^{j}\rfloor\delta 2^{-j}}{A}\int_{s\in\Ccal_{i}}|\inpr{s}{e_{i}}|\,\mu_{i,R^{A,i,0}_{n}}(\dd s)\right)\leq\exp\left(-\frac{(1-\theta)\delta}{A\sqrt{d}}\left\lfloor\frac{2^{j}}{\delta}\right\rfloor\right),
\]
for $\delta<A-\|t\|$. In order to derive the upper bound, we have used the fact that $1-z\leq\ee^{-z}$ for any real $z$ and that $|\inpr{s}{e_{i}}|\geq d^{-1/2}$ for any $s\in\Ccal_{i}$. As a result,
\[
\prob\left(\forall n\leq\left\lfloor\frac{2^{j}}{\delta}\right\rfloor \quad |\xi^{e_{i}}_{t}(H^{A,i,0}_{n})|\in[\delta,\infty]\cup\bigcup_{q\in\Qrm} \lambda_{q}\right)\leq\exp\left(-\frac{1-\theta}{A\sqrt{d}}\,2^{j-1}\right).
\]
Now, recall that the set $\Qrm$ is a subset of $\{0,\ldots,2^{j}-1\}$ with cardinality $\lfloor\theta 2^{j}\rfloor$. We finally infer that
\[
\prob(\#\Qrm^{e_{i}}_{t,\delta,j}(H^{A,i,0})\leq\theta 2^{j})\leq\binom{2^{j}}{\lfloor\theta 2^{j}\rfloor}\exp\left(-\frac{1-\theta}{A\sqrt{d}}\,2^{j-1}\right).
\]
Making use of Stirling's formula, we get
\begin{equation}\label{eq:majlogprobQalpha0}
\limsup_{j\to\infty}\frac{1}{2^{j}}\log\prob(\#\Qrm^{e_{i}}_{t,\delta,j}(H^{A,i,0})\leq\theta 2^{j}) \leq -\Gamma(\theta)-\frac{1-\theta}{2A\sqrt{d}},
\end{equation}
with $\Gamma(\theta)=\theta\log\theta+(1-\theta)\log(1-\theta)$. Clearly, there exists a unique $\theta_{0}\in(0,1)$ at which the right-hand side of~(\ref{eq:majlogprobQalpha0}) vanishes. Moreover, this right-hand side is negative for any $\theta\in(0,\theta_{0})$. Using the Borel-Cantelli lemma, we see that for any $\delta\in(0,A-\|t\|)$ and any such $\theta$, almost surely, $\underline{\Qrm}^{e_{i}}_{t,\delta}(H^{A,i,0})\geq\theta$. The result follows.

Now, if $\ell=1$, we may apply Lemma~\ref{lem:asymptPi} with the image under $(\rho,s,r)\mapsto r$ of the measure $\leb^{1}_{|(0,A)}\otimes\nu_{\alpha,i,1}$, where $\leb^{1}_{|(0,A)}$ is the restriction of $\leb^{1}$ to $(0,A)$. Consequently, we infer that $\Phi^{A}_{i,1}\sim_{0}A\,\ph_{i,1}$ almost surely, where
\[
\forall w\in(0,1] \qquad
\left\{\begin{array}{l}
\Phi^{A}_{i,1}(w)=\Nu^{A,i,1}_{\alpha}((0,A)\times\Ccal_{i}\times[w,1])\\[1mm]
\ph_{i,1}(w)=\nu_{\alpha,i,1}(\Ccal_{i}\times[w,1]).
\end{array}\right.
\]
This directly implies that with probability one,
\[
C=\sup_{w\in(0,1]}\frac{\Phi^{A}_{i,1}(w)}{A(1+\ph_{i,1}(w))}\in[1,\infty).
\]
Furthermore, let us consider an integer $k\geq 1$ and an integer $j$ large enough to ensure that $1+\ph_{i,1}(w)\leq 2^{j}/(A\delta k)$ for some $w\in(0,1]$. As $(R^{A,i,1}_{n})_{n\geq 1}$ is nonincreasing, we have $n\leq\Phi^{A}_{i,1}(w)\leq C2^{j}/(\delta k)$ for any $n\geq 1$ with $R^{A,i,1}_{n}\geq w$. Therefore, if $C\leq k$ and~(\ref{eq:xieiHail}) hold simultaneously, then no integer $n\geq 1$ can satisfy both $R^{A,i,1}_{n}\geq w$ and $|\xi^{e_{i}}_{t}(H^{A,i,1}_{n})|\in[0,\delta)\setminus\bigcup_{q\in\Qrm}\lambda_{q}$. This happens with probability $\ee^{-I(w)}$, where
\[
I(w)=\int_{(\rho,s)\in(0,A)\times\Ccal_{i}\atop r\in[w,1]}\ind_{\left\{\left|\frac{\rho-\inpr{s}{t}}{\inpr{s}{e_{i}}}\right|\in[0,\delta)\setminus\bigcup\limits_{q\in\Qrm}\lambda_{q}\right\}}\,\dd\rho\,\nu_{\alpha,i,1}(\dd s,\dd r)\geq \frac{(1-\theta)\delta}{\sqrt{d}}\ph_{i,1}(w).
\]
Here, the lower bound holds for $\delta<A-\|t\|$ and results from standard computations that exploit the symmetry of $\nu$. This leads to an upper bound on the probability that $C\leq k$ and~(\ref{eq:xieiHail}) hold simultaneously, which may be optimized by letting $w\downarrow w_{*}$, where $w_{*}$ is the infimum of all $w\in(0,1]$ with $1+\ph_{i,1}(w)\leq 2^{j}/(A\delta k)$. By definition of $w_{*}$ and $\Acal_{i,1}$, and as $\ph_{i,1}$ is left-continuous with right limits, we have
\[
\lim_{w\downarrow w_{*}}\ph_{i,1}(w)\geq\ph_{i,1}(w^{*})-1\geq\frac{2^{j}}{A\delta k}-2\geq\frac{2^{j-1}}{A\delta k},
\]
where the last inequality holds for $j$ large enough. We deduce that
\[
\prob\left(C\leq k\mbox{ and }\forall n\leq\left\lfloor\frac{2^{j}}{\delta}\right\rfloor \, |\xi^{e_{i}}_{t}(H^{A,i,1}_{n})|\in[\delta,\infty]\cup\bigcup_{q\in\Qrm} \lambda_{q}\right)\leq\exp\left(-\frac{1-\theta}{Ak\sqrt{d}}\,2^{j-1}\right).
\]
Just as in the previous case, this implies that for any $k\geq 1$ with $\prob(C\leq k)>0$,
\[
\limsup_{j\to\infty}\frac{1}{2^{j}}\log\prob(\#\Qrm^{e_{i}}_{t,\delta,j}(H^{A,i,1})\leq\theta 2^{j}\:|\:C\leq k) \leq -\Gamma(\theta)-\frac{1-\theta}{2Ak\sqrt{d}},
\]
where the right-hand side vanishes for a unique $\theta_{k}\in(0,1)$. Using the Borel-Cantelli lemma as above, we see that for any $\delta\in(0,A-\|t\|)$ and any $\theta\in(0,\theta_{k})$, conditionally on $C\leq k$, with probability one, $\underline{\Qrm}^{e_{i}}_{t,\delta}(H^{A,i,1})\geq\theta$. Therefore, with probability one, conditionally on $C\leq k$, we have $\limsup_{\delta\to 0}\underline{\Qrm}^{e_{i}}_{t,\delta}(H^{A,i,1})\geq\theta_{k}$. The result now follows from the fact that $C<\infty$ almost surely.
\end{proof}

%%%%%%%%%%%%%%%%%%%%%%%%%%%%%%%%%%%%%%%%%%%%%%%%%%%%%%%%%%%%%%%%%%%%%%%%%%%%%%
%%%%%%%%%%%%%%%%%%%%%%%%%%%%%%%%%%%%%%%%%%%%%%%%%%%%%%%%%%%%%%%%%%%%%%%%%%%%%%
\subsection{Ubiquity}
%%%%%%%%%%%%%%%%%%%%%%%%%%%%%%%%%%%%%%%%%%%%%%%%%%%%%%%%%%%%%%%%%%%%%%%%%%%%%%
%%%%%%%%%%%%%%%%%%%%%%%%%%%%%%%%%%%%%%%%%%%%%%%%%%%%%%%%%%%%%%%%%%%%%%%%%%%%%%

Last, the proof of Theorem~\ref{thm:sliKnualpha} calls upon the next lemma, which is a straightforward consequence of Theorem~3.6 in~\cite{Durand:2008jk}. (The hypotheses of that theorem are verified here because the diameter of the set of $t\in\R e_{i}$ such that $\dist(t,H^{A,i,\ell}_{n})<1$ is bounded above by $2d^{1/2}$, since $S^{A,i,\ell}_{n}\in\Ccal_{i}$.) For every $g\in\gauge_{1}$, let
\[
F^{A,i,\ell}(g)=\{ t\in\R^{d} \:|\: \dist(t,H^{A,i,\ell}_{n})<g(R^{A,i,\ell}_{n}) \mbox{ for i.m. } n\geq 1 \}.
\]

\begin{lem}\label{lem:ubiqFAilg}
Let $A\geq 1$, $i\in\{1,\ldots,d\}$ and $\ell\in\{0,1\}$ such that $\nu_{\alpha,i,\ell}$ has infinite total mass, and let $g\in\gauge_{1}$. If the set $F^{A,i,\ell}(\underline{g})$ has full Lebesgue measure in $\opball{0}{A}$ for some $\underline{g}\in\gauge_{1}$ with $g\prec\underline{g}$, then $F^{A,i,\ell}(r\mapsto r)\in\lic{d-1,g}{\opball{0}{A}}$.
\end{lem}

This lemma falls into the category of the ubiquity results obtained in~\cite{Durand:2007uq,Durand:2008jk,Durand:2008kx}, which enable one to deduce the large intersection properties of a set, such as $F^{A,i,\ell}(r\mapsto r)$ in the present situation, from the sole knowledge of the Lebesgue measure of a corresponding enlarged set, which is $F^{A,i,\ell}(\underline{g})$ here.

%%%%%%%%%%%%%%%%%%%%%%%%%%%%%%%%%%%%%%%%%%%%%%%%%%%%%%%%%%%%%%%%%%%%%%%%%%%%%%
%%%%%%%%%%%%%%%%%%%%%%%%%%%%%%%%%%%%%%%%%%%%%%%%%%%%%%%%%%%%%%%%%%%%%%%%%%%%%%
\subsection{End of the proof when $h_{\nu}(g)>\alpha$}
%%%%%%%%%%%%%%%%%%%%%%%%%%%%%%%%%%%%%%%%%%%%%%%%%%%%%%%%%%%%%%%%%%%%%%%%%%%%%%
%%%%%%%%%%%%%%%%%%%%%%%%%%%%%%%%%%%%%%%%%%%%%%%%%%%%%%%%%%%%%%%%%%%%%%%%%%%%%%

To begin with, let us recall that $\Rho_{n}$ and $X_{n}$ are defined in terms of the atoms of the Poisson random measure $\Nu$, see Section~\ref{sec:defnonGaussianLevy}. Then, let $R_{n}=|X_{n}|^{1/\alpha}$ for any $n\geq 1$ and, given an integer $A\geq 1$, let $\Ncal_{A}$ denote the set of all $n\geq 1$ such that $\Rho_{n}<A$ and $R_{n}\leq 1$.

For any $A\geq 1$, we may apply Lemma~\ref{lem:convPi} with the image under $(\rho,s,r)\mapsto r$ of the measure $\leb^{1}_{|(0,A)}\otimes\nu_{\alpha}$. Consequently, in view of~(\ref{eq:defhnug}), we deduce that with probability one, for any $A\geq 1$ and any $g\in\gauge_{1}$ with $h_{\nu}(g)>\alpha$, the series $\sum_{n\in\Ncal_{A}}g(R_{n})$ converges. In particular, for any $\delta>0$, there necessarily exists an integer $n_{0}\geq 1$ such that $R_{n}<\delta$ for any $n\in\Ncal_{A}$ with $n\geq n_{0}$. Then, for any $n_{1}\geq n_{0}$,
\[
K_{\nu}(\alpha)\cap\opball{0}{A-1}\subseteq\bigcup_{n\in\Ncal_{A} \atop n\geq n_{1}}\left\{ t\in\opball{0}{A} \:\bigl|\: \dist(t,H_{n})<R_{n} \right\},
\]
where $H_{n}$ is the hyperplane defined by~(\ref{eq:defHn}). Moreover, each set in the union above may be covered by $(3\lfloor 2A\sqrt{d}/R_{n}\rfloor)^{d-1}$ open balls with radius $2R_{n}$. Therefore,
\begin{align*}
\hau^{r\mapsto r^{d-1}g(r)}_{\delta}(K_{\nu}(\alpha)\cap\opball{0}{A-1}) &\leq \sum_{n\in\Ncal_{A} \atop n\geq n_{1}} \left(3\left\lfloor\frac{2A\sqrt{d}}{R_{n}}\right\rfloor\right)^{d-1} (4R_{n})^{d-1} g(4R_{n}) \\
&\leq 4(24A\sqrt{d})^{d-1} \sum_{n\in\Ncal_{A} \atop n\geq n_{1}} g(R_{n}).
\end{align*}
Letting $n_{1}\to\infty$ and $\delta\to 0$, we deduce that $\hau^{d-1,g}(K_{\nu}(\alpha)\cap\opball{0}{A-1})=0$. This holds for all integers $A\geq 1$, so the result follows.

%%%%%%%%%%%%%%%%%%%%%%%%%%%%%%%%%%%%%%%%%%%%%%%%%%%%%%%%%%%%%%%%%%%%%%%%%%%%%%
%%%%%%%%%%%%%%%%%%%%%%%%%%%%%%%%%%%%%%%%%%%%%%%%%%%%%%%%%%%%%%%%%%%%%%%%%%%%%%
\subsection{End of the proof when $h_{\nu}(g)<\alpha$}
%%%%%%%%%%%%%%%%%%%%%%%%%%%%%%%%%%%%%%%%%%%%%%%%%%%%%%%%%%%%%%%%%%%%%%%%%%%%%%
%%%%%%%%%%%%%%%%%%%%%%%%%%%%%%%%%%%%%%%%%%%%%%%%%%%%%%%%%%%%%%%%%%%%%%%%%%%%%%

By Lemma~\ref{lem:convPi}, the next statement holds almost surely: for all integers $A\geq 1$, $i\in\{1,\ldots,d\}$ and $\ell\in\{0,1\}$, and all $g\in\gauge_{1}$,
\begin{equation}\label{eq:intsumgnualphail}
\int_{s\in\S^{d-1} \atop r\in (0,1]} g(r) \,\nu_{\alpha,i,\ell}(\dd s,\dd r)=\infty \qquad\Longrightarrow\qquad \sum_{n=1}^{\infty} g(R^{A,i,\ell}_{n})=\infty.
\end{equation}
Moreover, the statement of Lemma~\ref{lem:HAilhomdist} holds with probability one as well. We shall work in what follows on the almost sure event on which these two statements hold. Let us consider a gauge function $g\in\gauge_{1}$ such that $h_{\nu}(g)<\alpha$. Due to~(\ref{eq:defhnug}) and the fact that $\S^{d-1}\times(0,1]$ is covered by the sets $\Ccal_{i}\times\Acal_{i,\ell}$, the integral in~(\ref{eq:intsumgnualphail}) is infinite for some $(i,\ell)\in\{1,\ldots,d\}\times\{0,1\}$.

Let us assume that $g\prec (r\mapsto r)$. Borrowing ideas from the proof of~\cite[Proposition~5]{Durand:2007fk}, we may build $\underline{g}\in\gauge_{1}$ such that $g\prec\underline{g}$ and the integral in~(\ref{eq:intsumgnualphail}) with $g$ replaced by $\underline{g}$ is infinite as well. Therefore, $\sum_{n}\underline{g}(R^{A,i,\ell}_{n})$ diverges, so that $F^{A,i,\ell}(\underline{g})$ has full Lebesgue measure in $\opball{0}{A}$, by virtue of Lemma~\ref{eq:homdistfull}. Thanks to Lemma~\ref{lem:ubiqFAilg}, we deduce that $F^{A,i,\ell}(r\mapsto r)\in\lic{d-1,g}{\opball{0}{A}}$. The same result holds if $g\not\prec (r\mapsto r)$, in view of the divergence of $\sum_{n=1}^{\infty} R^{A,i,\ell}_{n}$, combined with Lemma~\ref{eq:homdistfull} and Theorem~\ref{thm:sli}(\ref{item:slifull}).

To conclude, it is crucial to note that the set $F^{A,i,\ell}(r\mapsto r)$ is contained in $K_{\nu}(\alpha)$. Given that $F^{A,i,\ell}(r\mapsto r)$ is in $\lic{d-1,g}{\opball{0}{A}}$, the $G_{\delta}$-set $K_{\nu}(\alpha)$ belongs to the same class. Now, let $W$ be a bounded open subset of $\R^{d}$. Then, $W\subseteq\opball{0}{A}$ for $A$ large enough, so that $\netm^{f}_{\infty}(K_{\nu}(\alpha)\cap W)=\netm^{f}_{\infty}(W)$ for every $f\in\gauge_{d}$ with $f\prec (r\mapsto r^{d-1}g(r))$. Lemma~12 in~\cite{Durand:2007uq} finally ensures that $K_{\nu}(\alpha)\in\lic{d-1,g}{\R^{d}}$.

%%%%%%%%%%%%%%%%%%%%%%%%%%%%%%%%%%%%%%%%%%%%%%%%%%%%%%%%%%%%%%%%%%%%%%%%%%%%%%
%%%%%%%%%%%%%%%%%%%%%%%%%%%%%%%%%%%%%%%%%%%%%%%%%%%%%%%%%%%%%%%%%%%%%%%%%%%%%%
%%%%%%%%%%%%%%%%%%%%%%%%%%%%%%%%%%%%%%%%%%%%%%%%%%%%%%%%%%%%%%%%%%%%%%%%%%%%%%
\section{Proofs concerning the jump component}\label{sec:proofsecregnonGaussian}
%%%%%%%%%%%%%%%%%%%%%%%%%%%%%%%%%%%%%%%%%%%%%%%%%%%%%%%%%%%%%%%%%%%%%%%%%%%%%%
%%%%%%%%%%%%%%%%%%%%%%%%%%%%%%%%%%%%%%%%%%%%%%%%%%%%%%%%%%%%%%%%%%%%%%%%%%%%%%
%%%%%%%%%%%%%%%%%%%%%%%%%%%%%%%%%%%%%%%%%%%%%%%%%%%%%%%%%%%%%%%%%%%%%%%%%%%%%%

Throughout the section, we assume that the measure $\nu$ is admissible. We now establish Proposition~\ref{prp:existpathsLnu} and the results of Section~\ref{sec:regnonGaussian}.

%%%%%%%%%%%%%%%%%%%%%%%%%%%%%%%%%%%%%%%%%%%%%%%%%%%%%%%%%%%%%%%%%%%%%%%%%%%%%%
%%%%%%%%%%%%%%%%%%%%%%%%%%%%%%%%%%%%%%%%%%%%%%%%%%%%%%%%%%%%%%%%%%%%%%%%%%%%%%
\subsection{Proof of Proposition~\ref{prp:existpathsLnu}}\label{subsec:proofprpexistpathsLnu}
%%%%%%%%%%%%%%%%%%%%%%%%%%%%%%%%%%%%%%%%%%%%%%%%%%%%%%%%%%%%%%%%%%%%%%%%%%%%%%
%%%%%%%%%%%%%%%%%%%%%%%%%%%%%%%%%%%%%%%%%%%%%%%%%%%%%%%%%%%%%%%%%%%%%%%%%%%%%%

Let $A\in\N$ and $U_{A}=\clball{0}{A+1}\cap(d^{-1/2}\Z^{d})$. For any fixed $t\in\R^{d}$, the series $\sum_{j\geq 0} L_{\nu,j}(t)$ defining $L_{\nu}(t)$ converges almost surely, see Section~\ref{sec:defnonGaussianLevy}. Thus, the event $\event_{A}$ consisting in the fact that the series $\sum_{j\geq 0} L_{\nu,j}(u)$, for $u\in U_{A}$, converge simultaneously has probability one. Furthermore, it follows from Lemma~\ref{lem:oscLnuj} that the event $\event'_{A}=\{\Zeta_{\nu}(A+1)<\infty\}$ has probability one too.

Let us now assume that the almost sure event $\event_{A}\cap\event'_{A}$ happens, and let $\eps>0$ and $t\in\clball{0}{A}$. Then, $\|t-u\|\leq 1/2$ for some $u\in U_{A}$. The sum $\chi_{\nu}$ defined by~(\ref{eq:defchinu}) is finite and the series $\sum_{j\geq 0}L_{\nu,j}(u)$ converges, so that
\[
\sum_{j=j_{1}}^{j_{2}}2^{-j}(j+(j\nu_{j})^{1/2})\leq\eps \qquad\mbox{and}\qquad \left|\sum_{j=j_{1}}^{j_{2}} L_{\nu,j}(u)\right|\leq\eps
\]
for all integers $j_{1}$ and $j_{2}$ such that $j_{2}\geq j_{1}\geq j_{0}$, and some $j_{0}\geq 1$. Therefore,
\[
\left|\sum_{j=j_{1}}^{j_{2}} L_{\nu,j}(t)\right|\leq\left|\sum_{j=j_{1}}^{j_{2}} L_{\nu,j}(u)\right|+\sum_{j=j_{1}}^{j_{2}} \zeta_{\nu}(A+1,j,1)\leq (1+\Zeta_{\nu}(A+1))\,\eps.
\]
The partial sums of $\sum_{j\geq 0}L_{\nu,j}(t)$ form a Cauchy sequence. So, for any $A\in\N$, with probability one, $\sum_{j\geq 0}L_{\nu,j}(t)$ converges for any $t\in\clball{0}{A}$, and the result follows.

%%%%%%%%%%%%%%%%%%%%%%%%%%%%%%%%%%%%%%%%%%%%%%%%%%%%%%%%%%%%%%%%%%%%%%%%%%%%%%
%%%%%%%%%%%%%%%%%%%%%%%%%%%%%%%%%%%%%%%%%%%%%%%%%%%%%%%%%%%%%%%%%%%%%%%%%%%%%%
\subsection{Proof of Theorems~\ref{thm:sizepropEnuh} and~\ref{thm:lipropEnuh}}\label{subsec:proofslipropEnuh}
%%%%%%%%%%%%%%%%%%%%%%%%%%%%%%%%%%%%%%%%%%%%%%%%%%%%%%%%%%%%%%%%%%%%%%%%%%%%%%
%%%%%%%%%%%%%%%%%%%%%%%%%%%%%%%%%%%%%%%%%%%%%%%%%%%%%%%%%%%%%%%%%%%%%%%%%%%%%%

The statement of Corollary~\ref{cor:characLnuKnualpha} holds with probability one, and that of Theorem~\ref{thm:sliKnualpha} holds with probability one for all rationals $\alpha\in\Q\cap(0,\infty)$ simultaneously. Let us assume that the almost sure event on which these statements hold occurs. Theorems~\ref{thm:sizepropEnuh} and~\ref{thm:lipropEnuh} follow from a series of propositions that we now state and establish. Throughout, we consider a real $h\in[0,1/\beta_{\nu})$, a gauge function $g\in\gauge_{1}$ and a nonempty open set $W\subseteq\R^{d}$.

\begin{prp}\label{prp:hauEnuhzero}
If $h<h_{\nu}(g)$, then $\hau^{d-1,g}(E_{\nu}(h))=\hau^{d-1,g}(E'_{\nu}(h))=0$.
\end{prp}

\begin{proof}
There exists a rational $\alpha\in\Q\cap(h,h_{\nu}(g))$ such that $E'_{\nu}(h)\subseteq K_{\nu}(\alpha)$. Hence, $\hau^{d-1,g}(E'_{\nu}(h))\leq\hau^{d-1,g}(K_{\nu}(\alpha))=0$. Furthermore, if $h>0$, then $E_{\nu}(h)\subseteq E'_{\nu}(h)$, so $\hau^{d-1,g}(E_{\nu}(h))=0$ as well. This result still holds for $h=0$, because $E_{\nu}(0)=E'_{\nu}(0)\cup J_{\nu}$ and $J_{\nu}$ is a countable union of hyperplanes.
\end{proof}

\begin{prp}
If $h<h_{\nu}(g)$, then $E'_{\nu}(h)\not\in\lic{d-1,g}{W}$.
\end{prp}

\begin{proof}
Adapting the method developed in the proof of~\cite[Proposition~3]{Durand:2007fk}, we may build a gauge function $\overline{g}\in\gauge_{1}$ satisfying both $\overline{g}\prec g$ and $h_{\nu}(\overline{g})\geq h_{\nu}(g)$. Then, applying Proposition~\ref{prp:hauEnuhzero} with $\overline{g}$ instead of $g$, we infer that $\hau^{d-1,\overline{g}}(E'_{\nu}(h)\cap W)=0$. We conclude by Theorem~\ref{thm:sli}(\ref{item:slisize}).
\end{proof}

\begin{prp}\label{prp:Enuhinlic}
If $h\geq h_{\nu}(g)$, then $E'_{\nu}(h)\in\lic{d-1,g}{W}$.
\end{prp}

\begin{proof}
The mapping $\alpha\mapsto K_{\nu}(\alpha)$ being nondecreasing, we have
\[
E'_{\nu}(h)=(\R^{d}\setminus J_{\nu})\cap\bigcap_{h<\alpha\leq 1/\beta_{\nu} \atop \alpha\in\Q}K_{\nu}(\alpha).
\]
Furthermore, each set $K_{\nu}(\alpha)$ arising in this last intersection belongs to $\lic{d-1,g}{W}$. The set $\R^{d}\setminus J_{\nu}$ belongs to this class as well by virtue of Theorem~\ref{thm:sli}(\ref{item:slifull}), because $J_{\nu}$ is the union of countably many hyperplanes. We conclude using Theorem~\ref{thm:sli}(\ref{item:sliclosedint}).
\end{proof}

\begin{prp}\label{prp:hauEnuhinfty}
If $h\geq h_{\nu}(g)$, then $\hau^{d-1,g}(E'_{\nu}(h)\cap W)=\infty$.
\end{prp}

\begin{proof}
The assumption of the proposition implies that $g\prec(r\mapsto r)$. Indeed, otherwise, we would clearly have $h_{\nu}(g)\geq 1/\beta_{\nu}$, which is in contradiction with the fact that $h_{\nu}(g)\leq h<1/\beta_{\nu}$. Therefore, borrowing ideas from the proof of~\cite[Proposition~5]{Durand:2007fk}, we may build a gauge function $\underline{g}\in\gauge_{1}$ satisfying both $g\prec\underline{g}$ and $h_{\nu}(\underline{g})\leq h_{\nu}(g)$. Then, applying Proposition~\ref{prp:Enuhinlic} with the gauge function $\underline{g}$ instead of $g$, we see that $E'_{\nu}(h)$ belongs to $\lic{d-1,\underline{g}}{W}$. We conclude by Theorem~\ref{thm:sli}(\ref{item:slisize}).
\end{proof}

\begin{prp}
If $h=h_{\nu}(g)$, then $\hau^{d-1,g}(E_{\nu}(h)\cap W)=\infty$.
\end{prp}

\begin{proof}
In the case where $h=0$, the result follows directly from Proposition~\ref{prp:hauEnuhinfty}, because $E_{\nu}(0)$ contains $E'_{\nu}(0)$. In the case where $h>0$, it suffices to make use of Proposition~\ref{prp:hauEnuhinfty} again, together with the observation that
\[
E_{\nu}(h)=E'_{\nu}(h)\setminus\bigcup_{0<\alpha<h \atop \alpha\in\Q} K_{\nu}(\alpha),
\]
because the mapping $\alpha\mapsto K_{\nu}(\alpha)$ is nondecreasing, and that each set $K_{\nu}(\alpha)$ in the union above has Hausdorff measure zero for the gauge function $r\mapsto r^{d-1}g(r)$.
\end{proof}

%%%%%%%%%%%%%%%%%%%%%%%%%%%%%%%%%%%%%%%%%%%%%%%%%%%%%%%%%%%%%%%%%%%%%%%%%%%%%%
%%%%%%%%%%%%%%%%%%%%%%%%%%%%%%%%%%%%%%%%%%%%%%%%%%%%%%%%%%%%%%%%%%%%%%%%%%%%%%
\subsection{Proof of Corollary~\ref{cor:sizepropEnuh}}\label{subsec:proofcorsizepropEnuh}
%%%%%%%%%%%%%%%%%%%%%%%%%%%%%%%%%%%%%%%%%%%%%%%%%%%%%%%%%%%%%%%%%%%%%%%%%%%%%%
%%%%%%%%%%%%%%%%%%%%%%%%%%%%%%%%%%%%%%%%%%%%%%%%%%%%%%%%%%%%%%%%%%%%%%%%%%%%%%

Let us assume that the almost sure event on which the statement of Theorem~\ref{thm:sizepropEnuh} holds occurs. Let $h\in[0,1/\beta_{\nu})$, $g\in\gauge_{1}^{*}$ and $W$ be a nonempty open set. Most of the first part of the corollary, which gives $\hau^{d-1,g}(E_{\nu}(h)\cap W)$ and $\hau^{d-1,g}(E'_{\nu}(h)\cap W)$, follows directly from Theorem~\ref{thm:sizepropEnuh}. The only new property is that $\hau^{d-1,g}(E_{\nu}(h)\cap W)=\infty$ when $h>h_{\nu}(g)$. To prove this, let $\gamma_{g}=\lim_{r\to 0}(\log g(r))/\log r$. Then, let us assume that $h>h_{\nu}(g)$ and remark that $h_{\nu}(g)=\gamma_{g}/\beta_{\nu}$. It follows that $\gamma_{g}<\beta_{\nu}h$, so that $g(r)\geq r^{\beta_{\nu}h}$ for $r\geq 0$ small enough. Thus, the Hausdorff measure of $E_{\nu}(h)\cap W$ for the gauge function $r\mapsto r^{d-1}g(r)$ is larger than or equal to its $(d-1+\beta_{\nu}h)$-dimensional measure, which is infinite as a result of Theorem~\ref{thm:sizepropEnuh} and the fact that $h_{\nu}(r\mapsto r^{\beta_{\nu}h})=h$.

The proof of the second part of the corollary, which gives the value of the Hausdorff dimension of $E_{\nu}(h)\cap W$ and $E'_{\nu}(h)\cap W$, is a consequence of the first part, together with the following observations: for $h>0$, note that $h_{\nu}(r\mapsto r^{s})=s/\beta_{\nu}$ for any $s\in(0,1]$; for $h=0$, note that $h_{\nu}(g)=0$ for some $g\in\gauge_{1}$ (which may be built by borrowing ideas from the proof of~\cite[Proposition~5]{Durand:2007fk}).

%%%%%%%%%%%%%%%%%%%%%%%%%%%%%%%%%%%%%%%%%%%%%%%%%%%%%%%%%%%%%%%%%%%%%%%%%%%%%%
%%%%%%%%%%%%%%%%%%%%%%%%%%%%%%%%%%%%%%%%%%%%%%%%%%%%%%%%%%%%%%%%%%%%%%%%%%%%%%
\subsection{Proof of Proposition~\ref{prp:Enuhnufinite}}\label{subsec:proofprpEnuhnufinite}
%%%%%%%%%%%%%%%%%%%%%%%%%%%%%%%%%%%%%%%%%%%%%%%%%%%%%%%%%%%%%%%%%%%%%%%%%%%%%%
%%%%%%%%%%%%%%%%%%%%%%%%%%%%%%%%%%%%%%%%%%%%%%%%%%%%%%%%%%%%%%%%%%%%%%%%%%%%%%

As mentioned at the beginning of Section~\ref{sec:proofthmsliKnualpha}, if $\nu$ has finite total mass, then with probability one, $K_{\nu}(\alpha)=\emptyset$ for all $\alpha>0$. The result now follows from Corollary~\ref{cor:characLnuKnualpha}.

%%%%%%%%%%%%%%%%%%%%%%%%%%%%%%%%%%%%%%%%%%%%%%%%%%%%%%%%%%%%%%%%%%%%%%%%%%%%%%
%%%%%%%%%%%%%%%%%%%%%%%%%%%%%%%%%%%%%%%%%%%%%%%%%%%%%%%%%%%%%%%%%%%%%%%%%%%%%%
\subsection{Proof of Proposition~\ref{prp:sizeEnuhsimple}}\label{subsec:proofprpsizeEnuhsimple}
%%%%%%%%%%%%%%%%%%%%%%%%%%%%%%%%%%%%%%%%%%%%%%%%%%%%%%%%%%%%%%%%%%%%%%%%%%%%%%
%%%%%%%%%%%%%%%%%%%%%%%%%%%%%%%%%%%%%%%%%%%%%%%%%%%%%%%%%%%%%%%%%%%%%%%%%%%%%%

The case where $h>1/\beta_{\nu}$ follows from~(\ref{eq:sizeEnuhsimple}). To treat the case where $h=1/\beta_{\nu}$, let us use of Corollary~\ref{cor:characLnuKnualpha} in order to write that
\[
E'_{\nu}(1/\beta_{\nu})=\R^{d}\setminus J_{\nu} \qquad\mbox{and}\qquad E_{\nu}(1/\beta_{\nu})=\R^{d}\setminus\Biggl(J_{\nu}\cup\bigcup_{0<\alpha<1/\beta_{\nu} \atop \alpha\in\Q} K_{\nu}(\alpha)\Biggr).
\]
In the union above, we may restrict $\alpha$ to being rational, because $\alpha\mapsto K_{\nu}(\alpha)$ is nondecreasing. Now, applying Theorem~\ref{thm:sliKnualpha} with $g:r\mapsto r$, we infer that with probability one, the sets $K_{\nu}(\alpha)$, for $\alpha\in\Q\cap(0,1/\beta_{\nu})$, all have Lebesgue measure zero. Moreover, the set $J_{\nu}$ is a countable union of hyperplanes, thereby having Lebesgue measure zero as well. It follows that almost surely, the set $E_{\nu}(1/\beta_{\nu})$ has full Lebesgue measure in the whole space $\R^{d}$.

%%%%%%%%%%%%%%%%%%%%%%%%%%%%%%%%%%%%%%%%%%%%%%%%%%%%%%%%%%%%%%%%%%%%%%%%%%%%%%
%%%%%%%%%%%%%%%%%%%%%%%%%%%%%%%%%%%%%%%%%%%%%%%%%%%%%%%%%%%%%%%%%%%%%%%%%%%%%%
%%%%%%%%%%%%%%%%%%%%%%%%%%%%%%%%%%%%%%%%%%%%%%%%%%%%%%%%%%%%%%%%%%%%%%%%%%%%%%
\section{Estimates of the increments of $L_{\nu,j}$}\label{sec:prooflemoscLnuj}
%%%%%%%%%%%%%%%%%%%%%%%%%%%%%%%%%%%%%%%%%%%%%%%%%%%%%%%%%%%%%%%%%%%%%%%%%%%%%%
%%%%%%%%%%%%%%%%%%%%%%%%%%%%%%%%%%%%%%%%%%%%%%%%%%%%%%%%%%%%%%%%%%%%%%%%%%%%%%
%%%%%%%%%%%%%%%%%%%%%%%%%%%%%%%%%%%%%%%%%%%%%%%%%%%%%%%%%%%%%%%%%%%%%%%%%%%%%%

The purpose of this section is to establish Lemma~\ref{lem:oscLnuj}, that is, to prove the almost sure finiteness of $\Zeta_{\nu}(A)$, for any fixed integer $A\geq 1$.

%%%%%%%%%%%%%%%%%%%%%%%%%%%%%%%%%%%%%%%%%%%%%%%%%%%%%%%%%%%%%%%%%%%%%%%%%%%%%%
%%%%%%%%%%%%%%%%%%%%%%%%%%%%%%%%%%%%%%%%%%%%%%%%%%%%%%%%%%%%%%%%%%%%%%%%%%%%%%
\subsection{A net argument}
%%%%%%%%%%%%%%%%%%%%%%%%%%%%%%%%%%%%%%%%%%%%%%%%%%%%%%%%%%%%%%%%%%%%%%%%%%%%%%
%%%%%%%%%%%%%%%%%%%%%%%%%%%%%%%%%%%%%%%%%%%%%%%%%%%%%%%%%%%%%%%%%%%%%%%%%%%%%%

Given an integer $k\geq 1$, let $\sigma_{k}=2^{-k}/\lfloor d^{1/2}\rfloor$. There exists $U_{A,k}\subseteq\sigma_{k}\,\Z^d$ with cardinality at most $(2^{k+2}Ad^{1/2})^{d}$ such that $[-A-2^{-k},A+2^{-k}]^{d}$ is covered by the cubes $u+[0,\sigma_{k})^{d}$, for $u\in U_{A,k}$. Then, for any $t\in\Q^{d}\cap\clball{0}{A}$ and $\tau\in\Q^{d}\cap\clball{0}{2^{-k}}$, there are two points $u$ and $u'$ in $U_{A,k}$ such that $t-u$ and $t+\tau-u'$ belong to $[0,\sigma_{k})^{d}$. Moreover, writing $u=p\sigma_{k}$ and $u'=p'\sigma_{k}$ with $p,p'\in\Z^{d}$, we see that the $\ell^{1}$-norm of $p-p'$ is at most $3d$, so there is a finite sequence $(p_{i})_{0\leq i\leq n}$ in $\Z^d$ such that $n\leq 3d$, $p_{0}=p$, $p_{n}=p'$, and $\|p_{i+1}-p_{i}\|=1$ and $u_{i}=p_{i}\sigma_{k}\in U_{A,k}$ for all $i$. As a result, for any $j\geq 1$, the increment $|L_{\nu,j}(t+\tau)-L_{\nu,j}(t)|$ is at most
\[
|L_{\nu,j}(t+\tau)-L_{\nu,j}(u')| + |L_{\nu,j}(t)-L_{\nu,j}(u)|+\sum_{i=0}^{n-1}|L_{\nu,j}(u_{i+1})-L_{\nu,j}(u_{i})|,
\]
Making use of~(\ref{eq:defzetanuAjkalt}) and letting
\[
\zeta_{\nu}(t,j,k)=\sup_{\tau\in\Q^{d}\cap\clball{0}{2^{-k}}} |L_{\nu,j}(t+\tau)-L_{\nu,j}(t)|
\]
for any $j,k\geq 1$ and $t\in\R^{d}$, it follows that
\begin{equation}\label{eq:boundzetanuAjk}
\zeta_{\nu}(A,j,k)\leq (3d+2)\sup_{t\in U_{A,k}} \zeta_{\nu}(t,j,k).
\end{equation}

%%%%%%%%%%%%%%%%%%%%%%%%%%%%%%%%%%%%%%%%%%%%%%%%%%%%%%%%%%%%%%%%%%%%%%%%%%%%%%
%%%%%%%%%%%%%%%%%%%%%%%%%%%%%%%%%%%%%%%%%%%%%%%%%%%%%%%%%%%%%%%%%%%%%%%%%%%%%%
\subsection{Estimates of $\zeta_{\nu}(t,j,k)$}
%%%%%%%%%%%%%%%%%%%%%%%%%%%%%%%%%%%%%%%%%%%%%%%%%%%%%%%%%%%%%%%%%%%%%%%%%%%%%%
%%%%%%%%%%%%%%%%%%%%%%%%%%%%%%%%%%%%%%%%%%%%%%%%%%%%%%%%%%%%%%%%%%%%%%%%%%%%%%

Let us now derive an appropriate upper bound on $\zeta_{\nu}(t,j,k)$ and a control on the tail distribution of this bound. To this end, for any relatively compact Borel set $V\in\Bcal_{0}(\H_{d})$, let
\[
M_{\nu,j}(V)=\int_{(\rho,s)\in V\atop|x|\in\Ical_{j}} \Nu(\dd\rho,\dd s,\dd x) \qquad\mbox{and}\qquad m_{\nu,j}(V)=\int_{(\rho,s)\in V\atop|x|\in\Ical_{j}} \dd\rho\,\nu(\dd s,\dd x).
\]
Recall that $\Lfrak_{\nu,j}(V)$ is given by~(\ref{eq:defLfraknuj}), so that we clearly have
\begin{equation}\label{eq:boundLfraknujV}
|\Lfrak_{\nu,j}(V)|\leq 2^{-j+1}\left(M_{\nu,j}(V)+m_{\nu,j}(V)\right).
\end{equation}
Moreover, for any $t\in\R^{d}$ and any $\delta>0$, let $V^{\circ}_{t,\delta}=V^{+}_{t,\delta}\setminus V^{-}_{t,\delta}$, where
\[
V^{+}_{t,\delta}=\bigcup_{\tau\in\Q^{d}\cap\clball{0}{\delta}} V_{t+\tau} \qquad\mbox{and}\qquad V^{-}_{t,\delta}=\bigcap_{\tau\in\Q^{d}\cap\clball{0}{\delta}} V_{t+\tau}.
\]
Note that $V^{\circ}_{t,\delta}\in\Bcal_{0}(\H_{d})$ and that $|\rho-\inpr{s}{t}|\leq\delta$ for any $(\rho,s)\in V^{\circ}_{t,\delta}$. Therefore, exploiting the symmetry of $\nu$, we have
\begin{equation}\label{eq:boundmnujVcirc}
m_{\nu,j}(V^{\circ}_{t,\delta})\leq\frac{1}{2}\int_{s\in\S^{d-1}\atop|x|\in\Ical_{j}}\int_{\rho\in\R}\ind_{\{|\rho-\inpr{s}{t}|\leq\delta\}}\,\dd\rho\,\nu(\dd s,\dd x)=2\delta\nu_{j}.
\end{equation}
Our approach now depends on whether or not $2^{-k}\nu_{j}\leq\eta jk$, where $\eta\geq 1$ is a real constant to be fixed later.

%%%%%%%%%%%%%%%%%%%%%%%%%%%%%%%%%%%%%%%%%%%%%%%%%%%%%%%%%%%%%%%%%%%%%%%%%%%%%%
%%%%%%%%%%%%%%%%%%%%%%%%%%%%%%%%%%%%%%%%%%%%%%%%%%%%%%%%%%%%%%%%%%%%%%%%%%%%%%
\subsubsection{Case where $2^{-k}\nu_{j}\leq\eta jk$}
%%%%%%%%%%%%%%%%%%%%%%%%%%%%%%%%%%%%%%%%%%%%%%%%%%%%%%%%%%%%%%%%%%%%%%%%%%%%%%
%%%%%%%%%%%%%%%%%%%%%%%%%%%%%%%%%%%%%%%%%%%%%%%%%%%%%%%%%%%%%%%%%%%%%%%%%%%%%%

Here, the suitable bound on $\zeta_{\nu}(t,j,k)$ and an estimate of its tail distribution are given by the next two results.

\begin{lem}\label{lem:boundzetanutjk1}
For any $t\in\R^d$ and any $\eta,j,k\geq 1$ with $2^{-k}\nu_{j}\leq\eta jk$,
\[
\zeta_{\nu}(t,j,k)\leq 2^{-j+1} ( M_{\nu,j}(V^{\circ}_{t,2^{-k}}) + 2\eta jk ).
\]
\end{lem}

\begin{proof}
For every $\tau\in\Q^{d}\cap\clball{0}{2^{-k}}$, the increment $|L_{\nu,j}(t+\tau)-L_{\nu,j}(t)|$ is equal to
\begin{align*}
& \left|(\Lfrak_{\nu,j}(V_{t+\tau}\setminus V_{t})+\Lfrak_{\nu,j}(V_{t+\tau}\cap V_{t}))-(\Lfrak_{\nu,j}(V_{t}\setminus V_{t+\tau})+\Lfrak_{\nu,j}(V_{t}\cap V_{t+\tau}))\right|\\
\leq & 2^{-j+1}\left( M_{\nu,j}(V_{t+\tau}\Delta V_{t}) + m_{\nu,j}(V_{t+\tau}\Delta V_{t}) \right),
\end{align*}
where the last bound is due to~(\ref{eq:boundLfraknujV}). Furthermore, $V_{t+\tau}\Delta V_{t}\subseteq V^{\circ}_{t,2^{-k}}$, so that
\[
|L_{\nu,j}(t+\tau)-L_{\nu,j}(t)|\leq 2^{-j+1}( M_{\nu,j}(V^{\circ}_{t,2^{-k}}) + m_{\nu,j}(V^{\circ}_{t,2^{-k}}) )
\]
The result follows from~(\ref{eq:boundmnujVcirc}) and the assumption on $j$ and $k$.
\end{proof}

\begin{lem}\label{lem:tailboundzetanutjk1}
For any $t\in\R^d$ and any $\eta,j,k\geq 1$ with $2^{-k}\nu_{j}\leq\eta jk$,
\[
\prob(M_{\nu,j}(V^{\circ}_{t,2^{-k}})\geq 5\eta jk)\leq 2\,\ee^{-\eta jk}.
\]
\end{lem}

\begin{proof}
We may clearly assume that $m_{\nu,j}(V^{\circ}_{t,2^{-k}})$ is positive. In view of~(\ref{eq:boundmnujVcirc}), the fact that $M_{\nu,j}(V^{\circ}_{t,2^{-k}})\geq 5\eta jk$ implies that
\[
\int_{(\rho,s)\in V^{\circ}_{t,2^{-k}}\atop|x|\in\Ical_{j}} \Nu^{*}(\dd\rho,\dd s,\dd x)=M_{\nu,j}(V^{\circ}_{t,2^{-k}})-m_{\nu,j}(V^{\circ}_{t,2^{-k}})\geq 3\eta jk,
\]
which may happen with probability at most $2\,\ee^{-\eta jk}$, by virtue of Lemma~\ref{lem:Bernstein}.
\end{proof}

%%%%%%%%%%%%%%%%%%%%%%%%%%%%%%%%%%%%%%%%%%%%%%%%%%%%%%%%%%%%%%%%%%%%%%%%%%%%%%
%%%%%%%%%%%%%%%%%%%%%%%%%%%%%%%%%%%%%%%%%%%%%%%%%%%%%%%%%%%%%%%%%%%%%%%%%%%%%%
\subsubsection{Case where $2^{-k}\nu_{j}>\eta jk$}
%%%%%%%%%%%%%%%%%%%%%%%%%%%%%%%%%%%%%%%%%%%%%%%%%%%%%%%%%%%%%%%%%%%%%%%%%%%%%%
%%%%%%%%%%%%%%%%%%%%%%%%%%%%%%%%%%%%%%%%%%%%%%%%%%%%%%%%%%%%%%%%%%%%%%%%%%%%%%

Let $\delta_{j,k}=(2^{k}\nu_{j})^{-1/2}$ and $\Tau_{j,k}$ be the set of points in $\Q^{d}\cap\clball{0}{2^{1-k}}$ of the form $p/\lfloor (d\,2^{k}\nu_{j})^{1/2}\rfloor$ with $p\in\Z^{d}$. We clearly have $\#\Tau_{j,k}\leq (25d\,2^{-k}\nu_{j})^{d/2}$. Here are the analogs of Lemmas~\ref{lem:boundzetanutjk1} and~\ref{lem:tailboundzetanutjk1}.

\begin{lem}\label{lem:boundzetanutjk2}
For any $t\in\R^d$ and any $\eta,j,k\geq 1$ with $2^{-k}\nu_{j}>\eta jk$,
\begin{align*}
\zeta_{\nu}(t,j,k) \leq&\  2\sup_{\tau,\tau'\in \Tau_{j,k}} |\Lfrak_{\nu,j}(V^{+}_{t+\tau,\delta_{j,k}}\setminus V^{-}_{t+\tau',\delta_{j,k}})|\\
&+2^{3-j}\sup_{\tau\in \Tau_{j,k}} M_{\nu,j}(V^{\circ}_{t+\tau,\delta_{j,k}})+2^{4-j-k/2}\nu_{j}^{1/2}.
\end{align*}
\end{lem}

\begin{proof}
Given $\tau\in\Q^{d}\cap\clball{0}{2^{-k}}$, there clearly exists a $\tau'\in\Tau_{j,k}$ with $\|\tau-\tau'\|\leq\delta_{j,k}$, so that $V^{-}_{t+\tau',\delta_{j,k}}\subseteq V_{t+\tau}\subseteq V^{+}_{t+\tau',\delta_{j,k}}$. Moreover, we also have $V^{-}_{t,\delta_{j,k}}\subseteq V_{t}\subseteq V^{+}_{t,\delta_{j,k}}$. Then, just as in the proof of Lemma~\ref{lem:boundzetanutjk1},
\[
|L_{\nu,j}(t+\tau)-L_{\nu,j}(t)|\leq |\Lfrak_{\nu,j}(V_{t+\tau}\setminus V_{t})| + |\Lfrak_{\nu,j}(V_{t}\setminus V_{t+\tau})|.
\]
Splitting $V^{+}_{t+\tau',\delta_{j,k}}\setminus V^{-}_{t,\delta_{j,k}}$ into its subset $V_{t+\tau}\setminus V_{t}$ and the complement, we get
\[
|\Lfrak_{\nu,j}(V_{t+\tau}\setminus V_{t})| = |\Lfrak_{\nu,j}(V^{+}_{t+\tau',\delta_{j,k}}\setminus V^{-}_{t,\delta_{j,k}})-\Lfrak_{\nu,j}((V^{+}_{t+\tau',\delta_{j,k}}\setminus V^{-}_{t,\delta_{j,k}})\setminus (V_{t+\tau}\setminus V_{t}))|.
\]
Owing to the triangle inequality, the upper bound given by~(\ref{eq:boundLfraknujV}) and the observation that the complement of the set $V_{t+\tau}\setminus V_{t}$ in $V^{+}_{t+\tau',\delta_{j,k}}\setminus V^{-}_{t,\delta_{j,k}}$ is included in the union of $V^{\circ}_{t+\tau',\delta_{j,k}}$ and $V^{\circ}_{t,\delta_{j,k}}$, the right-hand side above is smaller than or equal to
\[
|\Lfrak_{\nu,j}(V^{+}_{t+\tau',\delta_{j,k}}\setminus V^{-}_{t,\delta_{j,k}})|+2^{-j+1}( M_{\nu,j}(V^{\circ}_{t+\tau',\delta_{j,k}}\cup V^{\circ}_{t,\delta_{j,k}}) + m_{\nu,j}(V^{\circ}_{t+\tau',\delta_{j,k}}\cup V^{\circ}_{t,\delta_{j,k}}) ).
\]
In addition,~(\ref{eq:boundmnujVcirc}) ensures that the sets $V^{\circ}_{t+\tau',\delta_{j,k}}$ and $V^{\circ}_{t,\delta_{j,k}}$ have $m_{\nu,j}$-measure at most $2\delta_{j,k}\nu_{j}$. Therefore, $|\Lfrak_{\nu,j}(V_{t+\tau}\setminus V_{t})|$ is at most
\[
|\Lfrak_{\nu,j}(V^{+}_{t+\tau',\delta_{j,k}}\setminus V^{-}_{t,\delta_{j,k}})|+2^{-j+1}( M_{\nu,j}(V^{\circ}_{t+\tau',\delta_{j,k}})+M_{\nu,j}(V^{\circ}_{t,\delta_{j,k}}) + 4\delta_{j,k}\nu_{j} ).
\]
Likewise, $|\Lfrak_{\nu,j}(V_{t}\setminus V_{t+\tau})|$ is smaller than or equal to
\[
|\Lfrak_{\nu,j}(V^{+}_{t,\delta_{j,k}}\setminus V^{-}_{t+\tau',\delta_{j,k}})|+2^{-j+1}( M_{\nu,j}(V^{\circ}_{t+\tau',\delta_{j,k}})+ M_{\nu,j}(V^{\circ}_{t,\delta_{j,k}}) + 4\delta_{j,k}\nu_{j} ),
\]
and the result follows.
\end{proof}

\begin{lem}\label{lem:tailboundzetanutjk2}
For any $t\in\R^d$, any $\eta,j,k\geq 1$ with $2^{-k}\nu_{j}>\eta jk$, and any $\tau,\tau'\in\Tau_{j,k}$,
\[\left\{
\begin{array}{l}
\prob(|\Lfrak_{\nu,j}(V^{+}_{t+\tau,\delta_{j,k}}\setminus V^{-}_{t+\tau',\delta_{j,k}})|\geq 2^{3-j-k/2} (\eta j k\nu_{j})^{1/2})\leq 2\,\ee^{-\eta j k} \\[2mm]
\prob( M_{\nu,j}(V^{\circ}_{t+\tau,\delta_{j,k}})\geq 2^{3-k/2} (\eta j k\nu_{j})^{1/2} )\leq 2\,\ee^{-\eta j k}.
\end{array}
\right.\]
\end{lem}

\begin{proof}
For the first bound, in view of~(\ref{eq:boundLfraknujV}), we may clearly assume that the $m_{\nu,j}$-measure of $V^{+}_{t+\tau,\delta_{j,k}}\setminus V^{-}_{t+\tau',\delta_{j,k}}$ is positive. Moreover, this set is included in $V^{\circ}_{t,3\cdot 2^{-k}}$, so its $m_{\nu,j}$-measure is at most $6\cdot 2^{-k}\nu_{j}$, owing to~(\ref{eq:boundmnujVcirc}). Thus, Lemma~\ref{lem:Bernstein} and the fact that $2^{-k}\nu_{j}>\eta jk$ imply that the probability under study is at most
\[
2\,\exp\left( -\frac{3\cdot 2^{6-2j-k}\eta j k\,\nu_{j}}{2^{5-2j-k/2}(\eta j k\nu_{j})^{1/2}+36\cdot 2^{2-2j-k}\nu_{j}} \right)\leq 2\,\ee^{-\eta j k}.
\]

The second inequality that we need to establish is an upper bound on the probability of an event which implies that
\begin{align*}
\int_{(\rho,s)\in V^{\circ}_{t+\tau,\delta_{j,k}} \atop |x|\in\Ical_{j}}\Nu^{*}(\dd\rho,\dd s,\dd x) &= M_{\nu,j}( V^{\circ}_{t+\tau,\delta_{j,k}} )-m_{\nu,j}( V^{\circ}_{t+\tau,\delta_{j,k}} ) \\
&\geq 2^{3-k/2} (\eta j k\nu_{j})^{1/2} - 2\delta_{j,k}\nu_{j}\geq 4(\eta j k 2^{-k} \nu_{j})^{1/2}.
\end{align*}
Thanks to Lemma~\ref{lem:Bernstein} and~(\ref{eq:boundmnujVcirc}) again, this may happen with probability at most
\[
2\,\exp\left( -\frac{3\eta j k 2^{4-k} \nu_{j}}{2^{3-k/2}(\eta j k\nu_{j})^{1/2}+12(2^{-k}\nu_{j})^{1/2}} \right)\leq 2\,\ee^{-\eta j k},
\]
where the last bound follows from the fact that $2^{-k}\nu_{j}>\eta jk\geq 1$.
\end{proof}

%%%%%%%%%%%%%%%%%%%%%%%%%%%%%%%%%%%%%%%%%%%%%%%%%%%%%%%%%%%%%%%%%%%%%%%%%%%%%%
%%%%%%%%%%%%%%%%%%%%%%%%%%%%%%%%%%%%%%%%%%%%%%%%%%%%%%%%%%%%%%%%%%%%%%%%%%%%%%
\subsection{End of the proof}
%%%%%%%%%%%%%%%%%%%%%%%%%%%%%%%%%%%%%%%%%%%%%%%%%%%%%%%%%%%%%%%%%%%%%%%%%%%%%%
%%%%%%%%%%%%%%%%%%%%%%%%%%%%%%%%%%%%%%%%%%%%%%%%%%%%%%%%%%%%%%%%%%%%%%%%%%%%%%

Let us consider a real $\eta\geq 1$. If $2^{-k}\nu_{j}\leq\eta jk$, let $\event_{j,k}$ denote the event consisting in the fact that the following does {\em not} hold:
\begin{equation}\label{eq:eventjk1}
\sup_{t\in U_{A,k}} M_{\nu,j}(V^{\circ}_{t,2^{-k}})\leq 5\eta jk.
\end{equation}
Thanks to Lemma~\ref{lem:tailboundzetanutjk1}, its probability satisfies
\[
\prob(\event_{j,k})\leq 2\,\ee^{-\eta jk}\,\# U_{A,k}\leq 2(4Ad^{1/2})^{d}\,2^{dk}\,\ee^{-\eta jk}.
\]
Now, if $2^{-k}\nu_{j}>\eta jk$, let $\event_{j,k}$ denote the event consisting in the fact that the following does {\em not} hold:
\begin{equation}\label{eq:eventjk2}
\left\{\begin{array}{l}
\sup\limits_{t\in U_{A,k} \atop \tau,\tau'\in \Tau_{j,k}} |\Lfrak_{j}(V^{+}_{t+\tau,\delta_{j,k}}\setminus V^{-}_{t+\tau',\delta_{j,k}})|\leq 2^{3-j-k/2} (\eta j k\nu_{j})^{1/2} \\[2mm]
\sup\limits_{t\in U_{A,k} \atop \tau\in \Tau_{j,k}} M_{\nu,j}(V^{\circ}_{t+\tau,\delta_{j,k}})\leq 2^{3-k/2} (\eta j k\nu_{j})^{1/2}.
\end{array}\right.
\end{equation}
Owing to Lemma~\ref{lem:tailboundzetanutjk2}, we have
\begin{align*}
\prob(\event_{j,k})&\leq 2\, \ee^{-\eta j k}\,\# U_{A,k} \, \#\Tau_{j,k}(1+\#\Tau_{j,k})\\
&\leq 4(100 A d^{3/2})^{d} \, \nu_{j}^{d} \, \ee^{-\eta j k}\leq 4(100 A d^{3/2}c_{\nu})^{d} \, 2^{2dj} \, \ee^{-\eta j k},
\end{align*}
where $c_{\nu}=\sum_{j\geq 1} 2^{-2j}\nu_{j}$, which is finite owing to~(\ref{eq:nuLevy}).

From now on, let us suppose that $\eta>2d\log 2$. Using the above bounds, it is easy to check that $\sum_{(j,k)\in\N^{2}}\prob(\event_{j,k})<\infty$. Letting $D_{n}$ be the set of $(j,k)\in\N^{2}$ with $\max\{j,k\}\geq n$, we deduce that $\prob(\bigcap_{n=1}^{\infty}\downarrow\bigcup_{(j,k)\in D_{n}}\event_{j,k})=0$. So, with probability one, there is an integer $n\geq 1$ such that~(\ref{eq:eventjk1}) holds for any $(j,k)\in D_{n}$ with $2^{-k}\nu_{j}\leq\eta jk$, and~(\ref{eq:eventjk2}) holds for any $(j,k)\in D_{n}$ with $2^{-k}\nu_{j}>\eta jk$. In the first case, it follows from Lemma~\ref{lem:boundzetanutjk1} and~(\ref{eq:boundzetanuAjk}) that $\zeta_{\nu}(A,j,k)$ is bounded by $14(3d+2)\eta\, 2^{-j} j k$. In the second case, it is bounded by $96(3d+2)2^{-j-k/2} (\eta j k \nu_{j})^{1/2}$, owing to Lemma~\ref{lem:boundzetanutjk2} and~(\ref{eq:boundzetanuAjk}) again. Letting $\eta=2d>2d\log 2$, we finally get
\[
\as \quad \exists n\geq 1 \quad \sup_{(j,k)\in D_{n}} \frac{\zeta_{\nu}(A,j,k)}{2^{-j}k(j+2^{-k/2}(j\nu_{j})^{1/2})}\leq 192(3d+2)d.
\]
In addition, for any $j,k\geq 1$, we deduce from~(\ref{eq:boundLfraknujV}) that the expectation of $\zeta_{\nu}(A,j,k)$ is at most $2^{4-j}(A+2^{-k})\nu_{j}$, so that $\zeta_{\nu}(A,j,k)<\infty$ almost surely. Thus,
\[
\as \quad \forall n\geq 2 \quad \sup_{(j,k)\in\N^{2}\setminus D_{n}} \frac{\zeta_{\nu}(A,j,k)}{2^{-j}k(j+2^{-k/2}(j\nu_{j})^{1/2})}<\infty.
\]
Lemma~\ref{lem:oscLnuj} now clearly follows.

%%%%%%%%%%%%%%%%%%%%%%%%%%%%%%%%%%%%%%%%%%%%%%%%%%%%%%%%%%%%%%%%%%%%%%%%%%%%%%
%%%%%%%%%%%%%%%%%%%%%%%%%%%%%%%%%%%%%%%%%%%%%%%%%%%%%%%%%%%%%%%%%%%%%%%%%%%%%%
%%%%%%%%%%%%%%%%%%%%%%%%%%%%%%%%%%%%%%%%%%%%%%%%%%%%%%%%%%%%%%%%%%%%%%%%%%%%%%
% References
%%%%%%%%%%%%%%%%%%%%%%%%%%%%%%%%%%%%%%%%%%%%%%%%%%%%%%%%%%%%%%%%%%%%%%%%%%%%%%
%%%%%%%%%%%%%%%%%%%%%%%%%%%%%%%%%%%%%%%%%%%%%%%%%%%%%%%%%%%%%%%%%%%%%%%%%%%%%%
%%%%%%%%%%%%%%%%%%%%%%%%%%%%%%%%%%%%%%%%%%%%%%%%%%%%%%%%%%%%%%%%%%%%%%%%%%%%%%

\end{document}